\documentclass[11pt]{article}\input{amssym.def}\input{amssym}
\usepackage[active]{srcltx}

\newenvironment{annotacia}{\centerline{\sc Abstract}\vspace{2mm}\narrower\narrower\sf}

\def\theequation{\thesection.\arabic{equation}}
\makeatletter\@addtoreset{equation}{section}\makeatother
\def\nn{\nonumber}\def\lb{\label}\def\nin{\noindent}\def\Mt{M \raisebox{1mm}{$\intercal$}}
\def\be{\begin{equation}}\def\ee{\end{equation}}\def\ba{\begin{eqnarray}}\def\ea{\end{eqnarray}}
\def\tr{{\rm Tr}\,}\def\Tr#1{{\rm Tr}_{\! R^{\mbox{\scriptsize$(#1)$}}}}
\def\TR#1#2{{\rm Tr}_{\! #2^{\mbox{\,\scriptsize$(#1)$}}}}\def\str#1{\rule[#1mm]{0pt}{1mm}}

\def\cR{{\cal R}}\def\cF{{\cal F}}\def\ot{\otimes}

\def\a{\alpha}\def\b{\beta}\def\g{\gamma}\def\d{\delta}\def\D{\Delta}\def\e{\epsilon}
\def\id{{\mathrm{id}}}\def\vf{v_{_\cF}}\def\ti{\tilde}

\newcounter{theorem}\makeatletter
\@addtoreset{theorem}{section}\makeatother

\newtheorem{prop}[theorem]{Proposition}
\newtheorem{rem}[theorem]{Remark}
\newtheorem{lem}[theorem]{Lemma}
\newtheorem{def-lem}[theorem]{Definition-Lemma}
\newtheorem{def-prop}[theorem]{Definition-Proposition}
\newtheorem{defin}[theorem]{Definition}
\newtheorem{theor}[theorem]{Theorem}
\newtheorem{cor}[theorem]{Corollary}

\begin{document}

\title{Quantum Matrix Algebras of BMW type:\\
Structure of the Characteristic Subalgebra}
\author{ \rule{0pt}{7mm} Oleg Ogievetsky\thanks{oleg@cpt.univ-mrs.fr~~ ORCID iD 0000-0003-2444-2789}\\
\\ {\small\it
Aix Marseille Universit\'e, Universit\'e de Toulon, CNRS,}\\
{\small\it CPT UMR 7332, 13288, Marseille, France}
\\[-3pt]
{\small \&}\\[-3pt]
{\small\it I.E.Tamm Department of Theoretical Physics, P.N. Lebedev Physical Institute,}\\[-2pt]
{\small\it Leninsky prospekt 53, 119991 Moscow, Russia}\\ \\ \\
\rule{0pt}{7mm} Pavel Pyatov\thanks{pyatov@theor.jinr.ru~~ ORCID iD 0000-0002-9773-3600}\\ \\
{\small\it National Research University Higher School of Economics}\\[-2pt]
{\small\it 20 Myasnitskaya street, Moscow 101000, Russia}\\[-2pt]
{\small \&}\\[-5pt]
{\small\it Bogoliubov Laboratory of Theoretical Physics}\\[-2pt]
{\small\it JINR, 141980 Dubna, Moscow region, Russia}}
\date{}
\maketitle

\begin{annotacia}
A notion of quantum matrix (QM-) algebra generalizes and unifies two
famous families of algebras from the theory of quantum groups: the RTT-algebras and the reflection equation (RE-) algebras. These algebras being generated by the components of a `quantum' matrix $M$ possess certain  properties which resemble structure theorems of the ordinary matrix theory.  It turns out that
such structure results are naturally derived in a more general framework of the QM-algebras.  In this work we consider a family of Birman-Murakami-Wenzl  (BMW) type QM-algebras. These algebras are defined with the use of R-matrix representations of the BMW algebras. Particular series of such algebras include orthogonal and symplectic types RTT- and RE- algebras, as well as their super-partners.

For a  family of BMW type QM-algebras, we
investigate the structure of their `characteristic subalgebras' ---
the subalgebras where the coefficients of characteristic polynomials take values. We define three sets
of generating elements of the characteristic subalgebra and derive recursive Newton and Wronski relations between them.
We also define  an associative
$\star$-product for the matrix $M$ of generators of the QM-algebra which is a proper generalization of the classical matrix multiplication. We
determine the set of all matrix `descendants' of the quantum matrix $M$,  and
prove the $\star$-commutativity of this set in the BMW type.
\end{annotacia}

\newpage

\tableofcontents

%\bigskip\bigskip\bigskip

%\newpage
\section{Introduction}\lb{sec1}

A notion of a quantum matrix group, also called the RTT-algebra, is implicit in
%structures of
the quantum inverse scattering method. A formal definition has been given in the works of
V. Drinfel'd, L. Faddeev, N. Reshetikhin and L. Takhtajan \cite{D1,FRT}. Since then, various aspects
of the quantum matrix group theory have been elaborated, especially in attempts to
define differential geometric structures on non-commutative spaces (see, e.g., \cite{Man,SchWZ}%\footnote{Our bibliography on the subject is not exhaustive here and below.
%%We give only those references, which are close in ideas and technique to our present work
%%and may be considered as our sources of motivation.
%}
). In particular, a different family of algebras
generated by matrix components, the so-called reflection equation (RE-) algebras \cite{C,KS},
has been brought into consideration. Soon it was realized that, for both the RTT- and the RE-algebras,
some of the basic concepts of the classical matrix algebra, like the notion of the spectral invariants
and the characteristic identity (the Cayley-Hamilton theorem) can be properly generalized (see
\cite{EOW,NT,PS,Zh}). So, it comes out that the matrix notation used for the definition of
the RTT- and the RE-algebras is not only technically convenient, but it dictates certain structure
properties for the algebras themselves. It is then natural to search for a possibly most general
algebraic setting for the matrix-type objects. Such family of algebras was introduced in refs.\cite{Hl}
and \cite{IOP1}, and in the latter case the definition was dictated by a condition that
the standard matrix theory statements should have their appropriate generalizations. These algebras were
called quantum matrix (QM-) algebras although one should have in mind that the QM-algebras
are generated by the matrix components rather than by the matrix itself.

\smallskip
The RTT- and the RE-algebras are probably the most important subfamilies in the variety of QM-algebras.
They are distinguished both from the algebraic point of view (the presence of additional
non-braided bi-algebra and bi-comodule structures) and from the geometric point of view (their interpretation
as, respectively, the algebras of quantized functions and of quantized invariant differential
operators on a group); also, the RE-algebras naturally appear in the representation theory, in the description of the diagonal reduction algebras \cite{KhO}.   However, for the generalization of the basic matrix algebra statements,
it is not only possible but often more clarifying to use a weaker structure settings of the QM-algebras.

\smallskip
So far, the program of generalizing the Cayley-Hamilton theorem was fully accomplished for the 'linear' (or Iwahori-Hecke) type QM-algebras. For the $GL(m)$-type algebras, the results were described in \cite{GPS1,IOPS,IOP1}
and for the $GL(m|n)$-type algebras in \cite{GPS2,GPS3}. These works generalize earlier results
on characteristic identities by A.J. Bracken,  H.S. Green, et. al.,
in the Lie (super)algebra case \cite{BGr,Gr,BCC,Gou,JGr} (for a review see \cite{IWG}) and in the quantized universal enveloping algebra case
\cite{GZB}, and by I. Kantor and I. Trishin in the matrix superalgebra case \cite{KT1,KT2}.

\smallskip
The similar investigation program for the QM-algebras of  Birman-Murakami-Wenzl (BMW) type  (for their definition see section \ref{subsec4.1}) was initiated in \cite{OP}.
In the present and forthconimng works we continue and complement this program.
The family of BMW type QM-algebras serves as a unifying set-up for the description of the orthogonal and symplectic
QM-algebras as well as for their supersymmetric partners.
Some partial results about specific examples of such algebras and their limiting cases were already derived.
In particular, the characteristic identities for the generators  of the orthogonal and symplectic Lie algebras have been considered at the representation theoretical and at the abstract algebraic levels in \cite{BGr,Gr} and in \cite{BCC,Gou,Mol}. The characteristic identities for the canonical Drinfeld-Jimbo quantizations of the orthogonal and symplectic universal enveloping algebras  were obtained  in \cite{MRS} and their images in the series of highest weight representations were discussed in details in \cite{Mudr}.
So, it is pretty clear that proper generalizations of the Cayley-Hamilton theorems do exist for the families of orthogonal and symplectic QM-algebras. However, in a derivation of these results one meets serious technical complications. The reason is that the structure of the
Birman-Murakami-Wenzl algebras
%which underly the construction of the BMW type QM-algebras
is substantially more sophisticated then that of the Iwahori-Hecke algebras (Iwahori-Hecke and Birman-Murakami-Wenzl algebras play similar roles in the construction of the  QM-algebras of linear and BMW types).
In the present work we develop an appropriate techniques to deal with these  complications.
%Using this techniques we study the so-called {\em characteristic} subalgebra of the BMW type QM-algebra.
%In the Iwahori-Hecke case, it is
%the subalgebra where the coefficients of the Cayley-Hamilton identities take their values.
\smallskip

In sections \ref{sec2} and \ref{sec3} we collect necessary results concerning the Birman-Murakami-Wenzl (BMW) algebras and their R-matrix representations.
In the beginning of section \ref{sec2} we define the BMW
algebras in terms of generators and relations, describe few helpful morphisms between these algebras, and introduce the baxterized elements.
These elements are used in  subsection \ref{subsec2.2} for the definition of three sets of idempotents called antisymmetrizers, symmetrizers and contractors. Necessary properties of these idempotents are
proved in  proposition \ref{proposition2.2}. All the material of this section, except the construction and  properties of the contractors
is fairly well known and we present it to make the presentation self-contained.

\smallskip
In  section \ref{sec3} we consider the R-matrix representations of the BMW algebras. We define
standard notions of the R-trace~\footnote{This operation is also called a quantum trace or, shortly,
a $q$-trace in the literature.},  skew-invertibility,   compatible pair of R-matrices
and R-matrix twist  (subsection \ref{subsec3.1}). In
subsection \ref{subsec3.2} we collect necessary formulas and statements relating the notions introduced before. To investigate the skew-invertibility of the
R-matrix after a twist, in  subsection \ref{subsec3.2a} we derive an expression for the twisted R-matrix, which is different from the standard one.
Next we describe the BMW type R-matrices (subsection \ref{subsec3.3}).
The major part of a technical preparatory work is done in
subsections \ref{subsec3.2}---\ref{subsec3.3}, and \ref{operatorG}, \ref{Twolinearmaps}. Here we develop the R-matrix
technique, which is later used in the main sections \ref{sec4}, \ref{sec5}.
\smallskip

In the beginning of  section \ref{sec4} we introduce the  QM-algebras of general and BMW types. We then define the characteristic subalgebra of the QM-algebra. In the Iwahori-Hecke case, it is
the subalgebra where the coefficients of the Cayley-Hamilton identity take their values. As it was shown
in \cite{IOP1}, the characteristic subalgebra is abelian. In subsection \ref{subsec4.2} we
describe three generating sets for the characteristic subalgebra of the BMW type QM-algebra.
As compared to the linear QM-algebras, all these generating sets contain  a single additional element --- the 2-contraction $g$ --- which at the classical level gives rise to bilinear invariant 2-forms for the orthogonal and symplectic groups.

\smallskip
Next, in  subsection \ref{subsec4.4}, we construct
%another important requisite of the Cayley-Hamilton theorem ---
a proper analogue of the matrix multiplication for the quantum matrices. We call it
the quantum matrix product `$\star $'.
%The definition is given for the general QM-algebra.
In general, the $\star \, $-product is different from the usual matrix product. It is worth noting that for the family
of RE-algebras, the $\star \, $-product coincides with the matrix product. The $\star \, $-product is proven to be associative and hence the $\star \,$-powers
of the same quantum matrix $M$ commute. We determine then the set of all  `descendants' of  the
quantum matrix $M$ in the BMW case and prove that this set is $\star \, $-commutative. It turns out that,
unlike the linear QM-algebra case, it is not possible to express all these descendants in terms of the
$\star \, $-powers of $M$ only. The expressions include also a new operation `\raisebox{1mm}{$\intercal$}',
which can be treated as a `matrix multiplication with a transposition'.
%This fact causes some technical
%complications in the proof of the Cayley-Hamilton theorem,
%but on the other hand, it gives rise to a
%diversity of the characteristic identities in the
%orthogonal and symplectic cases.

\smallskip
In subsection \ref{subsec4.5} we define an extension of the BMW type QM-algebra by the element $g^{-1}$ which is the inverse to the 2-contraction.
Then we construct in the extended algebra the inverse $\star \, $-power of the quantum matrix $M$.

\smallskip
The last section \ref{sec5} contains  the principal result of the present work, theorem  \ref{theorem6.1}, which establishes, for the BMW type QM-algebras, recursive relations between the elements of the three generating sets of their characteristic subalgebras.
% which we introduced in subsection \ref{subsec4.2}.
These formulas generalize the classical Newton and Wronsky relations for the sets of the power sums, elementary and complete symmetric polynomials (see \cite{Mac}) to the case of quantum matrices and simultaneously, to the situation where additional element of the characteristic subalgebra,
%independent symmetric polynomial,
the 2-contraction, is present.
To prove this result we first derive the matrix relations among the descendants of the  BMW type quantum matrix $M$ (see lemma \ref{lemma5.1}). These relations can be viewed as the matrix counterparts of the Newton relations, and they are expected to be important ingredients  in a future derivation of the characteristic identities for the  QM-algebras of the BMW type.

\smallskip
Some auxilliary results, which are interesting in themselves, although not necessary for considerations in the main text, are collected in the appendices.  In appendix \ref{primcontr}  we prove the primitivity of the contractors from subsection \ref{subsec2.4}. In appendix \ref{fuproco}  their further  properties are discussed. Appendix \ref{append1} is devoted to a discussion of universal counterparts of the matrix relations given in subsections \ref{subsec3.2},
\ref{subsec3.2a}.
\smallskip

In  forthcoming papers we are going to construct
the Cayley-Hamilton identities, and, more generally Cayley-Hamilton-Newton identities in the spirit of \cite{IOP}, for the series of orthogonal and symplectic QM-algebras and, further on, for their super-partners.
\smallskip

%\newpage
\section{Some facts about Birman-Murakami-Wenzl algebras}
\lb{sec2}

In this preparatory section we collect definitions and derive few results on the Birman-Murakami-Wenzl
algebras. We give a minimal information, which is required for the main part of the paper.
In particular, in sec.\ref{subsec2.2} we describe series of morphisms of the braid groups and their quotient BMW algebras; in sec.\ref{subsec2.3} we introduce baxterized elements which are then used
in the sec.\ref{subsec2.4} to define three series of idempotents in the BMW algerbas, the so called symmetrizers, antisymmetrizers and contractors.

The reader will find a more detailed presentation of the Birman-Murakami-Wenzl algebras in, e.g.,
papers \cite{W} and \cite{LR}.

\subsection{Definition}
\lb{subsec2.1}

The braid group ${\cal B}_{n}$, $n\geq 2$, in Artin presentation, is defined by generators
$\{\sigma_i\}_{i=1}^{n-1}$ and relations
\ba
\lb{braid}
\sigma_i \sigma_{i+1} \sigma_i & =&  \sigma_{i+1} \sigma_i \sigma_{i+1}\quad\, \forall\; i=1,2,\dots ,n-1,
\\
\lb{braid2}
\sigma_i \sigma_j & =&  \sigma_j \sigma_i\qquad  \qquad \forall\; i,j:\; |i-j|>1\ .
\ea
We put, by definition, ${\cal B}_{1}:=\{1\}$.

The {\em Birman-Murakami-Wenzl (BMW) algebra} ${\cal W}_{n}(q,\mu)$ \cite{BW,M1} is a finite
dimensional quotient algebra of the group algebra ${\Bbb C}{\cal B}_{n}$. It depends on
two complex parameters $q$ and $\mu$. Let
\be
\lb{kappa}
\kappa_{i} := {(q1-\sigma_i)(q^{-1} 1 +\sigma_i)\over \mu (q-q^{-1})} ,
\qquad i=1,2,\dots ,n-1 \ .
\ee
The quotient algebra ${\cal W}_{n}(q,\mu)$ is specified by conditions
\ba
\lb{bmw2a}
\sigma_i \kappa_i \, =\,  \kappa_i \sigma_i &=& \mu \kappa_i\, ,
\\
\lb{bmw2b}
\kappa_i \sigma_{i+1}^{\epsilon} \kappa_i &=&  \mu^{-\epsilon} \kappa_i\ ,
\ea
where $\epsilon$ is the sign~\footnote{If $\mu\neq q-q^{-1}$ then it is
enough to impose only one of the relations (\ref{bmw2b}), the relation with another sign follows (see \cite{IOP3}).}, $\epsilon =\pm 1$.

\medskip
Eqs.(\ref{kappa}) and (\ref{bmw2a}) imply that the characteristic polynomial for the generator $\sigma_i$ has degree three,
\be
(\sigma_i -q1)(\sigma_i +q^{-1}1)(\sigma_i -\mu 1)=0\ .
\lb{deg3-a}
\ee

The relations (\ref{bmw2a}) -- (\ref{bmw2b}) imply also
\ba
\lb{bmw7}
\sigma'_i\kappa_{i+1}\sigma'_i &=& \sigma'_{i+1}
\kappa_i\sigma'_{i+1}\, ,
\quad\;
\mbox{where}\quad\; \sigma' =\sigma-\raisebox{1pt}{$(q-q^{-1})$}1 \ ,
\\[5pt]
\lb{bmw3}
\kappa_i \sigma_{i+\pi}^{\epsilon}  &=&\kappa_i \kappa_{i+\pi}\sigma_i^{-\epsilon}\,,
\quad\;\;
%\phantom{\mbox{where}\quad\;}
\sigma_{i+\pi}^{\epsilon} \kappa_{i}\,  =\,
\sigma_{i}^{-\epsilon} \kappa_{i+\pi} \kappa_{i}\, ,
\\[5pt]
\lb{bmw5a}
\kappa_i \kappa_{i+\pi} \kappa_i &=& \kappa_i\, ,
\\
\lb{bmw5b}
\kappa_i^2 &=& \eta\, \kappa_i\, ,
\qquad\qquad\;
\mbox{where}\quad\;
\eta:= {\displaystyle {(q-\mu)(q^{-1}+\mu)\over\mu (q-q^{-1})}}\, .
\ea
Here $\epsilon$ and $\pi$ are the signs:
$\epsilon =\pm 1$ and $\pi =\pm 1$.

The parameters $q$ and $\mu$ of the BMW algebra are taken in domains\footnote{For particular values  $\mu=\pm q^{i}$, $i\in {\Bbb Z}$, the limiting cases $q\rightarrow \pm 1$ to the Brauer algebra
\cite{Br} can be consistently defined.}
\be
\lb{q-mu-restrict}
q\in{\Bbb C}\backslash \{0,\pm 1\}, \qquad
\mu\in{\Bbb C}\backslash\{0,q,-q^{-1}\},
\ee
so that the elements $\kappa_i$ are well defined and
non-nilpotent. Further restrictions on $q$ and $\mu$ will be imposed in subsection \ref{subsec2.3}.

\subsection{Natural morphisms}
\lb{subsec2.2}

\medskip
{\bf $\bullet$}~ The braid groups and their quotient BMW algebras admit a chain of monomophisms
\be
\begin{array}{l}
{\cal B}_2\hookrightarrow\dots\hookrightarrow {\cal B}_n\hookrightarrow
{\cal B}_{n+1}\hookrightarrow\dots \ \ ,\\[1em]{\cal W}_2\hookrightarrow\dots\hookrightarrow
{\cal W}_n\hookrightarrow {\cal W}_{n+1}\hookrightarrow\dots\end{array}
\lb{h-emb}
\ee
defined on the generators as
\be\lb{h-emb2}{\cal B}_{n}\ ({\mathrm or}\ \ {\cal W}_{n})\ni \sigma_i \mapsto
\sigma_{i+1}\in {\cal B}_{n+1}\ ({\mathrm or}\ \ {\cal W}_{n+1})\;\; \forall\; i=1,\dots ,n-1.\ee
We denote by $\alpha^{(n)\uparrow i}\in {\cal B}_{n+i}\ ({\mathrm{or}}\ \ {\cal W}_{n+i})$
an image of an element $\alpha^{(n)}\in {\cal B}_n\ ({\mathrm{or}}\ \ {\cal W}_{n})$
under a composition of the mappings (\ref{h-emb})--(\ref{h-emb2}). Conversely, if for some
$j<(n-1)$, an element $\alpha^{(n)}$ belongs to the image of ${\cal B}_{n-j}\ ({\mathrm{or}}\
\ {\cal W}_{n-j})$ in ${\cal B}_n\ ({\mathrm{or}}\ \ {\cal W}_{n})$ then by $\alpha^{(n)\downarrow j}$
we denote the preimage of $\alpha^{(n)}$ in ${\cal B}_{n-j}\ ({\mathrm{or}}\ \ {\cal W}_{n-j})$.

\medskip
This notation will be helpful in subsection \ref{subsec2.4} where we discuss
three distinguished sequences of idempotents in the BMW algebras.

\medskip
\noindent
{\bf $\bullet$}~ Consider series of elements $\tau^{(n)}\in{\cal B}_n$ defined inductively
\be
\lb{tau-n}
\tau^{(1)}:=1, \;\;\; \tau^{(j+1)}:=\tau^{(j)}\,\sigma_j\sigma_{j-1}\dots\sigma_1\, .
\ee
$\tau^{(n)}$ is the lift of the longest element of the symmetric group $S_n$.
% It follows that $(\tau^{(n)})^2$ is central in ${\cal B}_{n}$.
The inner  ${\cal B}_n$ (and, hence, ${\cal W}_n$) automorphism
\be
\lb{innalis}
\tau :\, \sigma_i\mapsto \tau^{(n)}\, \sigma_i\, (\tau^{(n)})^{-1}=\sigma_{n-i}\ ,
\ee
will be used below in derivations in sections \ref{subsec2.4} and \ref{sec4}.

\medskip
\noindent
{\bf $\bullet$}~ One has three algebra isomorphisms:
$$
\iota :\ {\cal W}_n (q,\mu) \rightarrow {\cal W}_n (-q^{-1},\mu)\; \mbox{,}\quad
\iota' :\ {\cal W}_n (q,\mu) \rightarrow {\cal W}_n (q^{-1},\mu^{-1})\quad \mbox{and}\quad
\iota'' :\ {\cal W}_n (q,\mu) \rightarrow {\cal W}_n (-q,-\mu)
$$
defined on generators by
\ba
\lb{homS-A}
\iota : && \sigma_i\mapsto\;\;\;\,\sigma_i\  ,
\\
\lb{homiota'}
\iota' : && \sigma_i\mapsto\;\;\;\,\sigma_i^{-1}\  ,
\\
\lb{homiota''}
\iota'' : && \sigma_i\mapsto -\,\sigma_i\  .
\ea
The map $\iota$ interchanges the two sets of baxterized elements $\sigma^\pm(x)$
and the series of symmetrizers  $a^{(n)}$ and antysimmetrizers $s^{(n)}$: $\iota(a^{(n)})=s^{(n)}$ (see subsections \ref{subsec2.3} and \ref{subsec2.4} below).
For the maps $\iota'$, $\iota''$ one has: $\iota'(\sigma^\pm(x))=x\sigma^\pm(x^{-1})$,~ $\iota''(\sigma^\pm(x))=\sigma^\pm(x)$. The series of (anti)symmetrizers  are stable under maps $\iota'$ and $\iota''$.
One also has $\iota(\kappa_i)=\iota' (\kappa_i)=\iota'' (\kappa_i)=\kappa_i$.

\medskip
\noindent
{\bf $\bullet$.}~ There exists an algebra antiautomorphism $\varsigma :\ {\cal W}_n (q,\mu) \rightarrow {\cal W}_n (q,\mu)$
($\varsigma (xy)=\varsigma (y)\varsigma (x)$), defined on generators as
\be
\varsigma :\
\sigma_i\mapsto\sigma_i\ .
\lb{antvs}
\ee
This morphism will be used later in the proofs of Propositions \ref{proposition2.2} and \ref{proposition4.14}.

\subsection{Baxterized elements}
\lb{subsec2.3}

A set of elements $\sigma_i(x)$, $i=1,2,\dots, n-1,$ depending on a complex parameter $x$, in a
quotient of the group algebra ${\Bbb C}{\cal B}_{n}$ is called a set of {\em baxterized elements} if
\ba\lb{bYBE}\sigma_i(x)\, \sigma_{i+1}(xy)\, \sigma_i(y)\, =\,
\sigma_{i+1}(y)\, \sigma_i(xy)\, \sigma_{i+1}(x)\, \ea
for $i=1,2,\dots, n-1$ and
\ba\lb{bcomm}
\sigma_i(x)\, \sigma_j(y)\, =\,\sigma_j(y)\, \sigma_i(x)\,  \ea
if $|i-j|>1$.

\begin{lem}
\lb{lemma2.1}
{\rm\bf \cite{J,I}} For the algebra ${\cal W}_{n}(q,\mu)$, the baxterized elements exist. There are two
sets of the baxterized elements $\{\sigma^{\varepsilon}_i\}$, $\varepsilon =\pm 1$, given by
\be
\lb{ansatz}
\sigma_i^{\varepsilon}(x)\, :=\,  1\, +\, {x-1\over q-q^{-1}}\, \sigma_i\, +\,
{x-1\over \alpha_{\varepsilon} x+1}\, \kappa_i\, ,
\ee
where $\alpha_{\varepsilon}\, :=\, -\varepsilon q^{-\varepsilon} \mu^{-1}$.
\end{lem}

The complex argument $x$, traditionally called {\em the spectral parameter},
is chosen in a domain $\ $ ${\Bbb C}\setminus\{-\alpha_{\varepsilon}^{-1}\}$.

\subsection{Symmetrizers, antisymmetrizers and contractors}
\lb{subsec2.4}

In terms of the baxterized generators we construct two series of elements $a^{(i)}$ and $s^{(i)}$,
$i=1,2,\dots ,n$, in the algebra ${\cal W}_n(q,\mu)$. They are defined iteratively in two ways:
\ba
\lb{ind1}
a^{(1)} := 1\ \ \ {\mathrm{and}}\ \ \ s^{(1)}\, :=\, 1\, ,\ \ \ \ \ \ \ \ \ \ \, &&
\\[1em]
\lb{a^k}
a^{(i+1)} :={q^i\over (i+1)_q} a^{(i)}\, \sigma^{-}_i(q^{-2i})\,
a^{(i)}\ \ &{\mathrm{or}}&\ \  a^{(i+1)}:=\, {q^i\over (i+1)_q} a^{(i)\uparrow 1}\, \sigma^{-}_1(q^{-2i})\,
a^{(i)\uparrow 1}\, ,
\\[1em]
\lb{s^k}
s^{(i+1)} :={q^{-i}\over (i+1)_q}\, s^{(i)}\, \sigma^{+}_i(q^{2i})\,
s^{(i)}\ \ \ \  &{\mathrm{or}}&\ \  s^{(i+1)}:=\, {q^{-i}\over (i+1)_q}\, s^{(i)\uparrow 1}\,
\sigma^{+}_1(q^{2i})\, s^{(i)\uparrow 1}\, ,
\ea
where $i_q$ are usual $q$-numbers, $i_q\, :=\, (q^i - q^{-i})/(q-q^{-1})$. 
Below we show that  in each of eqs. (\ref{a^k}), (\ref{s^k}) the two definitions coincide. 
We note that the factorized formula for the (anti)symmetrizers, in the spirit of the fusion procedure for the BMW algebra \cite{IMO}, follows 
from the eqs. (\ref{a^k}), (\ref{s^k}).

\smallskip
To avoid singularities in the definition of $a^{(i)}$ (respectively, $s^{(i)}$), $i=1,2,\dots ,n$,
we impose further restrictions on the parameters of ${\cal W}_n(q,\mu)$:
\be
\lb{mu}
j_q\,\neq\, 0\, , \quad \mu \neq -q^{-2j+3} \;\; (\mbox{respectively,~}\
\mu \neq q^{2j-3})\, \quad \forall\;j = 2,3,\dots , n\ .
\ee

\smallskip
The elements $a^{(i)}$ and $s^{(i)}$ are called an {\em $i$-th order antisymmetrizer} and
an {\em $i$-th order symmetrizer}, respectively.

\smallskip
The second order antisymmetrizer and symmetrizer
\be
\lb{as-2}
\hspace{-1mm}
a^{(2)}={q\over 2_q}\sigma_1^{-}(q^{-2})={(q1-\sigma_1)(\mu 1-\sigma_1)\over 2_q (\mu+q^{-1})} , \quad
s^{(2)}={q^{-1}\over 2_q}\sigma_1^{+}(q^2)={(q^{-1}1+\sigma_1)(\mu 1-\sigma_1)\over 2_q (\mu-q)}\,
\ee
are the idempotents participating in a resolution of unity in the algebra ${\cal W}_2(q,\mu)$
(c.f. with the property (\ref{deg3-a})$\,$),
\be
\lb{resolution}
1 \, =\, a^{(2)} + s^{(2)} + \eta^{-1}\kappa_1\, .
\ee
%the spectral decomposition of the generator $\sigma_1$ of ${\cal W}_2(q,\mu)$ is
%\be
%\lb{specdec}
%\sigma_1 \, =\, -q^{-1} a^{(2)} +q s^{(2)} +\mu \eta^{-1}\kappa_1\, ;
%\ee

\smallskip
Likewise for $a^{(2)}$ and $s^{(2)}$, one can introduce higher order analogues for the third idempotent
entering the resolution. Namely, define iteratively
\be
\lb{kappa-i}
c^{(2)}\, :=\, \eta^{-1}\kappa_1\, , \qquad c^{(2i+2)}\, :=\,
c^{(2i)\uparrow 1}\, \kappa_1 \kappa_{2i+1}\, c^{(2i)\uparrow 1}\, .
\ee
The element $c^{(2i)}$ is called an {\em $(2i)$-th order contractor}. Main properties of
the (anti)symmetrizers and contractors are summarized below.

\begin{prop}
\lb{proposition2.2}
Two expressions given for the antisymmetrizers and symmetrizers in eqs.(\ref{a^k}) and (\ref{s^k})
are identical. The elements $a^{(n)}$ and $s^{(n)}$ are central primitive idempotents
in the algebra ${\cal W}_n(q,\mu)$. One has
\ba
\lb{idemp-1}
a^{(n)}\sigma_i  =\sigma_i a^{(n)} = -q^{-1} a^{(n)},&&
s^{(n)}\sigma_i = \sigma_i s^{(n)} = q s^{(n)}\quad\quad\;\;\forall\; i=1,2,\dots ,n-1\,
\\[-3pt]
\nn
%\ea
\mbox{and}\hspace{57mm} &&
\\[-3pt]
%\ba
\lb{idemp-2}
a^{(n)} a^{(m)\uparrow i} = a^{(m)\uparrow i} a^{(n)} = a^{(n)} ,&&
s^{(n)} s^{(m)\uparrow i} = s^{(m)\uparrow i} s^{(n)} = s^{(n)} \qquad\, \mbox{if~~}\; m+i\leq n\, .
\hspace{12mm}
\ea
The antisymmetrizers $a^{(n)}$, for all $n=2,3,\dots$, are orthogonal to the symmetrizers
$s^{(m)}$, for all $m=2,3,\dots$ ,
\be
a^{(n)}s^{(m)}=0\ .
\lb{asort}
\ee

\smallskip
The element $c^{(2n)}$ is a primitive  idempotent in the algebra ${\cal W}_{2n}(q,\mu)$
and in the algebra ${\cal W}_{2n+1}(q,\mu)$. One has
\ba
\lb{idemp-c1}
c^{(2n)} c^{(2i)\uparrow n-i}\,=\, c^{(2i)\uparrow n-i} c^{(2n)}\, =\, c^{(2n)}\, &&
\forall\; i=1,2,\dots , n\, ;
\\[1em]
\lb{idemp-c2}
c^{(2n)} \sigma_i = c^{(2n)} \sigma_{2n-i}\, , \qquad
\sigma_i c^{(2n)} = \sigma_{2n-i}\, c^{(2n)} \, &&\forall\; i=1,2,\dots ,n-1\, ,
%\ea
\\[-3pt]
\nn
\mbox{and}\hspace{87mm} &&
\\[-3pt]
%\ba
\lb{idemp-c3}
c^{(2n)}\sigma_n = \sigma_n\, c^{(2n)} =\mu c^{(2n)}\, .
\hspace{31mm} &&
\ea
The contractors $c^{(2n)}$ are orthogonal to the antisymmetrizers $a^{(m)}$ and to the symmetrizers
$s^{(m)}$ for all $m>n$.
\end{prop}

\nin{\bf Proof.~} The explicit formula (\ref{a^k}) for idempotents, which we call  antisymmetrizers here, appears
in \cite{TW}, although without referring to the baxterized elements (see the proof of the lemma
7.6 in \cite{TW}).\footnote{Different expressions for the antisymmetrizers and symmetrizers, which are less suitable for
our applications, were derived in \cite{HSch}.} Our proof of the formulas (\ref{idemp-1}) and (\ref{idemp-2}) relies
on the relations (\ref{bYBE}) for the baxterized generators.

\medskip
We first check that the elements $a^{(i)}$ defined iteratively by the first formula in
(\ref{a^k}) satisfy the relations (\ref{idemp-1}) and (\ref{idemp-2}). The equalities
(\ref{idemp-1}) for the antisymmetrizers are equivalent to
$$
a^{(n)}s^{(2)\uparrow i-1}=s^{(2)\uparrow i-1}a^{(n)}=
a^{(n)}c^{(2)\uparrow i-1}=c^{(2)\uparrow i-1}a^{(n)}=0\ ,
\quad \forall\; i=1,2,\dots n-1\ ,
$$
which, in turn, are equivalent to
\be
\lb{divide}a^{(n)}
\,\sigma_i^{-}(q^2) = \sigma_i^{-}(q^2)\, a^{(n)} = 0\, .
\ee
Indeed, the spectral decomposition of $\sigma^{-}_i(q^2)$ contains (with nonzero coefficients) only two
idempotents, $s^{(2)\uparrow i-1}$ and $c^{(2)\uparrow i-1}$:
$$
\sigma^{-}_i(q^2)=\, q\, 2_q\,
(s^{(2)\uparrow i-1}+\frac{1+q\mu}{q^3+\mu}\, c^{(2)\uparrow i-1})\ .
$$
To avoid a singularity in the expression for $\sigma_i^-(q^2)$, we have to  assume additionally
$\mu \neq -q^3$ for the rest of the proof. However, the expressions entering the relations
(\ref{idemp-1}) and (\ref{idemp-2}) are well defined and continuous at the point $\mu=-q^3$
(unless $-q^3$ coincides with one of the forbidden by eq.(\ref{mu}) values of $\mu$), so
the validity of the relations (\ref{idemp-1}) and (\ref{idemp-2})
at the point $\mu=-q^3$ follows by the continuity.

\medskip
Notice that the equalities $a^{(n)} \sigma_i=-q^{-1} a^{(n)}$ are equivalent to the equalities
$\sigma_i a^{(n)} =-q^{-1} a^{(n)}$ due to the antiautomorphism
(\ref{antvs})
since $\varsigma(a^{(n)})=a^{(n)}$ by construction.

\medskip
We now prove the equalities (\ref{idemp-1}) and (\ref{idemp-2}) by induction on $n$.

\smallskip
{}For $n=2$,~ $a^{(2)} \sigma_1  = -q^{-1} a^{(2)}$,~ by (\ref{as-2}) and (\ref{deg3-a}).

\smallskip
Let us check the equalities for some fixed $n>2$ assuming that they are valid for all smaller values of $n$.
Notice that as a byproduct of the definition (\ref{a^k}) (the first equality)
and the induction assumption, the relations (\ref{divide}) and (\ref{idemp-2}) are satisfied,
respectively, for all $i=1,2,\dots ,n-2$ and for all $m,i:\; m+i\leq n-1$. It remains to check
the relation (\ref{divide}) for
$i=n-1$ and the relation (\ref{idemp-2}) for $m=n-i$. Respectively, we calculate
\ba
\nonumber
a^{(n)}\,\sigma^{-}_{n-1}(q^2) &\sim&a^{(n-1)}\sigma^{-}_{n-1}(q^{-2n+2})
a^{(n-1)}\,\sigma^{-}_{n-1}(q^2)\\[1em]\nonumber&\sim&(a^{(n-1)}a^{(n-2)})\,\sigma^{-}_{n-1}(q^{-2n+2})
\sigma^{-}_{n-2}(q^{-2n+4})\sigma^{-}_{n-1}(q^2)\, a^{(n-2)} \hspace{10mm}
\\[1em]
\nonumber
&=&
(a^{(n-1)}\sigma^{-}_{n-2}(q^2))\,\sigma^{-}_{n-1}(q^{-2n+4})
\sigma^{-}_{n-2}(q^{-2n+2})\, a^{(n-2)} = 0 \ ,
\ea
(`$\sim$' means `proportional') and
\ba
\nonumber
a^{(n)}\, a^{(n-i)\uparrow i}&=& {q^{n-i-1}\over (n-i)_q}\,
(a^{(n)} a^{(n-i-1)\uparrow i})\,\sigma^{-}_{n-1}(q^{-2(n-i-1)})\,a^{(n-i-1)\uparrow i}
\\[1em]
\nonumber
&=&{q^{n-i-1}\over (n-i)_q}\,a^{(n)}\, (1 + q^{i-n}(n-i-1)_q )\, a^{(n-i-1)\uparrow i} = a^{(n)}\ .
\ea
Here in both cases, the definition of antisymmetrizers (\ref{a^k}) (the first equality),
induction assumption and relation (\ref{bYBE}) were used.  The centrality and primitivity of
the idempotents $a^{(n)}\in {\cal W}_n(q,\mu)$ follow then from the relations (\ref{idemp-1}).

\medskip
To prove equivalence of the two expressions for the antisymmetrizers given in the formulas
(\ref{a^k}), notice that under
conjugation by $\tau^{(i+1)}$  (\ref{tau-n}) the first expression in the formulas (\ref{a^k}) gets transformed into
the second one.
However, the elements $a^{(i+1)}$ are central in ${\cal W}_{i+1}$, so they do not change under the
conjugation which proves the consistency of the equalities (\ref{a^k}).

\medskip
All the assertions concerning the symmetrizers follow from the relations for the antisymmetrizers by
an application of the map $\iota$ (\ref{homS-A})
$$
%\be
%\lb{iota}
\iota(a^{(n)})=s^{(n)},\qquad
\iota(s^{(n)})=a^{(n)},\qquad  \iota(c^{(2n)})=c^{(2n)}\, .
%\ee
$$
the latter formulas are direct consequences of the definitions.

\medskip
The orthogonality of the antisymmetrizers and the symmetrizers is a byproduct of the relations (\ref{idemp-1}):
$$
-q^{-1}a^{(n)}s^{(m)}= (a^{(n)}\sigma_1) s^{(m)}=a^{(n)}(\sigma_1 s^{(m)})=qa^{(n)}s^{(m)}\ .
$$

\medskip
The equalities (\ref{idemp-c1}) can be proved by induction on $n$. They are obvious in the case $n=1$. Let us
check them for some fixed $n\geq 2$, assuming they are valid for all smaller values of $n$. Notice that
the iterative definition (\ref{kappa-i}) together with the induction assumption approve
the relations (\ref{idemp-c1}) for all values of index $i$, except $i=n$. Checking the case
$i=n$ splits in two subcases: $n=2$ and $n>2$. In the subcase $i=n=2$, we have $c^{(4)}=\eta^{-2}\kappa_2\kappa_3\kappa_1
\kappa_2$ and
$$
\left(c^{(4)}\right)^2 =\eta^{-4} \kappa_2\kappa_3\kappa_1\kappa_2^2\kappa_3\kappa_1\kappa_2 =
\eta^{-3}\kappa_2\kappa_3(\kappa_1\kappa_2\kappa_1)\kappa_3\kappa_2 =
\eta^{-3}\kappa_2\kappa_3\kappa_1\kappa_3\kappa_2 =\eta^{-2}\kappa_2\kappa_3\kappa_1\kappa_2 = c^{(4)}\, ,
$$
while in the subcase $i=n>2$, the calculation is carried out as follows
\ba
\nonumber
\left(c^{(2n)}\right)^2 &=&c^{(2n-2)\uparrow 1}\kappa_1\kappa_{2n-1}c^{(2n-2)\uparrow 1}
\kappa_1\kappa_{2n-1}c^{(2n-2)\uparrow 1}
\\[2pt]
\nonumber&=&
\left(c^{(2n-2)\uparrow 1}c^{(2n-4)\uparrow 2}\right) (\kappa_1\kappa_2\kappa_1)
(\kappa_{2n-1}\kappa_{2n-2}\kappa_{2n-1}) \left(c^{(2n-4)\uparrow 2}c^{(2n-2)\uparrow 1}\right)
\\[3pt]
\nonumber
&=&c^{(2n-2)\uparrow 1}\kappa_1\kappa_{2n-1}c^{(2n-2)\uparrow 1}= c^{(2n)}\, .
\ea
Here in both calculations we used the definition (\ref{kappa-i}), the induction assumption and the
relations (\ref{bmw5a}) and (\ref{bmw5b}).

\medskip
Taking into account the relations (\ref{idemp-c1}), one can derive an alternative expression for the contractors
\be
\begin{array}{ccl}
c^{(2i)}&=&c^{(2i-2)\uparrow 1}\kappa_1\kappa_{2i-1}c^{(2i-2)\uparrow 1}\,
=\, c^{(2i-2)\uparrow 1}\kappa_1\kappa_{2i-1}c^{(2i-4)\uparrow 2}\kappa_2\kappa_{2i-2}c^{(2i-4)\uparrow 2}
\\[1em]
&=&
(c^{(2i-2)\uparrow 1}c^{(2i-4)\uparrow 2})\kappa_1\kappa_{2i-1}\kappa_2\kappa_{2i-2}c^{(2i-4)\uparrow 2}\,
=\, c^{(2i-2)\uparrow 1}\kappa_{2i-1}\kappa_{2i-2}\kappa_1\kappa_2c^{(2i-4)\uparrow 2}
\\[1em]
&=&\dots\;\,\, =\, c^{(2i-2)\uparrow 1}\left(\kappa_{2i-1}\kappa_{2i-2}\dots\kappa_{i+1}\right)
\left(\kappa_1\kappa_2\dots\kappa_{i-1}\right)c^{(2)\uparrow i-1}
\\[1em]
&=&
\eta^{-1} c^{(2i-2)\uparrow 1}\left(\kappa_{2i-1}\kappa_{2i-2}\dots\kappa_{i+1}\right)
\left(\kappa_1\kappa_2\dots\kappa_i\right) .
\end{array}
\lb{kappa-i2}
\ee
Now, using this expression and noticing that, by the relations (\ref{bmw3}),
$$
\kappa_{i+1} \kappa_{i-1} \kappa_i \sigma_{i-1} =\kappa_{i+1} \kappa_{i-1} \sigma^{-1}_i =
\kappa_{i-1} \kappa_{i+1} \sigma^{-1}_i =\kappa_{i+1} \kappa_{i-1} \kappa_i \sigma_{i+1}\ ,
$$
we conclude that the equality (\ref{idemp-c2}) is satisfied for $i=n-1$. In particular, the relations (\ref{idemp-c2}) hold for $n=2$ and $i=1$. It is enough (by induction on $n$) to prove
the relations (\ref{idemp-c2}) for $i=1$.
Then observe, again by the relation (\ref{bmw3}), that
%\be
%\lb{secoide}
$$
\kappa_i\kappa_{i\pm 1}\kappa_{i\pm 2}\,\sigma_i=\kappa_i\kappa_{i\pm 1}\sigma_i\,
\kappa_{i\pm 2}=\kappa_i\sigma^{-1}_{i\pm 1}\kappa_{i\pm 2}=\sigma_{i\pm 2}\, \kappa_i\kappa_{i\pm 1}
\kappa_{i\pm 2} \ .
$$
%\ee
Now, for $n>2$,
\ba\nonumber c^{(2n)}\sigma_1&=&\eta^{-1} c^{(2n-2)\uparrow 1}
\left(\kappa_{2n-1}\kappa_{2n-2}\dots\kappa_{n+1}\right)\left(\kappa_1\kappa_2\dots\kappa_n\right)\sigma_1
\\[1em]\nonumber &=&\eta^{-1} c^{(2n-2)\uparrow 1}\left(\kappa_{2n-1}\kappa_{2n-2}\dots\kappa_{n+1}\right)
\sigma_3\left(\kappa_1\kappa_2\dots\kappa_n\right)\\[1em]\nonumber
&=&\eta^{-1} c^{(2n-2)\uparrow 1}\sigma_3\left(\kappa_{2n-1}\kappa_{2n-2}\dots\kappa_{n+1}\right)
\left(\kappa_1\kappa_2\dots\kappa_n\right)\\[1em]\nonumber &=&\eta^{-1} c^{(2n-2)\uparrow 1}\sigma_{2n-3}
\left(\kappa_{2n-1}\kappa_{2n-2}\dots\kappa_{n+1}\right)\left(\kappa_1\kappa_2\dots\kappa_n\right)
\\[1em]\nonumber &=&\eta^{-1} c^{(2n-2)\uparrow 1}
\left(\kappa_{2n-1}\kappa_{2n-2}\dots\kappa_{n+1}\right)\sigma_{2n-1}
\left(\kappa_1\kappa_2\dots\kappa_n\right) =c^{(2n)}\sigma_{2n-1}\ .\ea

\medskip
The relation (\ref{idemp-c3}) follows from the property
(\ref{bmw2a})  and the expression (\ref{kappa-i2}) (with $i=n$) for the
contractor. Then, orthogonality of  the contractors $c^{(2n)}$ with the antisymmetrizers
and the symmetrizers $a^{(m)}$, $s^{(m)}$,  $m>n$ is a corollary of the relations
(\ref{idemp-1}) and (\ref{idemp-c3}).

\medskip
A statement of the primitivity of the idempotent $c^{(2n)}\in{\cal W}_i(q,\mu)$, $i=2n,2n+1$,
goes beyond the needs of the present paper, we mention it for a sake of completeness and
postpone a purely algebraic proof till the appendix \ref{primcontr}.\hfill$\blacksquare$

\smallskip
Since the family of higher contractors does not appear to have been previously discussed in the literature, we include
Appendix B, which contains their additional properties.

%\newpage
\section{R-matrices}\lb{sec3}

Let $V$ denote a finite dimensional ${\Bbb C}$-linear space, $\dim V = \mbox{\sc n}$. Fixing some
basis $\{v_i\}_{i=1}^{\mbox{\footnotesize \sc n}}$ in $V$ we identify elements $X\in
{\rm End}(V^{\otimes n})$ with  matrices $X_{i_1 i_2 \dots i_n}^{j_1 j_2 \dots j_n}$.

\medskip
In this section we investigate properties of certain elements in ${\rm Aut}(V^{\otimes 2})$
generating representations of the braid groups ${\cal B}_n$ or, more specifically, of the Birman-Murakami-Wenzl algebras ${\cal W}_{n}(q,\mu)$ on the spaces $V^{\otimes n}$.
Traditionally such operators are called R-matrices.

\smallskip
R-matrices and compatible pairs of R-matrices are introduced in subsection \ref{subsec3.1}.
We aslo discuss there the notions of the skew-invertibility and the R-trace. Some basic
technique, useful in the work with the R-matrices, is presented in subsection
\ref{subsec3.2}.

\smallskip
A twist operation which associates a new R-matrix to a compatible pair of R-matrices,
is discussed in  subsection \ref{subsec3.2a}. We derive there an alternative expression for the
twisted R-matrix and study its skew-invertibility.

\medskip
Starting from subsection \ref{subsec3.3}, we concentrate on the R-matrices of the
BMW type.
In subsections \ref{operatorG}, \ref{Twolinearmaps} important ingredients appear: a matrix $G$
and the linear maps $\phi$ and $\xi$. As it will be explained
in  section \ref{sec4}, the matrix $G$ is responsible for the commutation relation
of the quantum matrix with a special element, called 2-contraction, of the quantum matrix
algebra. The two maps $\phi$ and $\xi$, in turn, are necessary for the definition of the $\star\,$-product of the BMW type quantum matrices, which is a proper generalization of the usual matrix multiplication to the case of matrices with noncommuting entries.

\subsection{Definition and notation}
\lb{subsec3.1}

Let $X\in {\rm End}(V^{\otimes 2})$. For any $n=2,3,\dots$ and $1\leq m \leq n-1$, denote by $X_m$ an
operator whose action on the space $V^{\otimes n}$ is given by the matrix
$$
%\be
%\lb{X-k}
(X_m)_{i_1 \dots i_n}^{j_1 \dots j_n}\ :=\ I_{i_1\dots i_{m-1}}^{j_1\dots j_{m-1}}\
X_{i_m i_{m+1}}^{j_m j_{m+1}}\ I_{i_{m+2}\dots i_n}^{j_{m+2}\dots j_n}\ .
%\ee
$$
Here $I$ denotes the identity operator. In some formulas below (see, for instance, the
equations (\ref{s-inv})$\,$) we will also use
a notation $X_{mr}\in {\rm End}(V^{\otimes n})$, $1\leq m<r\leq n-1$, referring to an operator given by
a matrix
$$
%\be
%\lb{X-kl}
(X_{mr})_{i_1 \dots i_n}^{j_1 \dots j_n}\ :=\ X_{i_m i_r}^{j_m j_r}\
I_{i_1\dots i_{m-1}i_{m+1}\dots i_{r-1} i_{r+1}\dots i_n}^{j_1\dots j_{m-1}
j_{m+1}\dots j_{r-1} j_{r+1}\dots j_n}\ .
%\ee
$$
Clearly, $X_m= X_{m\,  m+1}$.

\medskip
We reserve the symbol $P$ for the permutation operator: $P(u\otimes v)= v\otimes u \;\;\; \forall\; u, v\in V\,$. Below we repeatedly make use of relations
$$P^2 = I\, ;\quad P_{12} X_{12} = X_{21} P_{12}\, \ \;\; \forall\;\  X\in {\rm End}(V\otimes V)\, ;\quad
\tr_{(1)} P_{12} = \tr_{(3)} P_{23} = I_2\, ,$$
where the symbol ${\rm Tr}_{(i)}$ stands for the trace over an $i$-th component space in the tensor power
of the space $V$.

\medskip
An operator $X\in {\rm End}(V^{\otimes 2})$ is called {\it skew invertible} if there exists an operator
${\Psi_X}\in {\rm End}(V^{\otimes 2})$ such that
\be\tr_{(2)} X_{12} {\Psi_X}_{23} =\tr_{(2)} {\Psi_X}_{12} X_{23} = P_{13}\, .\lb{s-inv}\ee
Define two elements of $\mbox{End}(V)$
\be\lb{CandD}C_X:={\rm Tr}_{(1)}{\Psi_X}_{12}\, ,\qquad D_X:={\rm Tr}_{(2)}{\Psi_X}_{12}\, .\ee
By (\ref{s-inv}),
\be\lb{traceCD-X}\tr_{(1)} {C_X}_1 X_{12} = I_2\, ,\qquad\tr_{(2)} {D_X}_2 X_{12} = I_1\, .\ee
A skew invertible operator $X$ is called {\em strict skew invertible} if one of the matrices,
$C_X$ or $D_X$, is invertible (by lemma \ref{lemma3.5} below, if one of the
matrices, $C_X$ or $D_X$, is invertible then they are both invertible).

\medskip
An equation
$$
%\be
%\lb{YBE}
R_{1}\, R_{2}\, R_{1}\, = \, R_{2}\, R_{1}\, R_{2}\ .
%\ee
$$
for an element $R\in {\rm Aut}(V^{\otimes 2})$ is  called the {\em Yang-Baxter equation}.

\smallskip
An element $R\in {\rm Aut}(V^{\otimes 2})$ that fulfills the Yang-Baxter equation
is called an {\em R-matrix}.

\smallskip
All R-matrices in this text are assumed to be invertible.

\smallskip
Clearly, the permutation operator $P$ is the R-matrix; $R^{-1}$ is the R-matrix iff $R$ is. Any
R-matrix $R$ generates
representations $\rho_R$ of the series of braid groups ${\cal B}_n$, $n=2,3,\dots$
\be
\lb{rhoR}
\rho_R:\, {\cal B}_n\rightarrow {\rm Aut}(V^{\otimes n})\ ,\quad
\sigma_i \mapsto \rho_R(\sigma_i) = R_i, \quad 1\leq i\leq n-1 .
\ee
If additionally the R-matrix $R$ satisfies a third order {\em minimal}
characteristic polynomial (c.f. with the relation (\ref{deg3-a})$\,$)
\be\lb{charR}(qI-R)(q^{-1}I+R)(\mu I-R)=0\ ,\ee
and an element
\be\lb{K}K := \mu^{-1} (q-q^{-1})^{-1}\, (qI-R)(q^{-1}I+R)\ee
fulfills  conditions
\be\lb{bmwRa}K_2\, K_1 \, = \,  R_1^{\pm 1}\, R_2^{\pm 1}\, K_1\,
\ee
and
\be\lb{bmwRb}
K_1\, K_2\, K_1 \, = \, K_1\, ,\ee
then we call $R$ an R-matrix of a {\em BMW type}~ (c.f. with eqs.(\ref{kappa})--(\ref{bmw5b});
we make a different but equivalent choice of defining relations).

\medskip
{}For an R-matrix of the BMWtype, the formulas (\ref{rhoR}) define representations of the algebras ${\cal W}_n(q,\mu)
\rightarrow {\rm End}(V^{\otimes n})$, $n=2,3,\dots $ . In particular, $\rho_R(\kappa_i)=K_i$.

\medskip
An ordered pair $\{ R, F\}$ of two operators $R$ and $F$ from ${\rm End}(V^{\otimes 2})$ is
called {\em a compatible pair} if conditions
\be R_1\, F_2\, F_1\, =\, F_2\, F_1\, R_2\, ,\qquad R_2\, F_1\, F_2\, =\, F_1\, F_2\, R_1\, ,\lb{sovm}\ee
are satisfied. If, in addition, $R$ and $F$ are R-matrices, the pair $\{ R, F\}$ is called
a compatible pair of R-matrices. The equalities (\ref{sovm}) are called {\em twist relations} (on the
notion of the twist see \cite{D2,Resh2,IOP2}). Clearly, $\{ R,P\}$ and $\{ R,R\}$ are compatible pairs
of R-matrices;
pairs $\{ R^{-1},F\}$ and $\{ R,F^{-1}\}$ are compatible iff the pair $\{R,F\}$ is.

\begin{defin}\lb{definition3.1}
Consider a  space of ~$\mbox{\sc n}\times \mbox{\sc n}$ matrices
${\rm Mat}_{\mbox{\footnotesize\sc n}}(W)$, whose entries belong to some $\Bbb C$-linear space $W$.
Let $R$ be a skew invertible R-matrix. A linear map
$$
%\be
%\lb{r-sled}
{\rm Tr\str{-1.3}}_R:\; {\rm Mat}_{\mbox{\footnotesize\sc n}}(W)\,\rightarrow \, W ,\qquad
{\rm Tr\str{-1.3}}_R(M) =\sum_{i,j=1}^{\mbox{\footnotesize\sc n}}{(D_R)}_i^jM_j^i\, , \qquad
M\in{\rm Mat}_{\mbox{\footnotesize\sc n}}(W)\,,
%\ee
$$
is called an  R-trace. \end{defin}

The relation (\ref{traceCD-X}) in this notation reads
\be
\lb{traceR}
\Tr{2} R_{12} = I_1\, .
\ee

\subsection{R-technique}
\lb{subsec3.2}

In this and the next subsections we develop a technique for dealing with the R-matrices, their
compatible pairs and the R-trace. Most of results reported here, like lemma \ref{lemma3.5}
and, in a particular case of a compatible pair $\{ R,R\}$ -- lemmas \ref{lemma3.2}
and \ref{lemma3.3} and the corollary \ref{corollary3.4} -- are rather well known (see,
e.g., \cite{I,O}). However, we often use them in a more general setting and so, when necessary,
we present sketches of proofs.

\smallskip
Proposition \ref{proposition3.6} contains new results.
Here we derive an expression, different from the standard one, for the
twisted R-matrix, which helps to investigate its skew-invertibility.

\medskip
A universal (i.e., quasi-triangular Hopf algebraic) content
of the matrix relations derived in this and the next subsections is discussed in the
appendix \ref{append1}.

\begin{lem}\lb{lemma3.2}
Let $\{X,F\}$ be a compatible pair,
where $X$ is skew invertible. Let ${\rm Mat}_{\mbox{\footnotesize\sc n}}(W)$ be as in the definition
\ref{definition3.1}. For any $M\in {\rm Mat}_{\mbox{\footnotesize\sc n}}(W)$, one has
\ba\lb{inv-trC}\tr_{(1)} \left( {C_X}_1 F_{12}^{\varepsilon}\, M_2\, F_{12}^{-\varepsilon}\right)
&=& I_2\, \tr (C_X M)\, ,\\[1em]\lb{inv-trD}\tr_{(2)} \left({D_X}_2 F_{12}^{-\varepsilon}\,
M_1\, F_{12}^{\varepsilon}\right) &=& I_1\ \tr(D_X M)\,
\ea
for $\varepsilon =\pm 1$.\end{lem}

\nin {\bf Proof.~} We use the twist relations  (\ref{sovm}) in a form
$$F_{23}^{\varepsilon}\, X_{34}\, F_{23}^{-\varepsilon}\, =\, F_{34}^{-\varepsilon}\,
X_{23}\, F_{34}^{\varepsilon}\, ,\quad \varepsilon=\pm 1\, .$$
Multiplying it by $({\Psi_X}_{12}{\Psi_X}_{45})$ and taking the traces in the spaces 2 and 4, we get
\be\lb{trpsifpf}\tr_{(2)}({\Psi_X}_{12}\ F_{23}^{\varepsilon}\ P_{35}\ F_{23}^{-\varepsilon})\ =\
\tr_{(4)}({\Psi_X}_{45}\ F_{34}^{-\varepsilon}\ P_{13}\ F_{34}^{\varepsilon})\ .\ee
Here the relation (\ref{s-inv}), defining the operator $\Psi_X$, was applied to calculate the
traces. Now taking the trace in the space number 1 or number 5, we obtain (after relabeling)
\ba
\lb{tr-0}
\tr_{(1)}({C_X}_1\ F_{12}^{\varepsilon}\ P_{23}\ F_{12}^{-\varepsilon})&=&
{C_X}_3\, I_2\ ,
\\[1em]
\lb{tr-4}
\tr_{(3)}({D_X}_3\ F_{23}^{-\varepsilon}\
P_{12}\ F_{23}^{\varepsilon})
&=&  {D_X}_1\, I_2\, .
\ea
These two relations are equivalent forms of the relations (\ref{inv-trC}) and (\ref{inv-trD}).
For example, the formula (\ref{inv-trC}) is obtained by multiplying the relation (\ref{tr-0}) by
the operator $M_3$ and taking the trace in the space 3.
\hfill$\blacksquare$

\begin{lem}\lb{lemma3.3}
Let $\{X,F\}$ be a compatible pair of skew invertible operators $X$ and $F$. Then the following relations
\ba\lb{psi-C} {C_X}_1\, {\Psi_F}_{12}\, =\, F_{21}^{-1}\, {C_X}_2\, , &&
{\Psi_F}_{12}\, {C_X}_1\, =\, {C_X}_2\, F_{21}^{-1}\, ,\\[1em]\lb{psi-D}
{\Psi_F}_{12}\, {D_X}_2\, = \, {D_X}_1\, F_{21}^{-1}\, , &&
{D_X}_2\, {\Psi_F}_{12}\, =\, F_{21}^{-1}\, {D_X}_1\, \ea
hold.\end{lem}

\nin {\bf Proof.~} For a skew invertible operator $F$, the relations (\ref{psi-C}) and
(\ref{psi-D}) are equivalent to the relations (\ref{tr-0}) and (\ref{tr-4}). Let us demonstrate
how  the left one of the relations (\ref{psi-C})
is derived from the relation (\ref{tr-0}) with $\varepsilon = 1$.

\smallskip
Multiply the relation (\ref{tr-0}) by a combination $(P_{23}{\Psi_F}_{24})$ from the right, take
the trace in the space 2 and simplify
the result using the relation (\ref{s-inv}) for $X=F$ and the properties of the permutation
$$\tr_{(1)}({C_X}_1\, P_{14}\, F_{13}^{-1})\, =\,
{C_X}_3\, \tr_{(2)}(P_{23}\, {\Psi_F}_{24})\, =\, {C_X}_3\, {\Psi_F}_{34}\, .$$
Then simplify the left hand side of the equality using the cyclic property of the trace
$$\tr_{(1)}({C_X}_1\, P_{14}\, F_{13}^{-1})\, =\, \tr_{(1)}( P_{14}\, F_{13}^{-1}\, {C_X}_1)\,
=\, F_{43}^{-1}\, {C_X}_4\, \tr_{(1)}P_{14}\, =\, F_{43}^{-1}\, {C_X}_4\, .$$
This proves the left relation in (\ref{psi-C}).\hfill$\blacksquare$

\begin{cor}\lb{corollary3.4}
Let $\{X,F\}$ and $\{Y,F\}$ be compatible pairs of skew invertible
operators  $X$,  $Y$ and $F$. Then the following relations
\ba\lb{FCC} &F_{12}\, {C_X}_1 {C_Y}_2  = {C_Y}_1 {C_X}_2 F_{12}\ ,
\qquad\quad  F_{12}\, {D_X}_1 {D_Y}_2 \, =\, {D_Y}_1 {D_X}_2 F_{12}\ ,&
\hspace{10mm}\\[1em]
\lb{FCD} &F_{12}\, (C_X D_Y)_2 = (C_X D_Y)_1\, F_{12}\ ,\qquad\qquad F_{12}\,
(D_Y C_X)_1 \, =\, (D_Y C_X)_2\, F_{12}\, ,&\ea
\ba\lb{CxDf}
&{\tr}_{(1)}({C_X}_1 F_{12}^{-1}) = (C_X D_F)_2 = (D_F C_X)_2\ ,&\\[1em]
\lb{DxCf}
&{\tr}_{(2)}({D_X}_2 F_{12}^{-1}) = (C_F D_X)_1 = (D_X C_F)_1\,\ &\ea
hold.\end{cor}

\nin {\bf Proof.~} A calculation
$( F_{12}^{-1} {C_Y}_1) {C_X}_2 ={C_Y}_2 ( {\Psi_F}_{21} {C_X}_2) =
{C_Y}_2 {C_X}_1 F_{12}^{-1}={C_X}_1 {C_Y}_2 F_{12}^{-1}$
proves the left one of the relations (\ref{FCC}). Here the relations (\ref{psi-C}) were applied.

\smallskip
A calculation $(F_{12}^{-1} {C_X}_1) {D_Y}_1 = {C_X}_2 (\Psi^F_{21} {D_Y}_1) =
{C_X}_2 {D_Y}_2 F_{12}^{-1}$
proves  the left one of the relations (\ref{FCD}). Here one uses subsequently the left equations
from (\ref{psi-C}) and (\ref{psi-D}).

\smallskip
The relations (\ref{CxDf}) follow by taking  $\tr_{(2)}$ of the equations (\ref{psi-C}).

\smallskip
The rest of the relations in (\ref{FCC})--(\ref{DxCf}) are derived in a similar way.\hfill$\blacksquare$

\begin{lem}\lb{lemma3.5}
Let $X$ be a skew invertible R-matrix. Then statements\par
a) the R-matrix $X^{-1}$ is skew invertible;\par
b) the R-matrix $X$ is strict skew invertible,\par\noindent
are equivalent.

Provided these statements are satisfied, both $C_X$ and $D_X$ are invertible and one has
\ba\lb{CDinv} &&C_{X^{-1}} =\,  D_X^{-1}\, , \qquad D_{X^{-1}} =\,  C_X^{-1}\, .\ea\end{lem}

\nin {\bf Proof.~} See \cite{O},  section 4.1, statements after eq.(4.1.77), or \cite{I},
proposition 2 in  section 3.1. \hfill$\blacksquare$

\medskip
Under an assumption of an existence, for an R-matrix $X$, of the operators $X^{-1}$, ${\Psi_{\! X\,}}$ and
${\Psi_{\! X^{-1}\,}}$, the relations (\ref{CDinv}) were proved in \cite{Resh}.

\medskip
Since, for a compatible pair $\{ X,F\}$, the pair $\{ X,F^{-1}\}$ is also compatible,
the formulas (\ref{CDinv}) together with the relations (\ref{CxDf}-\ref{DxCf}) imply that
$C_XC_F=C_FC_X$ and $D_XD_F=D_FD_X$.

\subsection{Twists}
\lb{subsec3.2a}

Let $\{R,F\}$ be a compatible pair of R-matrices. Define a {\em twisted} operator
\be
\lb{R_f}
R_f := F^{-1} R F\, .
\ee
It is well known that
$R_f$ is an R-matrix and the pair $\{ R_f , F\}$ is compatible. Therefore, one can twist again;
in \cite{IOP1} it was shown that if $F$ is skew invertible then
\be
\lb{rcdm}
{D_F}_1\, {D_F}_2\, ((R_f)_f)_{12}\, =\, R_{12}\, {D_F}_1\,  {D_F}_2\quad {\rm{and}}
\quad  {C_F}_1\, {C_F}_2\, ((R_f)_f)_{12}\, =R_{12}\, {C_F}_1\, {C_F}_2\ .
\ee
A comparison of two equalities in eq.(\ref{rcdm}) shows that
\be
\lb{rgrel}
[\,  R_{12}\, ,\, (C^{-1}_F\,  D_F)_1\,  (C^{-1}_F\,  D_F)_2\, ] =0\ .
\ee

\begin{prop}
\lb{proposition3.6}
Let $\{R,F\}$ be a compatible pair of R-matrices. The following statements hold:
\begin{itemize}
\item[{\em a)}] if $F$ is strict skew invertible then the twisted R-matrix $R_f$,
defined by the formula (\ref{R_f}), can be expressed in a form
\be
\lb{R_f-fin}
{R_f}_{12}\, =\,
\tr_{(34)}\left( F_{32}^{-1} {C_{F^{-1}}}_3 R_{34} {D_F}_4 F_{14}\right) \, ;
\ee
\item[{\em b)}]
if $R$ is skew invertible and $F$ is strict skew invertible then $R_f$ is skew invertible; its skew
inverse is
\be
\lb{Psi_R_f}
{\Psi_{\! R_f}}_{12}\, =\,
{C_{F^{-1}}}_2\, \tr_{(34)}\left( F_{23}^{-1} {\Psi_{\! R\,}}_{34} F_{41}\right) {D_F}_1\, ;
\ee
moreover, $\Psi_{\! R_f}$ can be expressed in a form
\be
\lb{Psi_R_f-another}
{\Psi_{\! R_f}}_{12}\, =\, {C_{F^{-1}}}_2 F_{21}{D_{F^{-1}}}_2
{\Psi_{\! R\,}}_{12}{C_F}_1 F^{-1}_{21}{D_F}_1\ ;
\ee
\item[{\em c)}] under the conditions in {\em b)},
\be
\lb{CDtwist}
C_{R_f} = C_{F^{-1}} D_R\, C_F\, , \quad D_{R_f} = D_{F^{-1}} C_R D_F\,
\ee
(thus, if, in addition to the conditions in {\em b)}, $R$ is strict skew invertible then
$R_f$ is strict skew invertible as well).
\end{itemize}
\end{prop}

\nin {\bf Proof.~} To verify the assertion a) we calculate
\be\!\!\begin{array}{ccl} {R_f}_{12}&=& (F^{-1} R F)_{12}
= F^{-1}_{12}\left(\tr_{(4)} F^{-1}_{41} {C_{F^{-1}}}_4\right)(R F)_{12}\\[1em]
&=&\tr_{(4)} \left( (R F)_{41} F^{-1}_{12} F^{-1}_{41} {C_{F^{-1}}}_4\,\right) =
\left(\tr_{(3)}P_{13}\right)\tr_{(4)} \left( (R F)_{41} F^{-1}_{12} {C_{F^{-1}}}_1
{\Psi_F}_{14}\right)\\[1em] &=&\tr_{(34)} \left( (R F)_{43} F^{-1}_{32} {C_{F^{-1}}}_3
P_{13} {\Psi_F}_{14}\right) =\tr_{(3)}\left( F^{-1}_{32} {C_{F^{-1}}}_3
\underline{P_{13} \tr_{(4)}{\Psi_F}_{14}(R F)_{43}}\right) \ ,\end{array}\lb{R_f-alt} \ee
where in the second equality we used the relation (\ref{traceCD-X}) for $X=F^{-1}$; in the third
equality we applied the twist relations for the compatible pairs $\{ R,F\}$ and $\{ F,F\}$;
in the fourth equality we applied the relations
(\ref{psi-C}) for $X=F^{-1}$ and inserted the identity operator $\tr_{(3)}P_{13}$; in the fifth
equality we permuted the operator $P_{13}$ rightwards and then, in the sixth equality, used the
cyclic property of the trace to move the combination $(RF)_{43}$ to the right.

\smallskip
To complete the transformation, we derive an alternative form for the underlined
expression in the last line in eq.(\ref{R_f-alt}). Multiplying the twist relation
$R_2 F_3 F_2 = F_3 F_2 R_3$ by a combination
$({\Psi_F}_{12}{D_F}_4)$ and taking the traces in the spaces 2 and 4, we obtain (using the formulas (\ref{s-inv}) and
(\ref{traceCD-X}) for $X=F$)
$$\tr_{(2)}\left( {\Psi_F}_{12}(R F)_{23}\right) \, =\,
\tr_{(4)} \left( {D_F}_4 F_{34} P_{13} R_{34}\right) \, ,$$
which is equivalent  (multiply by $P_{13}$ from the left and use the cyclic property of the trace) to
\be\lb{alternative} P_{13}\tr_{(2)}\left( {\Psi_F}_{12}(R F)_{23}\right) \, =\,
\tr_{(4)}\left( R_{34} {D_F}_4 F_{14}\right) \, .\ee
Now, substituting the equality (\ref{alternative}) into the last line of the
calculation (\ref{R_f-alt}), we finish the transformation and
obtain the formula (\ref{R_f-fin}).

\medskip
Given the formula for $R_f$, the calculation of ${\Psi_{R_f}}$ becomes straightforward
and one finds the formula (\ref{Psi_R_f}).

\smallskip
Thus, the skew invertibility of $R_f$ is established.

\medskip
Now we derive the expression (\ref{Psi_R_f-another}) for $\Psi_{\! R_f}$. Multiplying the
equality (\ref{trpsifpf}) with $\varepsilon =1$ by a combination $P_{35}{D_{F^{-1}}}_5$ from the
right and taking the trace in the space 5, we obtain
$$\begin{array}{rcl}\tr_{(2)}({\Psi_R}_{12}F_{23})&=&\tr_{(45)}({\Psi_R}_{45}F^{-1}_{34}P_{13}F_{34}
P_{35}{D_{F^{-1}}}_5)\\[1em] &=&\tr_{(4)}(F^{-1}_{34}P_{13}F_{34}{D_{F^{-1}}}_3{\Psi_R}_{43})\ .
\end{array}$$
Substituting this into the expression (\ref{Psi_R_f}), we find
\be\lb{lexfpsif}\begin{array}{rcl} {\Psi_{\! R_f}}_{12}\, &=&\,
{C_{F^{-1}}}_2\, \tr_{(34)}\left( F_{23}^{-1}F_{14}^{-1}P_{13} F_{14}{D_{F^{-1}}}_1 {\Psi_{\! R\,}}_{41}
\right) {D_F}_1\\[1em] &=&{C_{F^{-1}}}_2\, \tr_{(4)}\left(F_{14}^{-1}\tr_{(3)}(F_{23}^{-1}P_{13})
{}F_{14}{D_{F^{-1}}}_1 {\Psi_{\! R\,}}_{41}\right) {D_F}_1\\[1em]
&=&{C_{F^{-1}}}_2\, \tr_{(4)}\left(F_{14}^{-1}F_{21}^{-1}
{}F_{14}{D_{F^{-1}}}_1 {\Psi_{\! R\,}}_{41}\right) {D_F}_1\\[1em]
&=&{C_{F^{-1}}}_2\, F_{21}\, \tr_{(4)}\left(F_{14}^{-1}F_{21}^{-1}
{D_{F^{-1}}}_1 {\Psi_{\! R\,}}_{41}\right) {D_F}_1\\[1em]
&=&{C_{F^{-1}}}_2\, F_{21}\, {D_{F^{-1}}}_2\tr_{(4)}\left(F_{14}^{-1}{\Psi_F}_{12}
{\Psi_{\! R\,}}_{41}\right) {D_F}_1\ .\end{array}\ee
We used the Yang-Baxter equation for the operator $F$ in the fourth equality and the relations
(\ref{psi-C}) in the fifth equality.

\smallskip
Multiplying eq.(\ref{trpsifpf}) with $\varepsilon =-1$ by ${\Psi_F}_{01}P_{13}$ from the left and
by $P_{35}{\Psi_F}_{56}$ from the right and taking the traces in the spaces 1 and 5, we find
$$F_{14}^{-1}{\Psi_F}_{12}{\Psi_{\! R\,}}_{41}={\Psi_{\! R\,}}_{12}{\Psi_F}_{41} F_{21}^{-1}\ .$$
Substituting this into the last line of the calculation (\ref{lexfpsif}), we obtain the equality (\ref{Psi_R_f-another}).

\medskip
{}Finally, the expressions (\ref{CDtwist}) for the operators $C_{R_f}$ and $D_{R_f}$ are
obtained by taking the trace in
the space 1 or the space 2 of the expression (\ref{Psi_R_f}) for the skew inverse of the
twisted R-matrix and the subsequent use of the relations (\ref{traceCD-X}) for
$X=F^{\pm 1}$ and the relations (\ref{CandD}), (\ref{CxDf}) and (\ref{DxCf}) for $X=R$. \hfill$\blacksquare$

\begin{rem}{\rm If one uses the expression (\ref{R_f-fin}) for the twisted R-matrix then the
relation (\ref{rcdm}) becomes straightforward:
$$
%\be
\begin{array}{l} ((R_f)_f)_{12} =
\tr_{(3456)}\left( F_{32}^{-1} (D_F^{-1})_3 F_{54}^{-1} (D_F^{-1})_5 R_{56}
{D_F}_6 F_{36} {D_F}_4 F_{14}\right)
\\[1em]
\ \ \
= \tr_{(3456)}\left(
\underline{\underline{F_{54}^{-1}{D_F}_4F_{14}}}(D_F^{-1})_5 R_{56}{D_F}_6 \underline{F_{32}^{-1}(D_F^{-1})_3F_{36}}\right)
\\[1em]
\ \ \
=\tr_{(56)}\left( P_{15}{D_F}_1(D_F^{-1})_5 R_{56}{D_F}_6 (D_F^{-1})_2 P_{26}\right)
\\[1em]
\ \ \ =(D_F^{-1})_2\tr_{(56)}\left( P_{15}(D_F^{-1})_5 R_{56}{D_F}_6 P_{26}\right)
{D_F}_1%\,
\\[1em]
=(D_F^{-1})_1(D_F^{-1})_2 R_{12}{D_F}_1{D_F}_2\ .
\end{array}
%\ee
$$
In the first equality we applied the formula (\ref{R_f-fin}) twice and replaced the operators
$C_{F^{-1}}$ by $D_F^{-1}$ by the relation (\ref{CDinv}); in the second equality we
collected together the terms involving the space number 3 (they are underlined) and the
terms involving the space number 4 (they are underlined twice); in the third equality we
evaluated the traces in the spaces 3 and 4 using the relations from lemma \ref{lemma3.3};
in the fourth equality we moved the operator $(D_F^{-1})_2$ leftwards out of the trace
and the operator ${D_F}_1$ rightwards out of the trace; in the fifth equality we transported
the operator $P_{15}$ rightwards and the operator $P_{26}$ leftwards under the trace and
then evaluated the remaining traces in the spaces 5 and 6.}
\end{rem}

\subsection{BMW type R-matrices}
\lb{subsec3.3}

In this subsection we discuss the R-matrices of the BMW type in more detail.

\smallskip
In lemma \ref{lemma3.7} we collect additional relations specific to the BMW type
R-matrices. Based on these formulas, we will introduce later, in subsections \ref{operatorG}
and \ref{Twolinearmaps}, an invertible operator $G\in {\rm Aut}(V)$ and linear maps $\phi$ and
$\xi$, which will be used in  section \ref{sec4} for a definition of a product of
quantum matrices and for a quantum matrix inversion.

\begin{lem}\lb{lemma3.7}
Let $R$ be a skew invertible R-matrix of the BMW type. Then
\begin{itemize}
\item[$\bullet$] the operator $R$ is strict skew invertible;
\item[$\bullet$] the rank of the operator $K$ equals one, ${\rm rk}\, K = 1$;
\item[$\bullet$] the following relations
\be\lb{traceK}\tr_{(2)} K_{12} = \mu^{-1} {D_R}_1\ \  ,\ \ \tr_{(1)} K_{12} = \mu^{-1} {C_R}_2\ ,\ee
\be\lb{traceDK} \Tr{2} K_{12} = \mu\, I_1\, ,\ee
\be\lb{traceD}\tr_{\!\! R} I =
\mu\, \eta \,\equiv\,{\displaystyle {(q-\mu)(q^{-1}+\mu)\over (q-q^{-1})}}\, ,\ee
\be\lb{C*D}
C_R D_R = \mu^2 I\, ,\ee
\be\lb{KDD} K_{12} {D_R}_1 {D_R}_2 = {D_R}_1 {D_R}_2 K_{12}\
=\ \mu^2 K_{12}\, \ee
hold.
\end{itemize}\end{lem}

\nin {\bf Proof.~} The proof of all the statements in the lemma but the last one is given in \cite{IOP3}.

\smallskip
The last relation (\ref{KDD}) (which, in another form, figures in \cite{IOP3}, in  proposition 2) can
be established in the following way.

The first equality in (\ref{KDD}) is a consequence of a relation
\be\lb{RDD}R_{12}{D_R}_1 {D_R}_2\, =\, {D_R}_1 {D_R}_2 R_{12}\, ,\ee
which is just the equality (\ref{FCC}) written for the pair $\{R,R\}$. Then the conditions
$K^2 \sim K$ and ${\rm rk}\, K = 1$ together imply
$K_{12} {D_R}_1 {D_R}_2 \sim K_{12} {D_R}_1 {D_R}_2 K_{12} \sim K_{12}\,$.
A coefficient of proportionality in this relation is recovered by taking the trace of it in the space 2
and the subsequent use of the relations (\ref{traceK}) and (\ref{traceDK}). \hfill$\blacksquare$

\medskip
In \cite{IOP3}, a pair of mutually inverse matrices
\be
{E}_2 :=\tr_{(1)} (K_{12}P_{12})
\ \ {\mathrm{and}}\ \ E^{-1}_1 :=\tr_{(2)} (K_{12}P_{12})
\lb{defcaxy}
\ee
was introduced (see eqs.(32) and (33) and  proposition 2 in \cite{IOP3}).

\smallskip%%%syuda perestavil lemmu
We shall now collect several useful identities involving the operators $K$
and $E$.

\begin{lem}
\lb{lemma-K-prop}
{\rm (a)} The following relations
\ba
K_{12}K_{23}=E_3\; K_{12}P_{23}P_{12}& ,&
K_{23}K_{12}=E^{-1}_1\; K_{23}P_{12}P_{23}\ ,
\lb{k12k23}
\\[3pt]
%\ee
%\be
K_{13}K_{23}=\mu^{-1}{D_R}_2\; K_{13}P_{12}& ,&
K_{12}K_{13}=\mu^{-1}{C_R}_3\; K_{12}P_{23}\ ,
\lb{k13k23}
\\[3pt]
\nn
%\ee
%\be
K_{23}K_{14}P_{12}P_{34}=K_{23}K_{14}& ,&
%\lb{kkpp1}
%\ee
%\be
K_{23}K_{14}P_{13}P_{24}=K_{23}K_{14}P_{23}P_{14}\ ,
%\lb{kkpp2}
%\ee
\\[3pt]
\nn
%\be
K_{12}E_1^{-1}=\mu^{-1}K_{12}P_{12}{D_R}_1 & ,&
E_1 K_{12}=\mu^{-1}{D_R}_1P_{12} K_{12}\
%\lb{rekk21}
%\ee
\ea
hold.

\smallskip
{\rm (b)} We have
$$
%\be
K_{12}E_1 E_2= E_1 E_2K_{12}=K_{12}\ .
%\lb{kcax}
%\ee
$$

\smallskip
{\rm (c)} The operator $K$ is skew invertible, its skew inverse is
$$
%\be
{\Psi_K}_{12}=E_1 K_{12} E_2=\mu^{-2}{D_R}_1 K_{21}{D_R}_1\ .
%\lb{skinK}
%\ee
$$
\end{lem}

\nin {\bf Proof.~} (a) All these identities
follow from the rank one property of the operator $K$ (written explicitly, with indices, they become evident).

\smallskip
(b) To verify, for instance, that $K_{12}E_1^{-1} E_2^{-1}=K_{12}$, use the definition (\ref{defcaxy}) of the matrix
$E_2^{-1}$,  $E_2^{-1}=\tr_{(3)} (K_{23}P_{23})$, and then the relation (\ref{k12k23}) to remove the trace.

\smallskip
(c) This follows from the identities in (a) in the lemma.\hfill$\blacksquare$

\begin{rem}
{\rm The relations (\ref{k12k23}) admit the following generalizations:
$$
%\be
K_1K_2\dots K_j\, =\, E_3E_4\dots E_{j+1}\cdot (P_1P_2\dots P_j)^2 K_j\ ,
%\lb{strk1}
%\ee
$$
$$
%\be
K_j\dots K_2K_1\, =\, E_1^{-1}E_2^{-1}\dots E_{j-1}^{-1}K_j\cdot (P_j\dots P_2P_1)^2\ .
\lb{strk2}
%\ee
$$
The relations (\ref{k13k23}) admit the following generalizations:
$$
%\be
K_{10}K_{20}\dots K_{j0}\, =\, \mu^{1-j}({D_R}_2{D_R}_3\dots {D_R}_j)\cdot
(P_1P_2\dots P_{j-1})\ K_{j0}\ ,
%\lb{strk3}
%\ee
$$
$$
%\be
K_{01}K_{02}\dots K_{0j}\, =\, \mu^{1-j}({C_R}_2{C_R}_3\dots {C_R}_j)\cdot
(P_1P_2\dots P_{j-1})\ K_{0j}\ .
%\lb{strk4}
%\ee
$$
In all four formulas above
%(\ref{strk1})-(\ref{strk4}),
$j$ is an arbitrary positive integer. These relations
%(\ref{strk1})-(\ref{strk4})
can be proved by induction on $j$.
}\end{rem}

\subsection{Operator $G$}
\label{operatorG}

In the following lemma, we define analogues of the matrices $E$ and $E^{-1}$ for a compatible pair
$\{R,F\}$ of R-matrices. When the operator $F$ is the permutation operator, $F=P$, the
matrix $G$ of the definition-lemma \ref{lemma3.8} coincides with the matrix $E$.

\begin{def-lem}\lb{lemma3.8}
Let $\{R,F\}$ be a compatible pair of R-matrices, where
$R$ is skew-invertible of the BMW type and $F$ is strict skew-invertible.
Define an element $G\in {\rm End}(V)$ by
\be\lb{G}G_1 \, :=\, \tr_{(23)} K_2 F_1^{-1} F_{2}^{-1}\, .\ee
The operator $G$ is invertible, the inverse operator reads
\be\lb{G-inv}G_1^{-1}\, =\, \tr_{(23)} F_2 F_1 K_2\, .\ee

\noindent
{}The following relations
\ba
\lb{R-G}
 &&R_{12} G_1 G_2\, =\, G_1 G_2 R_{12}\, ,
 \\[1em]
\lb{F-G}
&&{}F_{12}^{\varepsilon} G_1 = G_2 F_{12}^{\varepsilon}\ \ \ \mbox{for}\;\;\varepsilon=\pm 1\, ,\\[1em]
\lb{comm-Ga}
&&[D_R, G]\, =\, 0\, ,
\\[1em]
\lb{comm-Gb}
&&[ C_F, G]\, =\, [D_F, G]\, =\, 0\, ,
\\[1em]
\nn
%\lb{comm-Gc}
&&[E, G]\, =\, 0\, ,
\\[1em]
\nn
%\lb{comm-EDf}
&&[ C_F, E]\, =\, [D_F, E]\, =\, 0\,
\ea
are satisfied.
\end{def-lem}

\nin {\bf Proof.~} A check of the invertibility of $G$ is a direct calculation
\be
\begin{array}{ccc}
G_1\, G_1^{-1} &=& (\tr_{(23)}K_2 F_1^{-1} F_2^{-1})
(\tr_{(23)}F_2 F_1 K_2)\, =\,\tr_{(23)} K_2 F_1^{-1} F_2^{-1} K_2 F_2 F_1\\[1em] &=& \tr_{(23)}K_2
{}F_1^{-1}{K_f}_2 F_1\, =\,\tr_{(23)} K_2 F_2 {K_f}_1 F_2^{-1}\, =\, \tr_{(23)} {K_f}_2 {K_f}_1 \, =\, I\, .
\end{array}
\lb{gginvpr}
\ee
Here in the first line we used the formulas (\ref{G}) and (\ref{G-inv}) and the property ${\rm rk}\, K=1$:
if $\Pi =|\zeta\rangle\langle\psi |$ is a rank one projector then $\tr (\Pi A)=\langle\psi |A|\zeta\rangle$
for any operator $A$ and
$$\tr (\Pi A)\,\tr (\Pi B)= \langle\psi |A|\zeta\rangle\, \langle\psi |B|\zeta\rangle=
\langle\psi |A\Pi B|\zeta\rangle =\tr (\Pi A\Pi B)\ $$
for any $A$ and $B$; in the second line of the calculation (\ref{gginvpr}) we passed from $K$ to
$K_f = F^{-1} K F$ and applied the twist relations (for the operators $K_f$ and $F$) and the
cyclic property of the trace. In the last equality
of (\ref{gginvpr}) we evaluated the traces using the relations (\ref{traceK}) and then the
relation
(\ref{traceDK}) for the operator $K_f$ (we are allowed to use these relations because
the operator $R_f$ is skew-invertible
by  proposition \ref{proposition3.6}).

\medskip
Notice that, in view of the relation (\ref{KDD}), we can rewrite the formula for the operator
$G$  using the R-traces instead of the ordinary ones
\be\lb{RG}G_1 \, =\, \mu^{-2}\,\Tr{23} K_2 F_1^{-1} F_{2}^{-1}\, .\ee
Applying the formula (\ref{inv-trD}) (written for $F^{\varepsilon}=X=R$) twice to this equality, we begin
our next calculation
\be
\begin{array}{ccl}
G_1 I_2 &=& \mu^{-2}\, \Tr{34} (R_2 R_3) K_2 F_1^{-1} F_2^{-1}
(R_3^{-1} R_2^{-1})\, =\,  \mu^{-2}\, \Tr{34} K_3 K_2 F_1^{-1} F_2^{-1}
R_3^{-1} R_2^{-1}\\[1em] &=& \mu^{-2}\, \Tr{34}K_2 F_1^{-1} F_2^{-1} K_2 K_3
\, =\, \mu^{-1}\, \Tr{3}K_2 F_1^{-1} F_2^{-1} K_2\, .\end{array}
\lb{GI}
\ee
Here we used the relation (\ref{bmwRa}) in the last equality of the first line. In the second
line we again applied the relation (\ref{bmwRa}) after moving the operator $K_3$ to the right
(for that we need the relation (\ref{KDD}) and the cyclicity
of the trace) and then we evaluated one R-trace with the help of the relation (\ref{traceDK}).

\medskip
Now we use the formula (\ref{GI}) for the product $G_1 G_2$ in a transformation
\ba
\nonumber
&& \hspace{-5mm} G_1 G_2 R_1\, =\, \mu^{-2}\, \Tr{34}(K_3 F_2^{-1} F_3^{-1} K_3)
(K_2 F_1^{-1} F_2^{-1} K_2) R_1\\[1em]\nonumber && =\, \mu^{-2}\,
\Tr{34}F_2^{-1} F_3^{-1} K_2 K_3 K_2 F_1^{-1} F_2^{-1} K_2 R_1\, =\,
\mu^{-2}\, \Tr{34}F_2^{-1} F_3^{-1} F_1^{-1} F_2^{-1} K_1 K_2 R_1\\[1em]\nonumber
&& =\,  \mu^{-2}\,\Tr{34}F_2^{-1} F_3^{-1} F_1^{-1} F_2^{-1} K_1 R_2^{-1}\, =\,
\mu^{-2}\,\Tr{34} K_3 F_2^{-1} F_3^{-1} F_1^{-1} F_2^{-1}  R_2^{-1}\\[1em]\nonumber
&& =\, \mu^{-2}\,\Tr{34}F_2^{-1} F_1^{-1} F_3^{-1} F_2^{-1}  R_3 K_2 K_3\, =\, \mu^{-2}
R_1\, \Tr{34}F_2^{-1} F_1^{-1} F_3^{-1} F_2^{-1} K_2 K_3\\[1em]\nonumber
&& =\, \mu^{-2} R_1\, \Tr{34}K_3 F_2^{-1} F_3^{-1} F_1^{-1} F_2^{-1}  K_2\, =\, R_1 G_1 G_2\, ,
\ea
which demonstrates the
relation (\ref{R-G}). While doing the above calculation, we repeatedly
used the twist relations for the pairs $\{K,F^{-1}\}$ and $\{R,F^{-1}\}$, applied the formulas (\ref{bmwRa}) and (\ref{bmwRa})
and exploited the cyclic property of the trace to move the operator $K_3$ to the right/left
in the fourth/fifth line, respectively.

\medskip
Due to the expression (\ref{RG}) for the operator $G$, we can write
$$
%\be
G_1 I_2\, =\,\mu^{-2}\,\Tr{34}\left( (F_2^{-\varepsilon}F_3^{-\varepsilon})
K_2 F_1^{-1}F_2^{-1}(F_3^{\varepsilon} F_2^{\varepsilon})\right)\
%\ee
$$
by the formula (\ref{inv-trD}).

\smallskip
The relation (\ref{F-G}) is now proved as follows
\be
\begin{array}{l}
G_1 F_1^{\varepsilon}\, =\,\mu^{-2}\,\Tr{34}\left( (F_2^{-\varepsilon}F_3^{-\varepsilon})
K_2 F_1^{-1}F_2^{-1}(F_3^{\varepsilon} F_2^{\varepsilon})\right) F_1^{\varepsilon}\hspace{3mm}
\\[1em]
\ \ \ \ \ =\mu^{-2}\,\Tr{34} (F_2^{-\varepsilon}F_3^{-\varepsilon}) F_3^{\varepsilon}
{}F_2^{\varepsilon}F_1^{\varepsilon} K_3 F_2^{-1} F_3^{-1}\, =\, F_1^{\varepsilon}\mu^{-2}\,
\Tr{34}K_3F_2^{-1}F_3^{-1}\, =\, F_1^{\varepsilon} G_2\ .
\end{array}
\lb{FG}
\ee
Here we subsequently used the twist relations for the
pair $\{K,F^{\varepsilon}\}$, the Yang-Baxter equations for $F$ and again the expression
(\ref{RG}) for the operator $G$.

\medskip
Vanishing of the commutators $[ C_F, G]$ and $[ D_F, G]$
in eq.(\ref{comm-Gb}) follow from the above proved equality. To
find these commutators, transform eq.(\ref{FG}) to
$$
%\be
%\lb{PsiFG}
G_1 {\Psi_F}_{12} \, =\, {\Psi_F}_{12} G_2\, , \qquad
G_2 {\Psi_F}_{12} \, =\, {\Psi_F}_{12} G_1\, ,
%\ee
$$
(multiply the relation (\ref{FG}) by a combination ${\Psi_F}_{41}{\Psi_F}_{23}$ and take $\tr_{(12)}$)
and then apply the trace in the space 1 or the space 2 to these relations and compare results.

\medskip
The relation (\ref{comm-Ga}) is approved by a calculation
\ba
\nonumber
G_1 {D_R}_1 &=& \mu^{-2}\, \Tr{23} K_2 F_1^{-1} F_2^{-1} {D_R}_1\, =\,
\mu^{-2}\, \tr_{(23)} K_2 F_1^{-1} F_2^{-1} {D_R}_1{D_R}_2{D_R}_3
\\[1em]
\nonumber
&=& \mu^{-2}\, \tr_{(23)} {D_R}_1{D_R}_2{D_R}_3 K_2 F_1^{-1} F_2^{-1}\, =\, {D_R}_1 G_1\, .
\ea
Here the expression (\ref{RG}) for the operator $G$, the relations (\ref{FCC}) for $X=Y=R$ and
the relation (\ref{KDD}) were used.

\medskip%%%dobavleno
To prove the relation $[ E, G] =0$, we rewrite the expression for $G$:
\be
\begin{array}{l}
G_1=\tr_{(23)}(K_2F_1^{-1} F_2^{-1})=\eta^{-1} \tr_{(23)}(K_2K_2F_1^{-1} F_2^{-1})=
\eta^{-1} \tr_{(23)}(K_2F_1^{-1} F_2^{-1}K_1)\\[1em]
\ \ \ \ \ =\eta^{-1} \tr_{(23)}(F_1^{-1} F_2^{-1}K_1K_2)=\eta^{-1}
\tr_{(23)}(F_1^{-1} F_2^{-1}K_1P_{23}P_{12})\, E_1\ .
\end{array}
\lb{GcomE}
\ee
In the second equality we used the relation $K^2=\eta K$; in the third equality
we used the twist relation; in the fourth equality we moved the operator $K_2$
cyclically under the trace; in the fifth equality we used the first of the
relations (\ref{k12k23}).

\smallskip
Due to the relation (\ref{F-G}), the combination
$\tr_{(23)}(F_1^{-1} F_2^{-1}K_1P_{23}P_{12})$ commutes with the operator $G_1$.
Therefore the operators $G$ and $E$ commute.

\medskip
We have already shown that the operators $C_F$ and $D_F$ commute with the
operator $G$. It follows then from the expression (\ref{GcomE}) for the operator $G$
that to prove that the operators $C_F$ and $D_F$ commute with the operator $E$ it is enough
to prove that the operators $C_F$ and $D_F$ commute with the combination
$\Xi_1 :=\tr_{(23)}(F_1^{-1} F_2^{-1}K_1P_{23}P_{12})$. We have
\be
\begin{array}{l}
\Xi_1{D_{F^{-1}}}_1=\tr_{(23)}(F_1^{-1} F_2^{-1}{D_{F^{-1}}}_3
K_1P_{23}P_{12})
=\tr_{(23)}(F_1^{-1} F_2^{-1}{D_{F^{-1}}}_3 {C_R}_3{C_R}_3^{-1}
K_1P_{23}P_{12})\\[1em]
\ \ \ \ \ \ \ \ \ \ \ \ \
=\tr_{(23)}(F_1^{-1} {D_{F^{-1}}}_2 {C_R}_2F_2^{-1}K_1P_{23}P_{12})
{C_R}_1^{-1}
={D_{F^{-1}}}_1 {C_R}_1\Xi_1{C_R}_1^{-1}
={D_{F^{-1}}}_1\Xi_1\ .
%\lb{EcommDF}
\end{array}
\lb{oves}
\ee
In the first equality we moved the operator $D_{F^{-1}}$ leftwards through
the permutation operators; in the second equality we inserted
${C_R}_3{C_R}_3^{-1}$; in the third equlity we used the relations
(\ref{FCD}) and moved the operator ${C_R}_1^{-1}$ rightwards out of the
trace; in the fourth equality we used again the relations
(\ref{FCD}). The operator $C_R$ commutes with the operators $G$ and $E$
by the already proved relation (\ref{comm-Ga}) for the compatible pairs $\{ R,F\}$ and
$\{ R,P\}$; therefore, due to the expression (\ref{GcomE}) for the operator $G$,
the operator $C_R$ commutes with the operator $\Xi$, which is used in the fifth equality.

\smallskip
The calculation (\ref{oves}) establishes the relation $[ C_F,E]=0$;
the proof of the relation $[ D_F,E]=0$ is similar, we do not repeat
details. \hfill$\blacksquare$

\begin{rem}{\rm
One can rewrite further the expression (\ref{RG}) for $G$:
\ba
\nonumber
G_1&=&\mu^{-2}\,\Tr{23}F_1^{-1}F_{2}^{-1}K_1=\mu^{-2}\Tr{2}F_1^{-1}{C_F}_2{D_R}_2K_1
\\[1em]
\nonumber
%\lb{anothG}
&=&\Tr{2}F_1^{-1}{C_F}_2{D^{-1}_R}_1K_1=\mu^{-2}\Tr{2}F_1^{-1}{C_F}_2{C_R}_1K_1
\\[1em]
\nonumber
&=&\mu^{-2}{C_F}_1\Tr{2}{C_R}_2F_1^{-1}K_1={C_F}_1\tr_{(2)}F_1^{-1}K_1\ .
\ea
Here we used subsequently: the twist relation, the relations (\ref{DxCf}), (\ref{KDD}), (\ref{C*D}), (\ref{FCC}) and then again (\ref{C*D}).

\smallskip
Similarly,
$$
%\be
G^{-1}_1=\tr_{(2)}(K_1F_1){D^{-1}_F}_1\ .
%\lb{anothGb}
%\ee
$$
}\end{rem}

\subsection{Two linear maps}
\label{Twolinearmaps}

The next lemma introduces two linear maps which will be important in the study
of the matrix $\star$-product.%%%fraza

\begin{def-lem}\lb{lemma3.9}
Let $\{R,F\}$ be a compatible pair of skew invertible R-matrices, where the operator $R$ is of the
BMW type and the operator $F$ is
strict skew invertible. Define two endomorphisms $\phi$ and $\xi$ of the space ${\rm Mat}_{\mbox{\footnotesize\sc n}}(W)$:
\ba\lb{phi} \phi(M)_1 &:=& \Tr{2} \left( F_{12}M_1 F^{-1}_{12} R_{12}\right),
\qquad M\in {\rm Mat}_{\mbox{\footnotesize\sc n}}(W),\ea
and
\ba\lb{xi}
\xi(M)_1 &:=& \Tr{2} \left( F_{12}
M_1 F^{-1}_{12} K_{12}\right), \qquad M\in {\rm Mat}_{\mbox{\footnotesize\sc n}}(W)\, .\ea
The mappings $\phi$ and $\xi$ are invertible; their inverse mappings read
\ba
\lb{phi-inv}
\phi^{-1}(M)_1 &=& \mu^{-2}\TR{2}{R_f} \left(  F_{12}^{-1}
M_1 R^{-1}_{12} F_{12}\right)
\ea
and
\ba\lb{xi-inv}\xi^{-1}(M)_1 &=& \mu^{-2} \TR{2}{R_f}\left( F^{-1}_{12} M_1
K_{12} F_{12}\right)\, .\ea
{}The following relations for the R-traces
\be\lb{phi-xi-traces}
{\rm Tr}_{\! R_f^{\rule{0pt}{6pt}}}\phi(M)\, =\, {\rm Tr}_{\! R^{\rule{0pt}{6pt}}} M\, , \qquad
{\rm Tr}_{\! R_f^{\rule{0pt}{6pt}}}\xi (M)\, =\, \mu\, {\rm Tr}_{\! R^{\rule{0pt}{6pt}}} M\, .\ee
are satisfied.
\end{def-lem}

\nin {\bf Proof.~} The expressions in the right hand sides of the formulas (\ref{phi-inv}) and
(\ref{xi-inv}) are well defined, since, by  proposition \ref{proposition3.6} b), the R-matrix
$R_f$ is skew invertible.

\smallskip
Let us check the relation $\phi^{-1}(\phi(M))=M$ directly.

\smallskip
Using the formulas (\ref{phi}) and (\ref{phi-inv}) and applying the relation (\ref{inv-trD}) for the
pair $\{R,F\}$ we begin a calculation
\ba\nn {\phi^{-1}(\phi(M))}_1 &=& \mu^{-2}\, \TR{2}{R_f}\left( F_{12}^{-1}
(\Tr{2'}  F_{12'} M_1 F^{-1}_{12'} R_{12'}) R^{-1}_{12} F_{12}\right)\\[1em]\nn &=&
\mu^{-2}\, \TR{2}{R_f}\Tr{3}\left( F_1^{-1} F_2^{-1} \underline{F_1} M_1
{}F_1^{-1} R_1 F_2 R_1^{-1} F_1\right)\ .\ea
In the next step we move the element $F_1$, underlined in the expression above, to the left and
it becomes $F_2$ due to the Yang-Baxter equation; then we
transport the operator to the right using the cyclic property of the trace
(when $F_2$ moves cyclically,
$\TR{2}{R_f}\Tr{3}$ becomes $\Tr{2}\TR{3}{R_f}$ due to the relations (\ref{FCC})$\,$).
Applying the Yang-Baxter equation for the operator $F$ and the relations (\ref{FCC}) in the case
$X=R$ and $Y=R_f$, we continue the calculation
\be\begin{array}{ccl} {\phi^{-1}(\phi(M))}_1&=& \mu^{-2}\, \Tr{2}
\TR{3}{R_f}\!\left(  F^{-1}_1 F^{-1}_2 M_1
\underline{F_1^{-1} R_1} F_2 \underline{R_1^{-1} F_1}F_2\right)\\[1em] &=&\mu^{-2}\, \Tr{2}\TR{3}{R_f}\!
\left(  F^{-1}_1M_1F_2^{-1} {R_f}_1\,
\underline{F_1^{-1} F_2 F_1} {R_f}^{-1}_1\, F_2\right)\\[1em] &=&
\mu^{-2}\, \Tr{2}\TR{3}{R_f}\!\left(  F^{-1}_1M_1
\underline{F_2^{-1}{R_f}_1F_2}\, F_1\,\underline{F_2^{-1}
{R_f}^{-1}_1\, F_2}\right)\\[1em] &=&\mu^{-2}\, \Tr{2}\!\left(  F^{-1}_1 M_1 F_1
\underline{(\TR{3}{R_f} {R_f}_2 F_1 {R_f}_2^{-1}\, )} F_1^{-1}\right)\ .\end{array}\lb{step}\ee
Here we consequently transformed the underlined expressions using the definition of the twisted
R-matrix $R_f$, the Yang-Baxter equation for the operator $F$ and the twist relations for the
compatible pair $\{R_f,F\}$. To calculate the trace
underlined in the last line of eq.(\ref{step}), we apply the relation (\ref{inv-trD}) for the
compatible pair $\{R_f, R_f\}$ and then use the relation (\ref{DxCf}) written for the
compatible pair $\{R_f,F^{-1}\}$. The result reads
\ba\nn {\phi^{-1}(\phi(M))}_1&=& \mu^{-2}\, \Tr{2}\!\left(  F^{-1}_1 M_1 F_1
\underline{(D_{R_f} C_{F^{-1}})_1 F^{-1}_1}\right)\ .\ea
Now, using the relations (\ref{FCD}), written for the compatible pairs $\{R_f, F\}$ and
$\{F^{-1},F\}$, the relations (\ref{CDtwist}) and (\ref{CDinv}) for $X=F$, the relations
(\ref{C*D}) and the (\ref{traceCD-X})
for $X=F^{-1}$, we complete the calculation
\ba
\nn
&&\hspace{-22mm}{\phi^{-1}(\phi(M))}_1\, =\, \mu^{-2}\, \tr_{(2)}\!\!\left((D_{R_f}
C_{F^{-1}} D_R)_2 F^{-1}_1\right)\! M_1
\\[1em]
\nn
&&=\, \mu^{-2}\, \tr_{(2)}\!\!\left((D_{F^{-1}}
C_R D_F C_{F^{-1}} D_R)_2 F^{-1}_1\right)\! M_1
\, =\, \tr_{(2)}({D_{F^{-1}}}_2 F^{-1}_{12})M_1 = M_1\, .
\ea

\smallskip
A proof of the equality $\xi^{-1}(\xi(M))=M$ proceeds quite similarly until
the line (\ref{step}), where one has to use a relation
$$
%\be
%\lb{traceKMK}
\Tr{2}(K_1 M_1 K_1) =  (\tr_{\!\! R} M) I_1\, \quad \forall\;
M\in {\rm Mat}_{\mbox{\footnotesize\sc n}}(W)\,
%\ee
$$
instead of the relation (\ref{inv-trD}). This in turn follows from the relations (\ref{traceK})
and (\ref{traceDK}) and the property ${\rm rk}\, K=1$.

\smallskip
The relations (\ref{phi-xi-traces}) can be directly checked starting from the definitions (\ref{phi}) and
(\ref{xi}), applying the relation (\ref{FCC}) in the case $X=R$ and $Y=R_f$ and then using the formulas (\ref{traceR})
and (\ref{traceDK}). \hfill$\blacksquare$

\begin{rem}\lb{remark3.10}
{\rm For the mapping $\phi$, the statement of lemma \ref{lemma3.9} remains valid if one weakens
the conditions, imposed on the R-matrix $R$, replacing the BMW type condition by the strict skew
invertibility. In this case, one should substitute the term $\mu^{-2} D_{R_f}$ by $D_{R_f^{-1}}$ in the
expression (\ref{phi-inv}) for the inverse mapping $\phi^{-1}$. The proof repeats
the proof of the formula (\ref{phi-inv}).}\end{rem}

%\newpage
\section{Quantum matrix algebra}\lb{sec4}

In this section we deal with the main objects of our study, the quantum matrix algebras,
and construct the $\star$-product for them. We mainly discuss the quantum matrix algebras of
the type BMW.

\medskip
In subsection \ref{subsec4.2} we introduce a {\em characteristic} subalgebra of
the quantum matrix algebra. In the theory of the polynomial identities, a ring, generated by
the traces of products of generic matrices, is known as the ring of matrix invariants (see,
{\em e.g.}, \cite{F}). The characteristic subalgebra can be understood as a generalization of the ring of matrix invariants (in the simplest case of a single matrix) to the setting of the
quantum matrix algebras and, simultaneously,
to a situation when the invariants can be formed not only by taking a trace (on the quantum
level, the invariants can be conveniently formed by taking the R-trace of a product of a
`string' $M_{\overline{1}}M_{\overline{2}}\dots M_{\overline{n}}$ by a matrix image of a
word in the braid group ${\cal B}_n$).

\smallskip
In propositions \ref{proposition4.7}, \ref{proposition4.8} we exhibit three generating sets of the characteristic subalgebra in the BMW case. Explicit relations between the generators of these sets will be constructed in  section \ref{sec5}. Some preparatory work for this constructions is performed in the rest of  section \ref{sec4}.

\medskip
In subsection \ref{subsec4.4} we introduce an algebra ${\cal P}(R,F)$ for the
quantum matrix algebras of the general type. The algebra ${\cal P}(R,F)$ has the same
relationship to the characteristic subalgebra as the trace ring (see, {\em e.g.}, \cite{F})
to the ring of matrix invariants.

\smallskip
In subsection \ref{subsec4.4a} we prove the commutativity of the algebra ${\cal P}(R,F)$
in the case of the quantum matrix algebras of the BMW type.

\smallskip
In subsection \ref{subsec4.5} we define an extended quantum matrix algebra of the BMW type
by adding an inverse of the quantum matrix.

\subsection{Definition}\lb{subsec4.1}

Consider a linear space ${\rm Mat}_{\mbox{\footnotesize\sc n}}(W)$, introduced in the definition
\ref{definition3.1}. For a fixed element $F\in {\rm Aut}(V\otimes V)$, we consider series of `copies'
$M_{\overline{i}}$, $i=1,2,\dots ,n,$ of a matrix $M\in {\rm Mat}_{\mbox{\footnotesize\sc n}}(W)$.
They are defined recursively by
\be
M_{\overline 1}:=M_1, \quad M_{\overline{i}}:=
{}F^{\phantom{-1}}_{i-1}M_{\overline{i-1}}F_{i-1}^{-1}\ .
\lb{kopii}
\ee
{}For $F=P$, these are usual copies, $M_{\overline{i}}=M_i$, but, in general,
$M_{\overline{i}}$ can be nontrivial in all the spaces $1,\dots ,i$.

\smallskip
We shall, slightly abusing notation, denote by the same symbol $M_{\overline{i}}$
an element in ${\rm Mat}_{\mbox{\footnotesize\sc n}}(W)^{\otimes k}$ for any $k\geq i$, which is
defined by an inclusion of the spaces
$$
%\be
%\lb{kopii2}
{\rm Mat}_{\mbox{\footnotesize\sc n}}(W)^{\otimes j}\hookrightarrow
{\rm Mat}_{\mbox{\footnotesize\sc n}}(W)^{\otimes (j+1)}:\quad
 {\rm Mat}_{\mbox{\footnotesize\sc n}}(W)^{\otimes j}\ni X\mapsto
X\otimes I\in{\rm Mat}_{\mbox{\footnotesize\sc n}}(W)^{\otimes (j+1)}\, .
%\ee
$$

From now on we specify $W$ to be
the associative $\Bbb C$-algebra freely generated by the unity and by $\mbox{\sc n}^2$ elements $M_a^b$,~~ $W := {\Bbb C}\langle 1,M_a^b\rangle$, $1\le a,b\le \mbox{\sc n}$.

\begin{defin}\lb{definition4.1}
Let $\{R,F\}$ be a compatible pair of strict skew invertible R-matrices (see  section \ref{subsec3.1}).
A quantum matrix algebra ${\cal M}(R,F)$ is a quotient algebra of the algebra $W={\Bbb C}\langle 1,M_a^b\rangle$ by a
two-sided ideal generated by entries of the matrix relation
\be R_1M_{\overline 1}M_{\overline 2} = M_{\overline 1}M_{\overline 2}R_1\, ,\label{qma}\ee
where $M = \|M_a^b\|_{a,b=1}^{\mbox{\footnotesize\sc n}}$ is a matrix of the generators
of ${\cal M}(R,F)$ and the matrix copies $M_{\overline i}$ are constructed with the help
of the R-matrix $F$ as in eq.(\ref{kopii}).

\medskip
If $R$ is an R-matrix of the BMW type  (see eqs.(\ref{charR})--(\ref{bmwRb})$\,$) then ${\cal M}(R,F)$ is
called a BMW type quantum matrix algebra.

%\smallskip
%If $R$ is a $Sp(2k)$- , respectively, an $O(k)$-type R-matrix (see the definition \ref{definition3.11}) then
%${\cal M}(R,F)$ is called a $Sp(2k)$- , respectively, an $O(k)$-type quantum matrix algebra.
\end{defin}

\begin{rem}\lb{remark4.2}
{\rm The quantum matrix algebras were introduced in Ref. \cite{Hl} under the name `quantized
braided groups'. In the context of the present paper they have been first investigated
in \cite{IOP1}. The matrix $M'$ of the generators of the algebra ${\cal M}(R,F)$ used in \cite{IOP1}
is different from the matrix $M$ that we use here. A relation between these two matrices is explained
in  section 3 of \cite{IOP2}: $M' =D_R M (D_F)^{-1}$.}\end{rem}

\begin{lem}{\rm \bf \cite{IOP1}} \lb{lemma4.3}
The matrix copies of the matrix $M= \|M_a^b\|_{a,b=1}^{\mbox{\footnotesize\sc n}}$ of the generators of
the algebra ${\cal M}(R,F)$ satisfy relations
\ba
F_i\, M_{\overline j}\, &=&\, M_{\overline j}\, F_i\,\;\;\quad\qquad \mbox{for}\;\; j\neq i,i+1,
\lb{fm-k}
\\[1em]
R_i\, M_{\overline j}\, &=&\,M_{\overline j}\, R_i\,\;\;\quad\qquad \mbox{for}\;\;
j\neq i,i+1,\lb{rm-k}\\[1em] R_j\,M_{\overline j}\,M_{\overline {j+1}}\, &=&\,  M_{\overline j}\,
M_{\overline {j+1}}\, R_j\,\quad \mbox{for}\;\; j=1,2,\dots \ ,
\lb{rmm-k}
\\[1em]
{}F_i\, F_{i+1}\dots F_k\cdot M_{\overline i}\, M_{\overline {i+1}}\dots
M_{\overline k}\, &=& \, M_{\overline {i+1}}\, M_{\overline {i+2}}\dots
M_{\overline{k+1}}\cdot F_i\, F_{i+1}\dots F_k\,\quad \mbox{for}\;\; i\leq k\ .
\lb{stringshift}
\ea
\end{lem}

\subsection{Characteristic subalgebra}\lb{subsec4.2}

{}From now on we assume that $M$ is the matrix of generators of the quantum matrix algebra
${\cal M}(R,F)$ and its copies $M_{\overline n}$ are calculated by the rule (\ref{kopii}).

\medskip
Denote by ${\cal C}(R,F)$ a vector subspace of the quantum matrix algebra ${\cal M}(R,F)$ linearly
spanned by the unity and elements
\be
\lb{char}
ch(\alpha^{(n)}) := \Tr{1,\dots ,n}(M_{\overline 1}\dots M_{\overline n}\,
\rho_R(\alpha^{(n)}))\ ,\quad n =1,2,\dots\ ,
\ee
where $\alpha^{(n)}$ is an arbitrary element
of the braid group ${\cal B}_n$.

\smallskip
Notice that elements of the space ${\cal C}(R,F)$ satisfy a {\em cyclic property}
\be
\lb{cyclic}
ch(\alpha^{(n)}\beta^{(n)}) = ch(\beta^{(n)}\alpha^{(n)})\,  \quad
\forall\; \alpha^{(n)}, \beta^{(n)}\in
%{\Bbb C}
{\cal B}_{n}\, ,\quad n=1,2,\dots\ ,
\ee
which is a direct consequence of the relations (\ref{rm-k}), (\ref{rmm-k}) and (\ref{RDD}) and the cyclic property of
the trace.

\begin{def-prop} {\rm\bf \cite{IOP1}}\lb{proposition4.4}
The space ${\cal C}(R,F)$ is a commutative subalgebra of the quantum matrix algebra ${\cal M}(R,F)$:
\be
\lb{multip-rule}
ch(\alpha^{(n)})\, ch(\beta^{(i)}) =ch(\alpha^{(n)}\, \beta^{(i)\uparrow n}) =
ch(\alpha^{(n)\uparrow i}\, \beta^{(i)})\, .
\ee
Recall that $\alpha^{(n)\uparrow i}$ denotes the image of an element $\alpha^{(n)}$ under the embedding
${\cal B}_n \hookrightarrow {\cal B}_{n+i}$ defined in (\ref{h-emb2}).
We shall call ${\cal C}(R,F)$ the characteristic subalgebra of ${\cal M}(R,F)$.\end{def-prop}

A proof of the proposition given in \cite{IOP1} is based in particular on the following lemma:

\begin{lem}{\rm\bf \cite{IOP1}}\lb{lemma4.5}
Consider an arbitrary element $\alpha^{(n)}$ of the braid group ${\cal B}_n$. Let $\{R,F\}$
be a compatible pair of R-matrices, where $R$ is skew invertible. Then relations
\be\lb{char1}\Tr{i+1,\dots ,i+n}(M_{\overline{i+1}}\dots
M_{\overline{i+n}}\ \rho_R(\alpha^{(n)\uparrow i}))\, =\, I_{1,2,\dots ,i}\  ch(\alpha^{(n)})\,\ee
hold for any matrix $M\in {\rm Mat}_{\mbox{\footnotesize\sc n}}(W)$\footnote{
Here there is no need to specify $M$ to be the matrix of the generators of the algebra ${\cal M}(R,F)$.}.\end{lem}

We will make use of lemma \ref{lemma4.5} several times below.
\medskip

Let us introduce a shorthand notation for certain elements of ${\cal C}(R,F)$
\ba\lb{P-01} p_0& :=& \tr_{\!\! R}\, I\; (=\mu \eta \mbox{~in the BMW case})\,  ,\qquad
p_1 := \tr_{\!\! R}\, M\, ,\\[1em]\lb{P_k} p_i& := & ch(\sigma_{i-1}\dots\sigma_2\sigma_1) =
ch(\sigma_1\sigma_2\dots\sigma_{i-1})\, ,  \quad i=2,3,\dots .\ea
The last equality in eq.(\ref{P_k}) is due to the inner automorphism (\ref{innalis}) and the
cyclic property (\ref{cyclic})$\,$.

The elements $p_i$ are called {\em traces of powers of $M$} or, shortly, {\em power sums}.

\medskip
From now on in this subsection we assume the R-matrix $R$ and, hence, the algebra  ${\cal M}(R,F)$
to be of the BMW type. Denote
\be\lb{tau}\textstyle g\, :=\, ch(c^{(2)}) \,\equiv\, \eta^{-1} ch(\kappa_1)\, \equiv\,  \eta^{-1}\,
\Tr{1,2} \left( M_{\overline{1}}M_{\overline{2}}\, K_1\right)\, .\ee
The notation used here was introduced in the formulas (\ref{kappa}), (\ref{bmw5b}),
(\ref{kappa-i}) and (\ref{K}). We call the element $g$  a {\em contraction of two matrices} $M$
or, simply, a {\em 2-contraction}.

\begin{lem}\lb{lemma4.6}
Let $M$ be the matrix of generators of the BMW type quantum matrix algebra ${\cal M}(R,F)$. Then
its copies, defined in eq.(\ref{kopii}), fulfill relations
\be\lb{tau2}K_{n}\, M_{\overline{n}}M_{\overline{n+1}}\, =\,
M_{\overline{n}}M_{\overline{n+1}}\, K_{n}\, =\,\mu^{-2} K_{n}\, g\, \quad
\forall\; n\geq 1\, .\ee\end{lem}

\nin {\bf Proof.~} We employ induction on $n$. Due to the property ${\rm rk}\, K=1$, one has
$$
%\be
%\lb{tau3}
K_1\, M_{\overline{1}}M_{\overline{2}}\, =\, M_{\overline{1}}M_{\overline{2}}\, K_1\, =\,
K_1\, t\, ,
%\ee
$$
where $t\in {\cal M}(R,F)$ is a scalar. Evaluating the R-trace of this equality in the spaces 1 and 2
and using the relations (\ref{traceDK}) and (\ref{traceD}), one finds  $t=\mu^{-2} g$, which
proves the relation (\ref{tau2}) in the case $i=1$. It remains to check the induction step
$n\rightarrow (n+1)$:
$$\begin{array}{ccl} K_{n+1}M_{\overline{n+1}}M_{\overline{n+2}}&=&
K_{n+1}(F_n M_{\overline{n}}
{}F^{-1}_{n})M_{\overline{n+2}}\, =\, K_{n+1} F_n M_{\overline{n}} (F_{n+1}M_{\overline{n+1}}F^{-1}_{n+1})
{}F^{-1}_{n} \\[1em] &=&(K_{n+1}F_n F_{n+1}) M_{\overline{n}} M_{\overline{n+1}}F^{-1}_{n+1}F^{-1}_{n}\, =\,
{}F_n F_{n+1} (K_{n}M_{\overline{n}} M_{\overline{n+1}})F^{-1}_{n+1}F^{-1}_{n}\\[1em] &=& \mu^{-2} F_n
{}F_{n+1} K_n F^{-1}_{n+1}F^{-1}_{n} g\, =\, \mu^{-2} K_{n+1}\, g\, .\end{array}$$
Here eqs.(\ref{kopii}) and (\ref{fm-k}), the twist relation (\ref{sovm}) for the pair $\{K,F\}$
and the induction assumption were used for the transformation. \hfill$\blacksquare$

\begin{prop}\lb{proposition4.7}
Let ${\cal M}(R,F)$ be the quantum matrix algebra of the BMW type. Its characteristic subalgebra
${\cal C}(R,F)$ is generated by the set $\{g,p_i\}_{i\geq 0}$.
\end{prop}

\nin {\bf Proof.~} Consider the chain of the BMW algebras monomorphisms (\ref{h-emb})--(\ref{h-emb2}).
We adapt, for $n\geq 3$, the  following presentation for  an element $\alpha^{(n)}\in {\cal W}_n$
\be
\lb{char8}
\alpha^{(n)} = \beta\sigma_1 \beta' + \gamma \kappa_1 \gamma' + \delta\, ,
\ee
where $\beta,\beta',\gamma,\gamma',\delta\in {\rm Im}({\cal W}_{n-1})\subset {\cal W}_n$. For $n=3$,
the formula (\ref{char8}) follows from the relations (\ref{braid})--(\ref{bmw7}). For  $n>3$, it can be proved
by induction on $n$ (one has to prove that the expressions of the form
(\ref{char8}) form an algebra, for which it is enough to show that
the products $\sigma_1\beta\sigma_1$, $\sigma_1\beta\kappa_1$,
$\kappa_1\beta\sigma_1$ and $\kappa_1\beta\kappa_1$ with $\beta\in
{\rm Im}({\cal W}_{n-1})\subset {\cal W}_n$ can be rewritten in the
form (\ref{char8}); this is done by further decomposing $\beta$,
using the induction assumption,
$\beta = \tilde{\beta}\sigma_2\tilde{\beta}'+\tilde{\gamma}
\kappa_2\tilde{\gamma}' +\tilde{\delta}$, where $\tilde{\beta},\tilde{\beta}',\tilde{\gamma},\tilde{\gamma}',
\tilde{\delta}\in {\rm Im}({\cal W}_{n-2})\subset {\cal W}_n$).

\medskip
Using the expression (\ref{char8}) for $\alpha^{(n)}$ and the cyclic property (\ref{cyclic}),
we conclude that, in the BMW case, any element (\ref{char}) of the characteristic subalgebra
can be expressed as a linear combination of terms
\be
\lb{char8a}
ch(\alpha_1\alpha_2\dots\alpha_{n-1})\, , \quad
\mbox{where}\quad \alpha_i\in\{1,\sigma_i,\kappa_i\}\, .
\ee

Let us analyze the expressions (\ref{char8a}) for different choices of $\alpha_i$.

\medskip\noindent
i)~ If $\alpha_i=1$ for some value of $i$,  then, applying the relation (\ref{char1}), we get
\be
\lb{char8b}
ch(\alpha_1\dots\alpha_{i-1}\alpha_{i+1}\dots\alpha_{n-1}) =
ch(\alpha_1\dots\alpha_{i-1})\, ch((\alpha_{i+1}\dots\alpha_{n-1})^{\downarrow i})\, ,
\ee
where $(\alpha_{i+1}\dots\alpha_{n-1})^{\downarrow i}\in {\cal W}_{n-i}$ is the preimage
of $(\alpha_{i+1}\dots\alpha_{n-1})\in {\cal W}_n$.

\medskip\noindent
ii)~ In the case when $\alpha_{n-1}=\kappa_{n-1}$, we apply the relation (\ref{tau2}) and then
the relations (\ref{traceR}), (\ref{traceDK}) or (\ref{traceD}) to reduce the expression (\ref{char8a}) to
\be
\lb{char8c}
ch(\alpha_1\dots\alpha_{n-2}\kappa_{n-1}) =f(\alpha_{n-2})\,
ch(\alpha_1\dots\alpha_{n-3})\, g\, ,
\ee
where $f(\sigma_{n-2})=\mu^{-1}$, $f(\kappa_{n-2})=1$ and $f(1)=\eta$.

\medskip\noindent
iii)~ In the case when $\alpha_i=\kappa_i$ for some $i$, and $\alpha_j=\sigma_j$ for all $j=i+1,\dots ,n-1$,
we perform the following transformations
\be
\begin{array}{ccl}
ch(\alpha_1\dots\alpha_{i-1}\underline{\kappa_i\sigma_{i+1}}\sigma_{i+2}\dots \sigma_{n-1}) =
ch(\alpha_1\dots\alpha_{i-2}\,\sigma_i^{-1}\alpha_{i-1}\kappa_i\underline{\kappa_{i+1}\sigma_{i+2}}\dots
\sigma_{n-1}) &&\\[1em] =\dots =ch(\alpha_1\dots\alpha_{i-2} (\sigma_{n-2}^{-1}\dots \sigma_i^{-1})
\alpha_{i-1}\kappa_i\kappa_{i+1}\dots\kappa_{n-1}) .\hspace{35mm}&&
\end{array}
\lb{char10}
\ee
Here the relations (\ref{bmw3}) and the cyclic property (\ref{cyclic}) are repeatedly used; expressions
suffering a transformation are underlined.

\smallskip
Now, depending on a value of $\alpha_{i-1}$, we proceed in different ways.

If $\alpha_{i-1}=\kappa_{i-1}$ then by  eqs.(\ref{bmw3}) and (\ref{char8c}) we have
$$
%\be
\begin{array}{ccl}
(\ref{char10}) &=& ch(\alpha_1\dots\alpha_{i-2}\,
\sigma_{i-1} \sigma_i\dots \sigma_{n-3} \kappa_{n-2}\kappa_{n-1})
\\[1em]
&=&
ch(\alpha_1\dots\alpha_{i-2}\,\sigma_{i-1} \sigma_i\dots \sigma_{n-3})\, g\, .
\end{array}
%\lb{char11}
%\ee
$$
If $\alpha_{i-1}=\sigma_{i-1} = \sigma_{i-1}^{-1}+(q-q^{-1})(1-\kappa_{i-1})$ then, using
the relations $\sigma_i^{-1}\sigma_{i-1}^{-1}\kappa_i=\kappa_{i-1}\kappa_i$ and applying
the previous results (\ref{char8c}) and (\ref{char8b}), we obtain
$$
%\be
\begin{array}{ccl}
(\ref{char10}) &=&ch(\alpha_1\dots\alpha_{i-2}\, \kappa_{i-1}\sigma_i\dots \sigma_{n-3})\, g
\\[1em]\nonumber &&+\, (q-q^{-1})\, \mu^{-1}\, ch(\alpha_1\dots\alpha_{i-2})\, p_{n-i-1}\, g\\[1em]
&&-\, (q-q^{-1})\, ch(\alpha_1\dots\alpha_{i-2}\,
\sigma_{i-1} \sigma_i\dots \sigma_{n-3})\, g\, .
\end{array}
%\lb{char12}
%\ee
$$
Iterating transformations i)---iii)  finitely many times, we eventually prove the assertion
of the proposition. \hfill$\blacksquare$

\medskip
We keep considering the BMW type quantum matrix algebra ${\cal M}(R,F)$ with the R-matrix $R$ generating
representations of the algebras ${\cal W}_n(q,\mu)$, $n=1,2,\dots$. Assume that the antisymmetrizers $a^{(i)}$
and symmetrizers $s^{(i)}$ in these latter algebras are consistently defined (see eqs.(\ref{a^k}),
(\ref{s^k}) and (\ref{mu})$\,$). In this case, we can introduce two following sets of elements in the characteristic
subalgebra ${\cal C}(R,F)$
\ba\lb{SA_0}a_0& :=&1\ \ \ {\mathrm{and}}\ \ \ s_0\, :=\, 1\, ;\\[1em]\lb{SA_k} a_i &:=& ch(a^{(i)})\ \ \
{\mathrm{and}}\ \ \ s_i\, :=\,
ch(s^{(i)})\, ,\quad i=1,2,\dots  .\ea

\begin{prop}\lb{proposition4.8}
Let ${\cal M}(R,F)$ be the quantum matrix algebra of the BMW type. Assume that $j_q\neq 0 ,\;\; \mu\neq -q^{-2j+3}$ (respectively,
$j_q\neq 0,$ $\mu\neq q^{2j-3}$) for all $j=2,3,\dots\; .$ Then the characteristic subalgebra
${\cal C}(R,F)$ is generated by the set $\{g,a_i\}_{i\geq 0}$ (respectively, $\{g,s_i\}_{i\geq 0}$).

%\smallskip
%Let the quantum matrix algebra ${\cal M}(R,F)$ be of the types either  $Sp(2k)$ or $O(k)$ (this implies restrictions on $q$:
%$j_q\neq 0$ for all $j=2,3,\dots ,k$). Then the characteristic subalgebra ${\cal C}(R,F)$ is
%generated by the set $\{g,a_i\}_{i=0}^k$.
\end{prop}

\nin {\bf Proof.~} These statements are byproducts of the previous proposition and the Newton relations,
which are proved in  section \ref{sec5}, theorem \ref{theorem6.1}. \hfill$\blacksquare$

\subsection{Matrix $\star\,$-product, general case}\label{subsec4.4}

Consider the quantum matrix algebra ${\cal M}(R,F)$ of the general type (no additional conditions on an R-matrix $R$).

\smallskip
Denote by ${\cal P}(R,F)$ a linear subspace of ${\rm Mat}_{\mbox{\footnotesize\sc n}}({\cal M}(R,F))$
spanned by ${\cal C}(R,F)$-multiples of the identity matrix, $I\, ch$ $\forall\, ch\in {\cal C}(R,F)$,
and by elements
\be
\lb{pow}
M^1 := M , \qquad (M^{\alpha^{(n)}})_{1} := \Tr{2,\dots ,n}(M_{\overline 1}
\dots M_{\overline n}\,\rho_R(\alpha^{(n)}))\ ,\quad n =2,3,\dots\ ,
\ee
where $\alpha^{(n)}$ belongs to
%the group algebra of
the braid group ${\cal B}_n$. The space ${\cal P}(R,F)$
inherits a structure of a right ${\cal C}(R,F)$--module
\be
\lb{r-module}
M^{\alpha^{(n)}} ch(\beta^{(i)}) =M^{(\alpha^{(n)}\beta^{(i)\uparrow n})}\, \quad
\forall\, \alpha^{(n)}\in {\cal B}_n, \; \beta^{(i)}\in {\cal B}_i\, ,\quad n,i=1,2,\dots\, ,
\ee
which is just a component-wise multiplication of the matrix $M^{\alpha^{(n)}}$ by the element
$ch(\beta^{(i)})$ (use the relation (\ref{char1}) to check this). The ${\cal C}(R,F)$--module
structure agrees with an R-trace map
$\tr_{\!\! R}$ (which means that $\tr_{\!\! R} (Xa)=\tr_{\!\! R} (X)a\ \ \forall\ X\in{\cal P}(R,F)$
and $\forall\ a\in {\cal C}(R,F)$)
\be
\lb{Rtrace-map}
{\cal P}(R,F)\stackrel{\tr_{\!\! R}}{\longrightarrow}
{\cal C}(R,F)\, : \quad\left\{ \begin{array}{rcl} M^{\alpha^{(n)}}&\mapsto& ch(\alpha^{(n)})\, ,
\\[1em] I\, ch(\alpha^{(n)})&\mapsto& (\tr_{\!\! R} I)\, ch(\alpha^{(n)})\, , \end{array}\right.
\ee
where $\alpha^{(n)}\in {\cal B}_n\, , \;\; n=1,2,\dots$

\smallskip
Besides, elements of the space ${\cal P}(R,F)$ satisfy a {\em reduced cyclic property}
\be\lb{red-cycl} M^{(\alpha^{(n)}\beta^{(n-1)\uparrow 1})} =
M^{(\beta^{(n-1)\uparrow 1}\alpha^{(n)})}\, \quad\forall\, \alpha^{(n)}\in {\cal B}_n , \;
\beta^{(n-1)}\in {\cal B}_{n-1} ,\quad n=2,3,\dots \, .\ee

\begin{def-prop}\lb{proposition4.12}
{}Formulas
\be
\lb{MaMb}
M^{\alpha^{(n)}}\! \star  M^{\beta^{(i)}} \,
:=\,  M^{(\alpha^{(n)}\star \beta^{(i)})}\, ,
\ee
where
\ba
\lb{a*b}
\qquad\alpha^{(n)}\star \beta^{(i)} &:=&
\alpha^{(n)}\beta^{(i)\uparrow n} (\sigma_n\dots \sigma_2 \sigma_1\sigma_2^{-1}\dots \sigma_n^{-1})\, ,
\\[1em]
\lb{MaI}
(I\, ch(\beta^{(i)})) \star  M^{\alpha^{(n)}} &:=&
M^{\alpha^{(n)}}\! \star  (I\, ch(\beta^{(i)}))\, :=\, M^{\alpha^{(n)}}  ch(\beta^{(i)})\, ,
\\[1em]
\lb{II}
(I\, ch(\alpha^{(i)}))\star (I\, ch(\beta^{(n)}))& :=& I\,( ch(\alpha^{(i)})\, ch(\beta^{(n)}))\ ,
\ea
define an associative multiplication on the space ${\cal P}(R,F)$, which agrees with the ${\cal C}(R,F)$--module
structure (\ref{r-module}).\footnote{In other words,
a map
~$ch(\alpha^{(n)})\mapsto I\, ch(\alpha^{(n)})$~ is an algebra monomorphism
~${\cal C}(R,F)\hookrightarrow {\cal P}(R,F)$.}\end{def-prop}

\nin {\bf Proof.~} To prove the associativity of the multiplication (\ref{MaMb}), it is enough to check
$$(\alpha^{(n)} \star  \beta^{(i)}) \star  \gamma^{(m)}\ =\alpha^{(n)} \star  (\beta^{(i)} \star
\gamma^{(m)})\, ,$$
which is a staightforward exercise in an application of the relations (\ref{braid}) and
(\ref{braid2}).

\smallskip
It is less trivial to prove a compatibility condition for the formulas (\ref{MaMb}) and
(\ref{MaI})
$$\left\{M^{\alpha^{(n)}}\! \star  (I\, ch(\beta^{(i)}))\right\} \star  M^{\gamma^{(m)}}\, =\, M^{\alpha^{(n)}} \star
\left\{(I\, ch(\beta^{(i)}))\star  M^{\gamma^{(m)}}\right\}\, ,$$
which, in terms of the matrix `exponents', amounts to
\be\begin{array}{l} \alpha^{(n)}\beta^{(i)\uparrow n} \gamma^{(m)\uparrow (i+n)}
(\sigma_{i+n}\dots \sigma_2\sigma_1\sigma_2^{-1}\dots \sigma_{i+n}^{-1})\hspace{50mm}\\[1em]
\ \ \ \ \ \ \ \ \ \ \ \ \ \ \ \stackrel{\rm mod\, (\ref{red-cycl})}{=}\
\alpha^{(n)}\gamma^{(m)\uparrow n} \beta^{(i)\uparrow (m+n)}
(\sigma_n\dots \sigma_2\sigma_1\sigma_2^{-1}\dots \sigma_n^{-1})\, .\end{array}\lb{lesstriv}\ee
Here the symbol $\stackrel{\rm mod\, (\ref{red-cycl})}{=}$ means the equality modulo the reduced cyclic property
(\ref{red-cycl}).

\smallskip
To check eq.(\ref{lesstriv}), we apply a technique, which was used in \cite{IOP1} to prove
the commutativity of the characteristic subalgebra. Consider an element
\be
\lb{u}
\begin{array}{rcl}
u_{i,m}^{(i+m)} &:=& (\sigma_i\dots \sigma_2\sigma_1)
(\sigma_{i+1}\dots \sigma_3\sigma_2) \dots (\sigma_{i+m-1}\dots \sigma_{m+1}\sigma_m)\\[1em]
&=&(\sigma_i\sigma_{i+1}\dots\sigma_{i+m-1})
(\sigma_{i-1}\sigma_i\dots \sigma_{i+m-2}) \dots (\sigma_{1}\sigma_2\dots\sigma_m)\, ,
\end{array}
\ee
which intertwines certain elements of the braid group ${\cal B}_{(i+m)}$:
\be
\lb{u1}
\beta^{(i)}\, u_{i,m}^{(i+m)}\, =\, u_{i,m}^{(i+m)}\, \beta^{(i)\uparrow m}\, ,\qquad
\gamma^{(m)\uparrow i}\, u_{i,m}^{(i+m)}\, =\, u_{i,m}^{(i+m)}\, \gamma^{(m)}\, .
\ee
Substitute an expression $(u_{i,m}^{(i+m)\uparrow n}\gamma^{(m)\uparrow n}\beta^{(i)\uparrow (n+m)}
(u_{i,m}^{(i+m)\uparrow n})^{-1})$ for the factor $(\beta^{(i)\uparrow n} \gamma^{(m)\uparrow (i+n)})$
in the left hand side of the equation (\ref{lesstriv}), move the element $u_{i,m}^{(i+m)\uparrow n}$
cyclically to the right and then use an equality
\be
\lb{u2}
(\sigma^{-1}_1\sigma^{-1}_2\dots \sigma^{-1}_i) u_{i,m}^{(i+m)}
= u_{i,m-1}^{(i+m-1)\uparrow 1}
\ee
to cancel it on the right hand side. Such transformation results in the right hand side of
the equation (\ref{lesstriv}).

\smallskip
Consistency of the multiplication and the ${\cal C}(R,F)$--module structures on ${\cal P}(R,F)$
follows  obviously from the last equality in
(\ref{MaI}). \hfill$\blacksquare$

\medskip
To illustrate the relation between the $\star \, $-product and the usual matrix multiplication,
we present formulas (\ref{MaMb}) and (\ref{a*b}) in the case $n=1$ ($\alpha^{(1)}\equiv 1$) in a form
\be
\lb{M*}
M \star  N = M\cdot \phi(N)  \quad \forall N\in {\cal P}(R,F)\, ,
\ee
where $\cdot$ denotes the usual matrix multiplication and the map $\phi$ is defined by the formula (\ref{phi}) in subsection \ref{Twolinearmaps}.

\smallskip
The noncommutative analogue of the matrix power is given by a repeated $\star \, \, $-multi\-pli\-ca\-tion
by the matrix $M$
\be
\lb{M^k}
M^{\overline{0}} := I\, , \qquad M^{\overline{n}}\, :=\,
\underbrace{M\star  M\star \dots \star  M}_{\mbox{\small $n$ times}}\, =\, M^{(\sigma_1\sigma_2\dots \sigma_{n-1})}\, =\,
M^{(\sigma_{n-1}\dots \sigma_2\sigma_1)}\, .
\ee
Here we introduce symbol $M^{\overline{n}}$ for the {\em $n$-th power of the matrix $M$}. The standard
matrix powers multiplication formula follows immediately from the definition
\be\lb{M^k* M^p}M^{\overline{n}} \star  M^{\overline{i}}\, =\, M^{\overline{n+i}}\, .\ee

\begin{prop}\lb{proposition4.13}
A ${\cal C}(R,F)$--module, generated by the matrix powers $M^{\overline{n}}$, $n=0,1,\dots$,
belongs to the center of the algebra ${\cal P}(R,F)$.\end{prop}

\nin {\bf Proof.~} It is sufficient to check a relation $M\star  M^{\alpha^{(i)}} = M^{\alpha^{(i)}}
\star  M$, which, in turn, follows from a calculation\\[3pt]
\phantom{a}\hfill
$ \hspace{15mm} \alpha^{(i)}\sigma_i\dots \sigma_2\sigma_1\sigma_2^{-1}\dots
\sigma_i^{-1}\, =\, \sigma_i\dots \sigma_2\sigma_1\alpha^{(i)\uparrow 1}\sigma_2^{-1}\dots
\sigma_i^{-1}\,\stackrel{\rm mod\, (\ref{red-cycl})}{=}\, \alpha^{(i)\uparrow 1} \sigma_1\, .$
\hfill$\blacksquare$

\subsection{Matrix $\star\,$-product, BMW case}\label{subsec4.4a}

It is natural to expect that the algebra ${\cal P}(R,F)$ is commutative as all of its elements are
generated by the matrix $M$ alone. We can prove the commutativity in the BMW case. Notice that
(in contrast to the Iwahori-Hecke case),
in the BMW case, the algebra ${\cal P}(R,F)$ cannot be generated by the $\star\,$-powers of $M$ only.

\smallskip
By an analogy with formula (\ref{M*}), we define a ${\cal C}(R,F)$--module map~
$\Mt : {\cal P}(R,F)$ $\rightarrow$ $ {\cal P}(R,F)$
\be
\lb{Mt}
\Mt (N) := M\cdot \xi(N), \qquad  N\in {\cal P}(R,F)\, ,
\ee
where the endomorphism $\xi$ is defined by  formula (\ref{xi}) in subsection
\ref{Twolinearmaps}. Equivalently, we can write
\be
\lb{Mt-2}
\Mt (M^{\alpha^{(n)}}) = M^{(\alpha^{(n)\uparrow 1}\kappa_1 )}
\qquad \forall\; \alpha^{(n)}\in {\cal W}_n ,\quad n=1,2,\dots \, .
\ee

\begin{prop}
\lb{proposition4.14}
Let the quantum matrix algebra ${\cal M}(R,F)$ be of the BMW type. Then the algebra ${\cal P}(R,F)$
is commutative. As a ${\cal C}(R,F)$--module, it is spanned by matrices
\be
\lb{P-gen}
M^{\overline{n}}\,  \quad \mbox{and}\quad
\Mt(M^{\overline{n+2}})\, , \quad n=0,1,\dots\, .
\ee
\end{prop}

\nin{\bf Proof.~} A proof of the last statement of the proposition goes essentially along
the same lines as the
proof of  proposition \ref{proposition4.7} and we will not repeat it. The only modification is a reduction
of the cyclic property (c.f., eqs.(\ref{cyclic}) and (\ref{red-cycl})$\,$), which finally leads to
an appearance of the additional elements $\{\Mt(M^{\overline{n}})\}_{n\geq 2}$ in the generating set.

\smallskip
To prove the commutativity of ${\cal P}(R,F)$, we derive an alternative expression for the exponent in
the matrix product formula (\ref{MaMb})
\be
\lb{a*b-2}
\alpha^{(n)}\star \beta^{(i)}\, =\, (\sigma_i^{-1}\dots \sigma_2^{-1} \sigma_1
\sigma_2\dots \sigma_i)\alpha^{(n)\uparrow i}\beta^{(i)}\, .
\ee
The calculation proceeds as follows
\ba\nonumber &&\hspace{-6mm}\alpha^{(n)}\star \beta^{(i)} =\alpha^{(n)}\beta^{(i)\uparrow n}
(\sigma_n\dots\sigma_1\sigma_2^{-1}\dots\sigma_n^{-1}) =
u_{n,i}^{(n+i)}\alpha^{(n)\uparrow i}\beta^{(i)}(u_{n,i}^{(n+i)})^{-1}
(\sigma_n\dots\sigma_1\sigma_2^{-1}\dots\sigma_n^{-1})\\[1em]\nonumber &&\hspace{-5mm}
\stackrel{\rm mod (\ref{red-cycl})}{=}\, (u_{n,i-1}^{(n+i-1)\uparrow 1})^{-1}(\sigma_2^{-1}\dots\sigma_n^{-1})
u_{n,i}^{(n+i)}\alpha^{(n)\uparrow i}\beta^{(i)} = (u_{n,i-1}^{(n+i-1)\uparrow 1})^{-1}\sigma_1
u_{n,i-1}^{(n+i-1)\uparrow 1}\alpha^{(n)\uparrow i}\beta^{(i)}\\[1em]\nonumber &&
= (\sigma_i^{-1}\dots\sigma_2^{-1}) (u_{n-1,i-1}^{(n+i-2)\uparrow 2})^{-1}\sigma_1
u_{n-1,i-1}^{(n+i-2)\uparrow 2}(\sigma_2\dots\sigma_i)\alpha^{(n)\uparrow i}\beta^{(i)} =
\mbox{right hand side of eq.(\ref{a*b-2})}.\ea
Here  we applied again the intertwining operators (\ref{u}) and used their properties (\ref{u1}) and
(\ref{u2}) and the reduced cyclicity. One more property
$$
%\be
%\lb{u3}
u_{n,i}^{(n+i)} = u_{n-1,i}^{(n+i-1)\uparrow 1} (\sigma_1\sigma_2\dots \sigma_i)
%\ee
$$
is used in the last line of the calculation.

\medskip
Due to  proposition \ref{proposition4.13}, to prove the commutativity of the algebra ${\cal P}(R,F)$, it
remains to check the commutativity of the set $\{\Mt(M^{\overline{n}})\}_{n\geq 2}$.

\smallskip
Notice that the
factors of the exponents of the matrices $\Mt(M^{\overline{n}})$ can be taken in an opposite order,
$\Mt(M^{\overline{n}}) = M^{(\kappa_1 \sigma_2 \sigma_3\dots\sigma_n)}
=M^{(\sigma_n\dots \sigma_3 \sigma_2\kappa_1)}\, .$
This observation, together with  formula (\ref{a*b-2}), allow us to choose the exponents of two matrices
$\Mt(M^{\overline{n}})\star \Mt(M^{\overline{i}})$ and $\Mt(M^{\overline{i}})\star \Mt(M^{\overline{n}})$
to be mirror (left-right) images of each other. Finally,~ $M^{\alpha^{(n)}}=M^{\varsigma(\alpha^{(n)})},\; \forall\, \alpha^{(n)}\in {\cal W}_n(q,\mu)$, where $\varsigma$ is the antiautomorpism (\ref{antvs}), since both sides of this equality can be expanded into linear combinations of the generators
(\ref{P-gen}), which are invariant with respect to the mirror reflection of their exponents, and since
the expansion rules (i.e. the defining relations for the BMW algebras) are mirror symmetric.
\hfill$\blacksquare$

\begin{lem}\lb{lemma4.15}
{}For the BMW type quantum matrix algebra ${\cal M}(R,F)$, one has
\ba
\lb{Mt-IM}\Mt (I)\, =\, \mu M\, , \qquad \Mt (M)\, =\, \mu^{-1} I\, g\, ,
\\[1em]
\lb{Mt2-N}
\Mt (\Mt (N))\, =\, N\, g\,  \quad\forall\; N\in {\cal P}(R,F)\, .
\ea
\end{lem}

\nin {\bf Proof.~} The relations (\ref{Mt-IM}) follow immediately from the relations
(\ref{tau2}) and (\ref{traceDK}) and the definitions (\ref{Mt}) and (\ref{xi}).

\smallskip
As for the equality (\ref{Mt2-N}), it is enough to check it in the case when the matrix $N$ is a power
of the matrix $M$.

\smallskip
To evaluate the expression
$\Mt(\Mt(M^{\overline n}))=M^{(\kappa_1\kappa_2\sigma_3\dots\sigma_{n+1})}$, we transform its exponent,
using the relations (\ref{bmw3}) in the BMW algebra and the reduced cyclic property, to
\be
\begin{array}{ccl}
\kappa_1\kappa_2\sigma_3\dots\sigma_{n+1}&=&\kappa_1(\kappa_2\kappa_3\sigma_2^{-1})
\sigma_4\dots\sigma_{n+1}\,\stackrel{\rm mod\, (\ref{red-cycl})}{=}\, (\sigma_2^{-1}\kappa_1\kappa_2)
\kappa_3\sigma_4\dots\sigma_{n+1}\\[1em] &=&\sigma_1\kappa_2\kappa_3\sigma_4\dots\sigma_{n+1}
\, =\; \dots\;\stackrel{\rm mod\, (\ref{red-cycl})}{=}\,\sigma_1\sigma_2\dots\sigma_{n-1}\kappa_n
\kappa_{n+1}\, .
\end{array}
\lb{trans-exp}
\ee
{}For the exponent (\ref{trans-exp}), the matrix power is easily calculated, again with the help
of the relations (\ref{tau2}) and (\ref{traceDK}), and gives the expression $M^{\overline{n}}g$. \hfill$\blacksquare$

\medskip
The last relation in (\ref{Mt-IM}) shows that to introduce the inverse matrix to the matrix $M$ it is
sufficient to add the inverse $g^{-1}$ of the 2-contraction $g$ to the algebra ${\cal M}(R,F)$.
This is realized in the next subsection.

\subsection{Matrix inversion}
\label{subsec4.5}

In this subsection we define an {\em extended} quantum matrix algebra, to which the
inverse of the quantum matrix belongs.

\begin{lem}\lb{lemma4.16}
Let ${\cal M}(R,F)$ be the BMW type quantum matrix algebra. Its 2-contraction $g$ fulfills a relation
\be\lb{Mj} M\, g\, =\, g\, ( G^{-1} M G)\, ,\ee
where $G$ is defined by  formula (\ref{G}).\end{lem}

\nin {\bf Proof.~} The proof consists of a calculation
\be
\begin{array}{ccl}
\phantom{a}\hspace{-10mm} M_1\left(g K_2\right) &=&\mu^2 M_{\overline{1}}
M_{\overline{2}}M_{\overline{3}} K_2\, =\,\mu^2 M_{\overline{1}}M_{\overline{2}}M_{\overline{3}}K_2K_1K_2\,
=\,\mu^2 K_2 \left( M_{\overline{1}}M_{\overline{2}} K_1\right) M_{\overline{3}} K_2\\[1em] &=& gK_2K_1
M_{\overline{3}} K_2 \, =\, \left(g K_2\right) \tr_{(2,3)}\left(K_2 K_1 M_{\overline{3}}\right)
\\[1em] &=& \left(g K_2\right)\tr_{(2,3)}\left(K_2 F_2 F_1 K_2 M_1 F_1^{-1} F_2^{-1}\right)
\\[1em] &=& \left(g K_2\right)\tr_{(2,3)}(F_2 F_1 K_2) M_1
\tr_{(2,3)}(K_2 F_1^{-1} F_2^{-1}) =\left(g K_2\right) (G^{-1} M G)_1 .
\lb{MjK}
\end{array}
\ee
Here the relations (\ref{tau2}) and (\ref{bmw5a}) were used in the first two lines; the property
${\rm rk}\, K=1$ was
used in the last/first equality of the second/fourth line; the definition of $M_{\overline{3}}$ was
substituted and the twist relation for the pair $\{K,F\}$ was used in the third line; the formulas (\ref{G})
and (\ref{G-inv}) for $G$ and $G^{-1}$ were substituted in the last equality. \hfill$\blacksquare$

\begin{def-prop}\lb{proposition4.17}
Let ${\cal M}(R,F)$ be the BMW type quantum matrix algebra. Consider an extension of
the algebra ${\cal M}(R,F)$ by a generator $g^{-1}$ subject to relations
\be\lb{j-inv} g^{-1}\, g\, =\, g\, g^{-1}\, =\, 1\, , \qquad g^{-1}\, M\, =\, (G^{-1}MG)\, g^{-1}\, .\ee
The extended algebra, which we shall further denote by ${\cal M^{^\bullet\!}}(R,F)$, contains
an inverse matrix to the matrix $M$
\be\lb{M-inv} M^{-1}\, :=\, \mu\, \xi(M)\, g^{-1}\, :\qquad
M\cdot M^{-1}\, =\, M^{-1}\cdot M\, =\, I\, .\ee\end{def-prop}

\nin {\bf Proof.~} Lemma \ref{lemma4.16} ensures the consistency of the relations (\ref{j-inv}). The equality
$M\cdot M^{-1}=I$ for the inverse matrix (\ref{M-inv}) follows immediately from the formulas
(\ref{Mt-IM}) and (\ref{Mt}).

\smallskip
To prove the equality $M^{-1}\cdot M = I$, consider a mirror partner of the map $\xi$:
\be
\lb{Pi}
\theta(M)\, :=\, \mu^{-2}\, \Tr{2} K_1 M_{\overline{2}}\, .
\ee
By the (left-right) symmetry arguments in the assumptions of lemma \ref{lemma3.9},
the map $\theta$ is invertible and the inverse map reads
\be
\lb{Pi-inv}
\theta^{-1}(M)\, =\,\TR{2}{R_f}\left( F^{-1}_1 K_1 M_1  F_1\right)\, .
\ee
Applying in a standard way the transformation formula (\ref{inv-trD}), we calculate a composition of
the maps $\xi$ and $\theta$,
\ba
\lb{composition}
\xi(\theta(M))_1 = \theta(\xi(M))_1 =
\mu^{-2}\, \Tr{2,3} K_2 K_1 M_{\overline{3}} =\tr_{(2,3)} K_2 K_1 M_{\overline{3}} =(G^{-1} M G)_1\, .
\ea
Here the relation (\ref{KDD}) was used to substitute the R-traces by the usual traces; the last
equality follows
from a comparison of the second and the last lines in the calculation (\ref{MjK}).

\smallskip
Now we observe that, in view of the relations (\ref{tau2}) and (\ref{traceDK}), a matrix
$(^{-1}M)\, :=\, \mu\, g^{-1}\, \theta^{-1}(M)$ fulfills the relation $(^{-1}M)\cdot M=I$. The identity
$(^{-1} M) = M^{-1}$ follows then from the relations (\ref{composition}) and (\ref{Mj}). \hfill$\blacksquare$

\begin{rem}\lb{remark4.18}
{\rm One can generalize the definitions of the characteristic subalgebra and of the matrix powers to the case of
the extended quantum matrix algebra  ${\cal M^{^\bullet\!}}(R,F)$. Not going into details, we just mention
that the {\em extended} characteristic subalgebra ${\cal C^{^\bullet\!}}(R,F)$ is generated by the set
$\{g,g^{-1},p_i\}_{i\geq 0}$ and the {\em extended} algebra ${\cal P^{^\bullet\!}}(R,F)$, as a
${\cal C^{^\bullet\!}}(R,F)$--module, is spanned by matrices
$$M^{\overline{n}}\, \quad \mbox{and}\quad\Mt(M^{\overline{n}})\, \quad \forall\; n\in{\Bbb Z}\, .$$
Here inverse powers of $M$ are defined through the repeated $\star\, $-multiplication by
$M^{\overline{-1}}$, which is given by
\be
\lb{Minv*N}
M^{\overline{-1}}\star  N \, :=\, N \star  M^{\overline{-1}}\, :=\,\phi^{-1}(M^{-1}\cdot N)
\,  \quad
\forall\; N\in{\cal P^{^\bullet\!}}(R,F)\, .
\ee
Explicitly, one has
$$M^{\overline{-n}}\, :=\,\underbrace{M^{\overline{-1}}\star \dots \star  M^{\overline{-1}}\star }_{n\ times}I\, =\,
{\tr_{\!\!R_f^{-1}}}_{(2,\dots,n+1)}\!\left( M^{-1}_{\underline{2}} M^{-1}_{\underline{3}}\dots
M^{-1}_{\underline{n+1}}\, \rho_{R_f^{-1}}(\sigma_n\dots\sigma_2\sigma_1)\right) ;$$
where the copies $M^{-1}_{\underline{i}}$ of the matrix $M^{-1}$  are defined as (c.f.
with eq.(\ref{kopii})$\,$)
\be
\lb{kopii-2}
M_{\underline{1}}\, :=\, M_1 ,\qquad M_{\underline{i+1}}\, :=\,
{}F_i^{-1}M_{\underline{i}}\, F_i , \qquad i=2,3,\dots .
\ee
Notice that in general $M^{\overline{-1}} =
\phi^{-1}(M^{-1})\neq M^{-1}$. Here are some particular examples of the  multiplication
by $M^{\overline{-1}}$
$$M^{\overline{-n}}\star M^{\overline{i}}\, =\, M^{\overline{i-n}}\, ,\qquad
M^{\overline{-1}}\star M^{\alpha^{(n)}\uparrow 1}\, =\, ch(\alpha^{(n)})\, I\, .$$ }\end{rem}

\section{Relations for generating sets of the characteristic subalgebra: BMW case.}
\lb{sec5}

In this last section we use the basic identities from subsection \ref{subsec5.1} to establish relations
between the three sets of elements in  the characteristic subalgebra ---
$\{g,a_i\}_{i\geq 0}$, $\{g,s_i\}_{i\geq 0}$ and the power sums $\{g,p_i\}_{i\geq 0}$. As a byproduct,
we prove the  proposition \ref{proposition4.8}.

\medskip
Before we proceed, let us recall the initial data of the construction.
\begin{itemize}
\item
Given a compatible pair of R-matrices $\{R,F\}$, in which the operator
$F$ is strict skew invertible and the operator $R$ is skew invertible of the BMW type (and,
hence, strict skew invertible), we introduce the BMW type quantum matrix algebra ${\cal M}(R,F)$ (see definition \ref{definition4.1});
\item
Assuming additionally that the eigenvalues $q$ and $\mu$ of the R-matrix $R$ (i.e., the parameters of the BMW algebras, whose
representations are generated by the matrix $R$) satisfy conditions
$i_q\neq 0,\;\mu\neq -q^{3-2i}\;\;\forall
\; i=2,3,\dots ,n$ (see (\ref{mu})) we can consistently define the antisymmetrizers $a^{(i)}$
and introduce skew powers of the quantum matrix $M$: $M^{a^{(i)}}$, $0\leq i\leq n$.
\end{itemize}

\subsection{Basic identities}\lb{subsec5.1}

In this subsection we establish relations between `descendants' of the matrices $M^{a^{(i)}}$ in the
algebra ${\cal P}(R,F)$. These relations are used later in a derivation of  the Newton relations.

\medskip
{}For $1\leq i\leq n$ ~and~ $m\geq 0$, we consider two series of descendants of $M^{a^{(i)}}$:
\be
\lb{La}
A^{(m,i)}\ :=\ i_q\, M^{\overline{m}}  \star  M^{a^{(i)}}\,\ , \quad B^{(m+1,i)}\ :=\
i_q\, M^{\overline{m}} \star  \Mt (M^{a^{(i)}})\,  .
\ee

It is  suitable to define $A^{(m,i)}$ and $B^{(m,i)}$ for boundary values of their indices
\ba\lb{AB-boundary} & A^{(-1,i)}\, :=\, i_q\, \phi^{-1}\left(
\Tr{2,3,\dots i} M_{\overline{2}}M_{\overline{3}}\dots M_{\overline{i}}\, \rho_R(a^{(i)})
\right)\ , \quad B^{(0,i)}\ :=\ i_q\, \phi^{-1}\bigl(\xi\bigl(M^{a^{(i)}}\bigr)\bigr)\, &\ea
and
\ba\lb{LT0}
&A^{(m,0)}\ :=0\ \ \ {\mathrm{and}}\ \ \ \ B^{(m,0)}\ :=\ 0\,\qquad \forall\; m\geq 0\, .&\ea
Notice that although the elements $A^{(-1,i)}$ and $B^{(0,i)}$ do not, in general, belong to the algebra ${\cal P}(R,F)$,
their descendants  $A^{(-1,i)} g$ and $B^{(0,i)} g$ do (see eqs.(\ref{rek1}) and (\ref{rek2})  in the case $m=0$).

\smallskip
In the case when the contraction $g$ (and, hence, the matrix $M$) is invertible, the formulas (\ref{La}), with $m$ now an arbitrary integer,
can be used  to define descendants of $M^{a^{(i)}}$ in the extended algebra  ${\cal P^{^\bullet\!}}(R,F)$
(see the remark \ref{remark4.18}). In this case, the matrices $A^{(-1,i)}$ and $B^{(0,i)}$ are expressed uniformly:
$A^{(-1,i)} = i_q M^{\overline{-1}} \star  M^{a^{(i)}}$,~~ $B^{(0,i)} = i_q M^{\overline{-1}} \star  \Mt (M^{a^{(i)}})$.

\begin{lem}\lb{lemma5.1}
{}For ~$0\leq i\leq n-1$~ and ~$m\geq 0$, the matrices $A^{(m-1,i+1)}$ and $B^{(m+1,i+1)}$
satisfy recurrent relations
\ba\lb{rek1} A^{(m-1,i+1)} &=& q^i M^{\overline{m}}\, a_i\, -\,
A^{(m,i)}\, -\, { \mu q^{2i-1}(q-q^{-1})\over 1+\mu q^{2i-1}}\ B^{(m,i)}\, ,\\[1em]\lb{rek2}
B^{(m+1,i+1)} &=&\Bigl( \mu^{-1}q^{-i} M^{\overline{m}}\, a_i\, +\,
{q-q^{-1}\over 1+\mu q^{2i-1}}\ A^{(m,i)}\, -\,B^{(m,i)}\Bigr) g\,  .\ea\end{lem}

\noindent {\bf Proof.}~
For $i=0$ relations (\ref{rek1}) and (\ref{rek2}) by (\ref{LT0}) simplify to
$$
%\be
%\lb{rek-i=0}
A^{(m-1,1)} = M^{\overline{m}}\, , \quad B^{(m+1,1)} = \mu^{-1} M^{\overline{m}} g\, .
%\ee
$$
They follow from eqs. (\ref{M^k* M^p}), (\ref{Mt-IM}).

\smallskip
Let us check (\ref{rek1})  for $i>0$.
For $m \geq 0$, we calculate
\ba\nonumber A^{(m,i+1)} = (i+1)_q M^{(a^{(i+1)\uparrow m}\, \sigma_m\dots \sigma_2\sigma_1)}
\, =\, q^i M^{(a^{(i)\uparrow(m+1)}\,\sigma^{-}_{m+1}(q^{-2i})\,\sigma_m\dots \sigma_2\sigma_1)}
\\[1em]\nonumber =\, q^i M^{\overline{m+1}}\, a_i\, -\, A^{(m+1,i)}\, -\,
{\mu q^{2i-1}(q-q^{-1})\over 1+\mu q^{2i-1}}\, B^{(m+1,i)}\, .\hspace{20mm}\ea
Here in the first line we used the second formula from (\ref{a^k}) for $a^{(i+1)\uparrow m}$ and applied
the reduced cyclic property (\ref{red-cycl}) and the relations (\ref{idemp-2}) to cancel one of two terms
$a^{(i)\uparrow (m+1)}$. In the second line we substituted the formula (\ref{ansatz}) for the baxterized
elements $\sigma_{m+1}^-(q^{-2k})$ and applied the relation (\ref{char1}) to simplify the first term in the sum.

\smallskip
For $A^{(-1,i+1)}$, the relations (\ref{rek1}) are verified similarly
$$
\begin{array}{ccl}
A^{(-1,i+1)} &=& q^{i}\, \phi^{-1}\left( \Tr{2,3,\dots i+1} M_{\overline{2}}
M_{\overline{3}}\dots M_{\overline{i+1}}\rho_R(a^{(i)\uparrow 1}\sigma^{-}_1(q^{-2i}))\right)\,
\\[1em]
&=&
q^i\, \phi^{-1}(I)\,  a_i\, -\, i_q\, \phi^{-1}(\phi(M^{a{(i)}}))\, -\, {\displaystyle \frac{\mu q^{2i-1}(q-q^{-1})}
{1+\mu q^{2i-1}}}\, i_q\, \phi^{-1}(\xi(M^{a^{(i)}}))\,
\\[1em]
&=&
q^i\, I\, a_i \, -\, A^{(0,i)}\, -\,
{\displaystyle \frac{\mu q^{2i-1}(q-q^{-1})}
{1+\mu q^{2i-1}}}\, B^{(0,i)}\, .
\end{array}
$$
Here the definitions (\ref{phi}) and (\ref{xi}) of the endomorphisms $\phi$ and $\xi$ were additionally taken into account.

\smallskip
To prove (\ref{rek2})  for $i>0$ we proceed in the same way
\be
\lb{kusok}
\!\!\!\!\!\!\!\begin{array}{l} B^{(m+1,i+1)} = (i+1)_q
M^{(a^{(i+1)\uparrow m+1}\, \kappa_{m+1}\sigma_{m}\dots \sigma_2\sigma_1)}\, =\, q^i
M^{(a^{(i)\uparrow m+2}\, \sigma^{-}_{m+2}(q^{-2i})\,\kappa_{m+1}\sigma_{m}\dots \sigma_2\sigma_1)}
\hspace{13mm}\\[1em] = q^{-i} M^{\overline{m}}\!\star \! \Mt(M) a_i-
i_q M^{(a^{(i)\uparrow m+2 }\,\sigma^{-1}_{m+2}\kappa_{m+1}\sigma_{m}\dots\sigma_1)} +
\displaystyle {q^i-q^{-i}\over 1+\mu q^{2i-1}}\, M^{\overline{m}}\! \star \!\Mt(\Mt(M^{a^{(i)}})).\hspace{3mm}\end{array}
\ee
Here in the second line we used another expression for the baxterized generators
$$\sigma_{i}^{\varepsilon}(x)\, =\, x 1\, +\, {x-1\over q-q^{-1}}\, \sigma^{-1}_{i}\, -\,
{\alpha_{\varepsilon}x(x-1)\over \alpha_{\varepsilon}x+1}\,\kappa_{i}\, ,$$
which follows by a substitution $\sigma_i = \sigma_i^{-1} + (q-q^{-1})(1-\kappa_i)$ into the original
expression (\ref{ansatz}).

\smallskip
Now, notice that
\be\lb{kakaka}\sigma^{-1}_3\kappa_2\sigma_1\, =\,\sigma^{-1}_3\kappa_2\kappa_1\sigma_2^{-1}\,
\stackrel{\rm mod\, (\ref{red-cycl})}{=}\,\sigma_2^{-1}\sigma^{-1}_3\kappa_2\kappa_1\, =\,
\kappa_3\kappa_2\kappa_1\, ,\ee
and, hence, in the case $m\geq 1$, the second term in the last line of the equality
(\ref{kusok}) can be expressed as
\be\lb{kusok2} -\, i_q M^{(a^{(i)\uparrow m+2}\,\sigma^{-1}_{m+2}\kappa_{m+1}\sigma_{m}\dots\sigma_1)}\,
=\, - i_q M^{\overline{m-1}}\star \Mt(\Mt(\Mt(M^{a^{(i)}})))\, .\ee
Applying then the formulas (\ref{Mt-IM}) and (\ref{Mt2-N}) to the expressions (\ref{kusok}) and
(\ref{kusok2}), we complete verification of (\ref{rek2}) for $m\geq 1$.

\smallskip
For the case $m=0$, the transformation of the second term in (\ref{kusok}) should be slightly modified.
Notice that by eq.(\ref{kakaka}),
$$\phi(M^{a^{(i)\uparrow 2}\sigma_2^{-1}\kappa_1})\, =\, \xi(\Mt(\Mt(M^{a^{(i)}})))\, .$$
Inverting the endomorphism $\phi$ in this formula
and using the relation (\ref{Mt2-N}) and the definition of $B^{(0,i)}$  (\ref{AB-boundary}), we complete
the transformation of the second term in (\ref{kusok}) and, again, get the equality (\ref{rek2}).
\hfill$\blacksquare$

\subsection{Newton and Wronski relations}
\lb{newtonbmw}

\begin{theor}\lb{theorem6.1}
Let ${\cal M}(R,F)$ be a BMW type quantum matrix algebra. Assume that its two parameters $q$
and $\mu$ satisfy the conditions (\ref{mu}), which allow to introduce either the set
$\{a_i\}_{i=0}^n$ or, respectively, the set $\{s_i\}_{i=0}^n$ in the characteristic subalgebra
${\cal C}(R,F)$ (see the definitions (\ref{SA_0}) and (\ref{SA_k})$\,$). Then the following
Newton recurrent formulas
relating, respectively, the sets $\{a_i,g\}_{i=0}^n$, or $\{s_i,g\}_{i=0}^n$ to the set of the  power sums (see the definitions (\ref{P-01}) and (\ref{P_k})$\,$)
\ba\lb{Newton-a}\sum_{i=0}^{n-1} (-q)^i a_i\, p_{n-i} &=& (-1)^{n-1} n_q\, a_n \, +\, (-1)^n
\sum_{i=1}^{\lfloor {n/2}\rfloor}\Bigl( \mu q^{n-2i} -q^{1-n+2i}\Bigr)\, a_{n-2i}\, g^i\, \ea
and
\ba\lb{Newton-s}
\sum_{i=0}^{n-1} q^{-i} s_i\, p_{n-i} &=&  n_q\, s_n \, +\,\sum_{i=1}^{\lfloor {n/2}\rfloor}\Bigl(
\mu q^{2i-n} + q^{n-2i-1}\Bigr)\, s_{n-2i}\, g^i\, \ea
are fulfilled.

\smallskip
In the case, when both sets $\{a_i,g\}_{i=0}^n$ and $\{s_i,g\}_{i=0}^n$
are consistently defined, they satisfy the Wronski relations
\be\lb{Wronski}\sum_{i=0}^n (-1)^i a_i\, s_{n-i}\, =\, \delta_{n,0} -\delta_{n,2}\, g\, ,\ee
where $\delta_{i,j}$ is a Kronecker symbol.\end{theor}

\begin{rem}\lb{remark6.2}
{\rm One can use the formulas (\ref{Newton-a}) and (\ref{Newton-s}) for an iterative definition of the
elements $a_i$ and  $s_i$ for $i\geq 1$, with initial conditions $a_0=s_0=1$. In this case, the elements
$a_n$ and $s_n$ are well defined, assuming that $i_q\neq 0 \;\; \forall\; i=2,3,\dots n$. The additional
restrictions on the parameter $\mu$, which appeared in their
initial definition (\ref{SA_k}), are artifacts of
the use of the antisymmetrizers and symmetrizers $a^{(n)}, s^{(n)}\in {\cal W}_n(q)$.}\end{rem}

\noindent {\bf Proof.}~ We prove the  relation (\ref{Newton-a}). Denote
$$
%\be
%\lb{J_i}
J^{(0)} := 0, \quad J^{(i)} := {\displaystyle \sum_{j=0}^{i-1}} (-q)^j M^{\overline{i-j}} a_j,
\quad i=1,2,\dots ,n\ .
%\ee
$$
We are going to find an expression for the matrix $J^{(n)}$ in terms of
the matrices $A^{(0,i)}$ and $B^{(0,i)}$, $1\leq i\leq n$.

\medskip
As we shall see, there exist matrices $H^{(i)}$, which fulfill  equations
\be
\lb{H-1}
(1-q^2) H^{(i)}\, g\, =\,\Bigl(J^{(i)} + (-1)^i A^{(0,i)} \Bigr)\,  , \quad
i=0,1,\dots ,n.
\ee
To calculate the matrices $H^{(i)}$,
we substitute repeatedly the relations (\ref{rek1}) for the elements
$A^{(0,i)}$, $A^{(1,i-1)}, \dots , A^{(i-1,1)}$ in the right hand side of
eq.(\ref{H-1}). It then transforms to
\be
\lb{H-2}
H^{(i)}\, g\, =\, -\mu q^{-1}\sum_{j=1}^{i-1} (-1)^j
{q^{2j-1}\over 1+\mu q^{2j-1}} B^{(i-j,j)}\,  , \quad
i=0,1,\dots ,n.
\ee
Now, using the expressions (\ref{rek2}) for the elements  $B^{(i-j,j)}$, one can check that matrices
\ba
\lb{H}
H^{(0)}& := & H^{(1)}\, :=\, 0,
\\[2pt]
\nn
H^{(i)}& :=&\sum_{j=0}^{i-2} {(-q)^{j}\over 1+\mu q^{2j+1}}\Bigl(
M^{\overline{i-j-2}} a_j\, +\, {\mu q^{j}(q-q^{-1})\over 1+\mu q^{2j-1}} A^{(i-j-2,j)}\, -\,
\mu q^j B^{(i-j-2,j)}\Bigr) , \;\;
i=2,\dots ,n.
\ea
satisfy eq.(\ref{H-2}).

\smallskip
Next, consider a combination $(H^{(i+2)}- H^{(i)} g)$. Using eq.(\ref{H})
for the first term and eq.(\ref{H-2}) for the second term, we calculate
\ba
\nonumber
H^{(i+2)} - H^{(i)} g &=&\sum_{j=0}^{i-1} {(-q)^j\over 1 +\mu q^{2j+1}}\Bigl(
M^{\overline{i-j}} a_j\, +\, {\mu q^j (q-q^{-1})\over 1+ \mu q^{2j-1}} (A^{(i-j,j)}-q^{-1}B^{(i-j,j)})
\Bigr)\,
\\[2pt]
\nonumber
&&\hspace{15mm} + \,{(-q)^i\over 1+ \mu q^{2i+1}}\Bigl(
I a_i + {\mu q^i (q-q^{-1})\over 1+ \mu q^{2i-1}} A^{(0,i)}\, -\,\mu q^i B^{(0,i)}\Bigr)\, , \;\; \forall\, i=0,\dots , n.
\ea
To continue, we need the  following auxiliary result:

\begin{lem}
\lb{lemma6.3}
{}For $1\leq i\leq n$, one has
\be
\lb{A-tri}
{(-1)^{i-1} A^{(0,i)}\over 1+\mu q^{2i-1}}\, =\,\sum_{j=0}^{i-1}
{(-q)^j\over 1+\mu q^{2j+1}}\Bigl( M^{\overline{i-j}} a_j\, +\, {\mu q^j(q-q^{-1})\over 1+\mu q^{2j-1}}
(A^{(i-j,j)}-q^{-1} B^{(i-j,j)})\Bigr) .
\ee
\end{lem}

\noindent {\bf Proof.}~ Use the recursion (\ref{rek1}) for $A^{(i-j-1,j+1)}$ to calculate
$$
{A^{(i-j-1,j+1)}\over 1+\mu q^{2j+1}} +{A^{(i-j,j)}\over 1+\mu q^{2j-1}} =
{q^{j}\over 1+\mu q^{2j+1}}\Bigl( M^{\overline{i-j}}a_j +{\mu q^j (q-q^{-1})\over 1+\mu q^{2j-1}}
(A^{(i-j,j)}-q^{-1} B^{(i-j,j)}) \Bigr)\, .
$$
Compose an alternating sum of the above relations for $0\leq j \leq i-1$
and take into account the condition $A^{(i,0)}=0$. \hfill$\blacksquare$

\medskip
Using the relation (\ref{A-tri}), we finish the calculation
\be
\lb{H-H}
H^{(i+2)} - H^{(i)} g \, =\, (1+\mu q^{2i+1})^{-1}\Bigl(
(-q)^i I a_i\, +\, (-1)^{i+1}(A^{(0,i)}+\mu q^{2i} B^{(0,i)})\Bigr)\,\;\; \forall\; i=0,\dots ,n-2.
\ee

Now it is straightforward to get
\be
\lb{H-otvet}
H^{(i)}\, =\,\sum_{j=1}^{[i/2]}{(-1)^{i-1}\over
1+\mu q^{2(i-2j)+1}}\Bigl( A^{(0,i-2j)} + \mu q^{2(i-2j)} B^{(0,i-2j)} - q^{i-2j} I a_{i-2j}
\Bigr)\, g^{j-1}\, ,\;\; \forall\, i=0,\dots ,n,
\ee
where $[k]$ denotes the integer part of the number $k$.
Finally, substituting the expression (\ref{H-otvet}) back into eq.(\ref{H-1}), we obtain a formula
\be
\lb{J-otvet}
J^{(i)}\, =\, (-1)^{i-1} A^{(0,i)}\, +\,
\sum_{j=1}^{[i/2]}{(-1)^{i-1}(1-q^2)\over 1+\mu q^{2(i-2j)+1}}\Bigl(
A^{(0,i-2j)} + \mu q^{2(i-2j)} B^{(0,i-2j)} - q^{i-2j} I a_{i-2j}\Bigr) g^j ,
\ee
which is valid for $0\leq i\leq n$.

\smallskip
Taking the R-trace of eq.(\ref{J-otvet}), we obtain the Newton relations (\ref{Newton-a}). Here, in
the calculation of the R-trace of $B^{(0,i-2j)}$, we took into account the formulas (\ref{phi-xi-traces}).

\medskip
The formulas (\ref{Newton-s}) can be deduced from the relations (\ref{Newton-a}) by a substitution
$q\rightarrow -q^{-1}$,~ $a_j\rightarrow s_j$. This is justified by the existence of the BMW algebras
homomorphism (\ref{homS-A}) $\iota:\;{\cal W}_n(q,\mu)\rightarrow {\cal W}_n(-q^{-1},\mu)$ and a fact that one and
the same R-matrix $R$ generates representations of both algebras ${\cal W}_n(q,\mu)$ and
${\cal W}_n(-q^{-1},\mu)$.

\medskip
The relation (\ref{Wronski}) is proved by induction on $n$. The cases $n=0,1,2$ are easily checked with the use of
eqs.(\ref{Newton-a}) and (\ref{Newton-s}). Then, making an induction assumption, we derive the Wronski relations
for arbitrary $n>2$. To this end, we take a difference of  eqs.(\ref{Newton-s}) and (\ref{Newton-a})
$${\displaystyle \sum_{i=0}^{n-1}} \Bigl( q^{-i} s_i\, p_{n-i}\, -\, (-q)^i a_i\, p_{n-i}
\Bigr)\, =\, n_q (s_n\, +\, (-1)^n a_n )\, +\,\mbox{terms proportional to $g$}$$
and substitute for $p_{n-i}$ in the first/second term of the left hand side its expression from the Newton
relation (\ref{Newton-a})/(\ref{Newton-s}) (with  $n$ replaced by $n-i$). As a result, all terms,
containing the power sums, cancel and, after rearranging the summations, we get
$$n_q \sum_{i=0}^n (-1)^i a_i s_{n-i}\, =\, -\sum_{i=1}^{[n/2]}
(q^{1-n+2i}+q^{n-1-2i})g^i\sum_{j=0}^{n-2i}(-1)^j a_j s_{n-2i-j}\, .$$
By the induction assumption, the double sum in the right hand side of this relation vanishes identically:
when $n$ is odd, the second sum vanishes for all values of the index $i$; when $n$ is even, the second sum
is different from zero only for two values of the index $i$, $i=n/2$ and $i=n/2-1$, and
these two summands cancel.
\hfill$\blacksquare$

\appendix
\section{Primitivity of contractors}
\lb{primcontr}

\def\theequation{\thesection.\arabic{equation}}
\makeatletter\@addtoreset{equation}{section}\makeatother

In this appendix we return to the consideration of the contractors in the BMW algebra.
We shall establish useful properties of the contractors in lemmas \ref{zapo}, \ref{fuproco1}
and then use it to demonstrate their primitivity (announced in proposition
\ref{proposition2.2} in subsection \ref{subsec2.4})
in proposition \ref{primco}.

\medskip
In this appendix we shall denote by ${\cal W}(\sigma_i,\sigma_{i+1},\dots ,\sigma_j)$,
where $i\leq j$, the BMW algebra with the generators $\sigma_i,\sigma_{i+1},\dots ,\sigma_j$
(the values of the parameters $q$ and $\mu$ are fixed).

\begin{lem}
\lb{zapo}
~~Let~ $\alpha\in {\cal W}(\sigma_1,\sigma_{2},\dots ,\sigma_j)$,~ where
$j\geq n$.~ Then there exists an element
$\tilde{\alpha}\in {\cal W}(\sigma_{n+1},\sigma_{n+2},\dots ,\sigma_j)$ such that
$$
%\be
%\lb{zapof}
c^{(2n)}\alpha =c^{(2n)}\tilde{\alpha}\ .
%\ee
$$
\end{lem}

\noindent {\bf Proof.}~ Assume that $\alpha\in {\cal W}(\sigma_i,\sigma_{i+1},\dots ,\sigma_j)$
and $\alpha\notin {\cal W}(\sigma_{i+1},\dots ,\sigma_j)$. If $i>n$ then there is nothing to prove.

\smallskip
{}For $i\leq n$, we shall prove that there exists an element
$\alpha'\in {\cal W}(\sigma_{i+1},\dots ,\sigma_j)$ such that
$$
%\be
%\lb{zapod1}
c^{(2n)}\alpha =c^{(2n)}\alpha'\ .
%\ee
$$
Given this statement, the proof follows by induction on $i$.

\medskip
Due to the formula (\ref{char8}), we can express the element $\alpha$ as a linear combination of
elements of the form $xu_i \bar{x}$, where $x,\bar{x}\in {\cal W}(\sigma_{i+1},\dots ,\sigma_j)$
and $u_i$ is equal to 1, $\sigma_i$ or $\kappa_i$. The terms with $u_i=1$ belong already to
${\cal W}(\sigma_{i+1},\dots ,\sigma_j)$ so we may assume that the element $u_i$ is non-trivial
(that is, equals $\sigma_i$ or $\kappa_i$).

\smallskip
We express now the element $x$ as a linear combination of the elements of the form
$yu_{i+1}\bar{y}$, where $y,\bar{y}\in {\cal W}(\sigma_{i+2},\dots ,\sigma_j)$
and $u_{i+1}$ is equal to 1, $\sigma_{i+1}$ or $\kappa_{i+1}$. Each element $\bar{y}$ commutes
with the element $u_i$ thus the element $\alpha$ becomes a linear combination of elements of
the form $yu_{i+1}u_i\bar{\bar{x}}$ with $y\in {\cal W}(\sigma_{i+2},\dots ,\sigma_j)$ and
$\bar{\bar{x}}\in {\cal W}(\sigma_{i+1},\dots ,\sigma_j)$. In the terms with $u_{i+1}=1$ we move
the element $y$ rightwards through the element $u_i$ and continue the process for the terms
with $u_{i+1}$ equal to $\sigma_{i+1}$ or $\kappa_{i+1}$. After a finite number of steps the
process terminates and we will have an expression for the element $\alpha$ as a linear
combination of terms
\be
u_{i+k}\dots u_{i+1}u_iz\ ,
\lb{zapod2}
\ee
where the element $z$ belongs to ${\cal W}(\sigma_{i+1},\dots ,\sigma_j)$ and each of the
elements $u_{i+s}$, $s=0,1,\dots k$, is equal to $\sigma_{i+s}$ or $\kappa_{i+s}$.

\medskip
Let us first analyze expressions (\ref{zapod2}) with $i+k>n$. The contractor $c^{(2n)}$ is
divisible by the element $\kappa_n$ from the right due to the relation (\ref{idemp-c1}).
The element $\kappa_n$ can move rightwards in the product $c^{(2n)}u_{i+k}\dots u_{i+1}u_iz$
until it reaches the element $u_{n+1}$ and we arrive at the expression
$\dots\kappa_n u_{n+1}u_n\dots$. For all four possibilities
($\sigma_{n+1}\sigma_n$, $\sigma_{n+1}\kappa_n$, $\kappa_{n+1}\sigma_n$ or
$\kappa_{n+1}\kappa_n$) for the product $u_{n+1}u_n$, the expression $\kappa_n u_{n+1}u_n$
can be rewritten, with the help of the relations (\ref{bmw2b})--(\ref{bmw5a}), in a form
$\kappa_n v_{n+1}$, where $v_{n+1}$ is a polynomial in $\sigma_{n+1}$. Moving the element
$\kappa_n$ back to the contractor $c^{(2n)}$, we obtain
$$
%\be
c^{(2n)}u_{i+k}\dots u_{i+1}u_iz =c^{(2n)}u_{i+k}\dots u_{n+2}v_{n+1}\cdot u_{n-1}\dots
u_iz=c^{(2n)}u_{n-1}\dots u_i\bar{z}
%\lb{zapod3}
%\ee
$$
with some other $\bar{z}\in {\cal W}(\sigma_{i+1},\dots ,\sigma_j)$.

\medskip
Thus we can rewrite the product of the contractor $c^{(2n)}$ by an expression (\ref{zapod2})
with $i+k>n$ as a product of $c^{(2n)}$ with an expression of the same form (\ref{zapod2})
but with $i+k<n$.

\smallskip
Now using the relations (\ref{idemp-c2}) we remove the elements $u_{i+k}$ one by one to the
right:
$$
%\be
c^{(2n)}u_{i+k}\dots u_{i+1}u_i=c^{(2n)}u_{n-i-k}u_{i+k-1}\dots u_{i+1}u_i=
c^{(2n)}u_{i+k-1}\dots u_{i+1}u_iu_{n-i-k}\ .
%\lb{zapod4}
%\ee
$$
At the end we will obtain for the product $c^{(2n)}\alpha$ an expression of the form
$c^{(2n)}\alpha'$, where the element $\alpha'$ belongs to
${\cal W}(\sigma_{i+1},\dots ,\sigma_j)$, as stated.
\hfill$\blacksquare$

\medskip
%The identities in the lemma below have several versions obtained by an application of the automorphisms (\ref{homS-A}) and %(\ref{innalis}) and the antiautomorphism (\ref{antvs}). For an identity of each type we
%present one version.

\begin{lem}\label{fuproco1}
Relations  (\ref{bmw2b}) and (\ref{bmw5a}) involving the elements $\kappa_i$ have the following analogues for the higher contractors:
\ba
\lb{fup3}
c^{(2i)}\,\sigma_{2i}\, c^{(2i)}&=&\eta^{-1}\mu^{-1}c^{(2i)}\ ,
\\[2pt]
\lb{fup1}
c^{(2i)}\,\kappa_{2i}\,c^{(2i)} &=&\eta^{-1}c^{(2i)}\, .
\ea
\end{lem}

\nin{\bf Proof.~}
We prove the identity (\ref{fup1}) by induction on $i$ (the base of
induction, $i=1$, is the relation (\ref{bmw5a}) itself):
$$
%\be
\begin{array}{l}
c^{(2i+2)}\kappa_{2i+2}c^{(2i+2)}=c^{(2i)\uparrow 1} \kappa_{2i+1}\kappa_{1}
c^{(2i)\uparrow 1}\kappa_{2i+2}c^{(2i+2)}=
c^{(2i)\uparrow 1}\kappa_{2i+1}\kappa_1\kappa_{2i+2}c^{(2i+2)}
\\[1em]
\ \ \ \ \ \ \ \ \ \
=c^{(2i)\uparrow 1}\kappa_{2i+1}\kappa_{2i+2}\kappa_{2i+1}c^{(2i+2)}
=c^{(2i)\uparrow 1}\kappa_{2i+1}c^{(2i+2)}=\eta^{-1}c^{(2i+2)}\ .
\end{array}
%\lb{fupdo1}
%\ee
$$
In the first equality we used the definition (\ref{kappa-i}); in the
second equality we used the property (\ref{idemp-c1}); in the third equality we moved the element $\kappa_1$ rightwards to the contractor
$c^{(2i+2)}$ and used the property (\ref{idemp-c2}); in the fourth
equality we used the relation (\ref{bmw5a}); the fifth equality
is the induction assumption.

\medskip
The identity (\ref{fup3}) is proved again by induction on $i$ (the base
of induction, $i=1$, is now the relation (\ref{bmw2b})\,):
$$
%\be
\begin{array}{l}
c^{(2i+2)}\sigma_{2j+2}c^{(2i+2)}=c^{(2i)\uparrow 1} \kappa_{2i+1}\kappa_{1}
c^{(2i)\uparrow 1}\sigma_{2i+2}c^{(2i+2)}=
c^{(2i)\uparrow 1}\kappa_{2i+1}\kappa_1\sigma_{2i+2}c^{(2i+2)}
\\[1em]
\ \ \ \ \ \ \ \ \ \
=c^{(2i)\uparrow 1}\kappa_{2i+1}\sigma_{2i+2}\kappa_{2i+1}c^{(2i+2)}
=\mu^{-1}c^{(2i)\uparrow 1}\kappa_{2i+1}c^{(2i+2)}=\mu^{-1}
\eta^{-1}c^{(2i+2)}\ .
\end{array}
%\lb{fupdo3}
%\ee
$$
In the first equality we used the definition (\ref{kappa-i}); in the
second equality we used the property (\ref{idemp-c1}); in the third equality we moved the element $\kappa_1$ rightwards to the contractor
$c^{(2i+2)}$ and used the property (\ref{idemp-c2}); in the fourth
equality we used the relation (\ref{bmw2b}); the fifth equality
is the identity (\ref{fup1}).

\smallskip
The proof is finished.\hfill$\blacksquare$

\begin{prop}\lb{primco} The contractor $c^{(2n)}$ is a primitive idempotent in the algebra
${\cal W}_{2n}(q,\mu )$ and in the algebra ${\cal W}_{2n+1}(q,\mu )$.
\end{prop}

\noindent {\bf Proof.}~ To prove both statements about the primitivity, one has to check
that a combination $c^{(2n)}\alpha^{(2n+1)}c^{(2n)}$ is proportional to the contractor
$c^{(2n)}$ for an arbitrary element $\alpha^{(2n+1)}$ from the algebra
${\cal W}_{2n+1}(q,\mu)$.

\smallskip
Let $\alpha$ be an arbitrary element from the algebra ${\cal W}(\sigma_{1},\dots ,\sigma_j)$,
where $j\geq 2n+1$. Due to lemma \ref{zapo}, we have
$c^{(2n)}\alpha=c^{(2n)}\beta$ with
$\beta\in {\cal W}(\sigma_{n+1},\dots ,\sigma_{j})$.

\smallskip
Let $i$ ($i>0$) be such that
$\beta\in {\cal W}(\sigma_{n+i},\sigma_{n+i+1},\dots ,\sigma_{j})$
and $\beta\notin {\cal W}(\sigma_{n+i+1},\dots ,\sigma_{j})$. We shall demonstrate that
there exists an element $\bar{\beta}\in {\cal W}(\sigma_{n+i+1},\dots ,\sigma_{j})$
for which
$$
%\be
c^{(2n)}\beta c^{(2n)}=c^{(2n)}\bar{\beta} c^{(2n)}\ .
%\lb{primcod1}
%\ee
$$
Given this statement, the proof follows by induction on $i$.

\smallskip
The element $\beta$ is a linear combination of elements of the form $xu_{n+i}y$, where
the elements $x$ and $y$ belong to ${\cal W}(\sigma_{n+i+1},\dots ,\sigma_{j})$ and
$u_i$ is equal to $\sigma_{n+i}$ or $\kappa_{n+i}$. We have
$$
%\be
c^{(2n)}xu_{n+i}yc^{(2n)}=c^{(2n)}xc^{(2i)\uparrow n-i}u_{n+i}c^{(2i)\uparrow n-i}yc^{(2n)}
\sim c^{(2n)}xc^{(2i)\uparrow n-i}yc^{(2n)}=c^{(2n)}xyc^{(2n)}\ .
%\lb{primcod2}
%\ee
$$
In the first equality we used the relations (\ref{idemp-c1}); the proportionality follows from
the relations (\ref{fup1}) and (\ref{fup3}). Then we used again the relations (\ref{idemp-c1})
to absorb the contractor $c^{(2i)\uparrow n-i}$ into $c^{(2n)}$.

%\smallskip
The proof is finished.
\hfill$\blacksquare$

\section{Further properties of contractors}\label{fuproco}

The relations, involving the elements $\kappa_i$, for the generators of the BMW algebras have analogues for the higher contractors.
Two examples of such relations are proved in lemma \ref{fuproco1}. In proposition \ref{fuproco2} we prove further analogues.

\medskip
The identities in the lemma below have several versions obtained by an application of the automorphisms (\ref{homS-A}) and (\ref{innalis})
and the antiautomorphism (\ref{antvs}). For an identity of each type we
present one version.

\begin{prop}\label{fuproco2}
Another analogue of the identity (\ref{bmw5a}):
\be
\kappa_{2j} c^{(2j)}\kappa_{2j}=\eta^{-1}\kappa_{2j}
c^{(2j-2)\uparrow 1}\ .\lb{fup2}\ee

\smallskip
More general than (\ref{fup3}) analogues of the identity (\ref{bmw2b}):
\be c^{(2j)}\sigma_{j+k}\sigma_{j+k+1}\dots\sigma_{2j}c^{(2j)}=
(\eta^{-1}\mu^{-1})^{j+1-k}c^{(2j)}\quad {\mathrm{for}}\ \ 0<k\leq j
\label{fup4}\ee
and
\be c^{(2j)}\sigma_{j-k}\sigma_{j-k+1}\dots\sigma_{2j}c^{(2j)}=
\eta^{-j}(\mu^{-1})^{j-1-k}c^{(2j)}\quad {\mathrm{for}}\ \ 0\leq k<j\ .
\label{fup5}\ee

\smallskip
An analogue of the identities (\ref{bmw3}):
\be c^{(2j)}c^{(2j)\uparrow 1}=\eta^{-j} c^{(2j)}
\sigma_{2j}^{-1}\sigma_{2j-1}^{-1}\dots\sigma_1^{-1}\ .\label{fup6}\ee

\smallskip
An analogue of the identity (\ref{bmw7}):
\be \sigma_{j}'\sigma_{j-1}'\dots\sigma_1' c^{(2j)\uparrow 1}
\sigma_1'\dots\sigma_{j-1}'\sigma_{j}'=\sigma_{j+1}'\sigma_{j+2}'
\dots\sigma_{2j}'c^{(2j)}\sigma_{2j}'\dots\sigma_{j+2}'\sigma_{j+1}'
\ .\lb{fup7}\ee

\smallskip
Another analogue of the identity (\ref{bmw5a}):
\be c^{(2j)\uparrow 1}c^{(2j)}c^{(2j)\uparrow 1}=\eta^{-2j}c^{(2j)\uparrow 1}\ .\lb{fup8}\ee

\smallskip
An analogue of the identity (\ref{idemp-c3}):
\be c^{(2j)}\tau^{(2k)\uparrow j-k}=\mu^k c^{(2j)}\quad {\mathrm{for}}\ \ k\leq j\ .\lb{fup9}\ee
where the elements $\tau^{(i)}$ are defined in eq.(\ref{tau-n}).
\end{prop}

\nin{\bf Proof.~}
The identity (\ref{fup2}) is proved by induction on $j$
(the base of
induction, $j=1$, is the relation (\ref{bmw5a})\,):
$$
\begin{array}{l}
\kappa_{2j+2}c^{(2j+2)}\kappa_{2j+2}=\kappa_{2j+2}c^{(2j)\uparrow 1}
\kappa_{2j+1}\kappa_{1}c^{(2j)\uparrow 1}\kappa_{2j+2}=
c^{(2j)\uparrow 1}\kappa_{2j+2}\kappa_{2j+1}\kappa_{2j+2}\kappa_{1}
c^{(2j)\uparrow 1}\\[1em]
\ \ \ \ \ \ \ \ \ \
=c^{(2j)\uparrow 1}\kappa_{2j+2}\kappa_{1}c^{(2j)\uparrow 1}=
\eta^{-1}\kappa_{2j+2}c^{(2j)\uparrow 1}\ .
\end{array}
%\lb{fupdo2}
$$
In the first equality we used the definition (\ref{kappa-i}); in the
second equality we formed the combination
$\kappa_{2j+2}\kappa_{2j+1}\kappa_{2j+2}$; in the third equality we used
the relation (\ref{bmw5a}); the fourth equality
is the induction assumption.

\medskip
The identity (\ref{fup4}) is proved by induction on $k$ down; the base of induction, when $k=j$, is the identity (\ref{fup3}).
$$
\begin{array}{l}
c^{(2j)}\sigma_{j+k}\sigma_{j+k+1}\dots\sigma_{2j}c^{(2j)}=
c^{(2j)}c^{(2k)\uparrow j-k}\sigma_{j+k}\sigma_{j+k+1}\dots\sigma_{2j}
c^{(2k)\uparrow j-k}c^{(2j)}\\[1em]
\ \ \ \ \ \ \ \ \ \
=c^{(2j)}c^{(2k)\uparrow j-k}\sigma_{j+k}
c^{(2k)\uparrow j-k}\sigma_{j+k+1}\dots\sigma_{2j}c^{(2j)}\\[1em]
\ \ \ \ \ \ \ \ \ \ \
=\eta^{-1}\mu^{-1}c^{(2j)}c^{(2k)\uparrow j-k}\sigma_{j+k+1}
\dots\sigma_{2j}c^{(2j)}\\[1em]
\ \ \ \ \ \ \ \ \ \ \
=\eta^{-1}\mu^{-1}c^{(2j)}\sigma_{j+k+1}
\dots\sigma_{2j}c^{(2j)}
=(\eta^{-1}\mu^{-1})^{j+1-k}c^{(2j)}\ .\end{array}
%\lb{fupdo4}
$$
In the first equality we used the property (\ref{idemp-c1}); in the
second equality we formed the combination $c^{(2k)\uparrow j-k}\sigma_{j+k}
c^{(2k)\uparrow j-k}$; in the third equality we used the identity
(\ref{fup3}); in the fourth equality we used again the property
(\ref{idemp-c1}); the fifth equality is the induction assumption.

\medskip
The identity (\ref{fup5}) is proved by induction on $k$. We
have $c^{(2j)}\sigma_j=\mu c^{(2j)}$ by the relation (\ref{idemp-c3}),
so the identity (\ref{fup5}) with $k=0$ follows from the identity
(\ref{fup4}) with $k=1$.

\smallskip
Next, we have, for $i<j$:
\be \begin{array}{l}
c^{(2j)}\underline{\sigma}_i \sigma_{i+1}\dots\sigma_{2j}c^{(2j)}=
c^{(2j)}\underline{\sigma}_{2j-i} (\sigma_{i+1}\dots\sigma_{2j})c^{(2j)}\\[1em]
\ \ \ \ \ \ \ \ \ \ \
=c^{(2j)}(\sigma_{i+1}\dots\sigma_{2j})\underline{\sigma}_{2j-i-1}c^{(2j)}=
c^{(2j)}(\sigma_{i+1}\dots\sigma_{2j})\underline{\sigma}_{i+1}c^{(2j)}\ .
\end{array}\lb{fupdo5}\ee
Here we used the property (\ref{idemp-c2}) and the defining relation (\ref{braid}).

The last expression in eq.(\ref{fupdo5}) can be rewritten in a form
$$
c^{(2j)}\underline{\sigma}_{i+2}(\sigma_{i+1}\dots\sigma_{2j})c^{(2j)}
\ ,
%\lb{fupdo6}
$$
again by the braid relation (\ref{braid}).

\smallskip
If $i+2$ is still smaller than $j$, we continue in the same manner:
\be \begin{array}{l}
c^{(2j)}\underline{\sigma}_{i+2} (\sigma_{i+1}\dots\sigma_{2j})c^{(2j)}=
c^{(2j)}\underline{\sigma}_{2j-i-2} (\sigma_{i+1}\dots\sigma_{2j})c^{(2j)}\\[1em]
\ \ \ \ \ \ \ \ \ \ \
=c^{(2j)}(\sigma_{i+1}\dots\sigma_{2j})\underline{\sigma}_{2j-i-3}c^{(2j)}=
c^{(2j)}(\sigma_{i+1}\dots\sigma_{2j})\underline{\sigma}_{i+3}c^{(2j)}\
\end{array}\lb{fupdo7}\ee
and the last expression in eq.(\ref{fupdo7}) can again be rewritten in a
form
$$
c^{(2j)}\underline{\sigma}_{i+4}(\sigma_{i+1}\dots\sigma_{2j})c^{(2j)}
\ .
%\lb{fupdo8}
$$
We repeat this process till the moment when the index of the underlined
$\sigma$ becomes equal to $j$. Then we use the property (\ref{idemp-c3})
and conclude
$$
c^{(2j)}\sigma_i \sigma_{i+1}\dots\sigma_{2j}c^{(2j)}
=\mu c^{(2j)}\sigma_{i+1}\dots\sigma_{2j}c^{(2j)}\ ,
%\lb{fupdo9}
$$
which, due to the induction assumption, finishes the proof of the
identity (\ref{fup5}).

\medskip
The proof of the identity (\ref{fup6}) consists of a calculation
$$
c^{(2j)}c^{(2j)\uparrow 1}=c^{(2j)} \sigma_1\sigma_2\dots
\sigma_{2j}c^{(2j)}\sigma_{2j}^{-1}\dots\sigma_2^{-1}\sigma_1^{-1}
=\eta^{-j}c^{(2j)}c^{(2j)}\sigma_{2j}^{-1}\dots\sigma_2^{-1}\sigma_1^{-1}
\ .
%\lb{fupdo10}
$$
The first equality here is valid due to the defining relations
(\ref{braid}); in the second equality we used the identity (\ref{fup5})
with $k=j-1$.

\smallskip
Using a combination of the isomorphisms (\ref{innalis}) and (\ref{homiota'}), we can rewrite the
identity (\ref{fup6}) in forms
\be c^{(2j)\uparrow 1}c^{(2j)}=\eta^{-j} c^{(2j)\uparrow 1}
\sigma_{1}\sigma_{2}\dots\sigma_{2j}\ ,\label{fupdo15}\ee
\be c^{(2j)}c^{(2j)\uparrow 1}
=\eta^{-j} c^{(2j)}\sigma_{2j}\dots\sigma_{2}\sigma_{1}\label{fupdo16}\ee
and
\be c^{(2j)\uparrow 1}c^{(2j)}=\eta^{-j} c^{(2j)\uparrow 1}
\sigma_{1}^{-1}\sigma_{2}^{-1}\dots\sigma_{2j}^{-1}\ .\label{fupdo17}\ee

\medskip
We now turn to the proof of the identity (\ref{fup7}). First,
we prove by induction on $i$ the following identity:
\be \sigma_1'(\kappa_2\kappa_3\dots \kappa_{j+1})\sigma_1'\sigma_2'
\dots\sigma_j' =\sigma_2'\sigma_3'\dots\sigma_{j+1}'
(\kappa_1\kappa_2\dots \kappa_{j})\sigma_{j+1}'\ .\lb{fupdo11}\ee
The base of induction ($j=1$) is the identity (\ref{bmw7}).
The induction step is
$$
\begin{array}{l}
\sigma_1'(\kappa_2\kappa_3\dots \kappa_{j+2})\sigma_1'\sigma_2'
\dots\sigma_{j+1}'= \sigma_1'(\kappa_2\kappa_3\dots \kappa_{j+1})
(\sigma_1'\sigma_2'\dots\sigma_j')\kappa_{j+2}\sigma_{j+1}'\\[1em]
\ \ \ \ \ \ \ \ \ \ \
=\sigma_2'\sigma_3'\dots\sigma_{j+1}'(\kappa_1\kappa_2\dots \kappa_{j})
\sigma_{j+1}'\kappa_{j+2}\sigma_{j+1}'=
\sigma_2'\sigma_3'\dots\sigma_{j+1}'(\kappa_1\kappa_2\dots \kappa_{j})
\sigma_{j+2}'\kappa_{j+1}\sigma_{j+2}'\\[1em]
\ \ \ \ \ \ \ \ \ \ \
=\sigma_2'\sigma_3'\dots\sigma_{j+2}'(\kappa_1\kappa_2\dots \kappa_{j+1})
\sigma_{j+2}'\ ,\end{array}
%\lb{fupdo12}
$$
where we used the identity (\ref{bmw7}) in the third equality.

\smallskip
The image of the identity (\ref{fupdo11}) under the antiautomorphism
(\ref{antvs}) reads
\be \sigma_{j}'\sigma_{j-1}'\dots\sigma_{1}'
(\kappa_{j+1}\kappa_{j}\dots\kappa_{2})\sigma_{1}' =
\sigma_{j+1}' (\kappa_{j}\kappa_{j-1}\dots\kappa_{1})
\sigma_{j+1}'\sigma_{j}'\dots\sigma_{2}'\ .\lb{fupdo13}\ee

\smallskip
The proof of the identity (\ref{fup7}) is again by induction on $j$
(the base of induction is the identity (\ref{bmw7})\,):
$$
\begin{array}{l}
(\sigma_{j+1}'\sigma_{j}'\dots\sigma_1') c^{(2j+2)\uparrow 1}
(\sigma_1' \dots\sigma_{j}'\sigma_{j+1}')\\[1em]
\ \ \ \ \ \ \ \ \ \ \
=\eta^{-1}
(\sigma_{j+1}'\dots\sigma_1') c^{(2j)\uparrow 2}
(\kappa_{2j+2}\dots\kappa_{j+3})
(\kappa_2\dots\kappa_{j+2})
(\sigma_1' \dots\sigma_{j+1}')\\[1em]
\ \ \ \ \ \ \ \ \ \ \
=\eta^{-1}
(\sigma_{j+1}'\dots\sigma_2') c^{(2j)\uparrow 2}
(\kappa_{2j+2}\dots\kappa_{j+3})\underline{\sigma_1'
(\kappa_2\dots\kappa_{j+2})
(\sigma_1' \dots\sigma_{j+1}')}\\[1em]
\ \ \ \ \ \ \ \ \ \ \
=\eta^{-1}(\sigma_{j+1}'\dots\sigma_{2}') c^{(2j)\uparrow 2}
(\kappa_{2j+2}\dots\kappa_{j+3})
(\sigma_2' \dots\sigma_{j+2}')
(\kappa_1\dots\kappa_{j+1})\sigma_{j+2}'\\[1em]
\ \ \ \ \ \ \ \ \ \ \
=\eta^{-1}\underline{(\sigma_{j+1}'\dots\sigma_{2}') c^{(2j)\uparrow 2}
(\sigma_2' \dots\sigma_{j+1}')}
(\kappa_{2j+2}\dots\kappa_{j+3})\sigma_{j+2}'
(\kappa_1\dots\kappa_{j+1})\sigma_{j+2}'
\\[1em]
\ \ \ \ \ \ \ \ \ \ \
=\eta^{-1}(\sigma_{j+2}'\dots\sigma_{2j+1}')
c^{(2j)\uparrow 1}\underline{(\sigma_{2j+1}'\dots\sigma_{j+2}')
(\kappa_{2j+2}\dots\kappa_{j+3})
\sigma_{j+2}'}(\kappa_1\dots\kappa_{j+1})\sigma_{j+2}'\\[1em]
\ \ \ \ \ \ \ \ \ \ \
=\eta^{-1}(\sigma_{j+2}'\dots\sigma_{2j+1}')
 c^{(2j)\uparrow 1}
\sigma_{2j+2}' (\kappa_{2j+1}\dots\kappa_{j+2})
(\sigma_{2j+2}'\dots\sigma_{j+3}')
(\kappa_1\dots\kappa_{j+1})\sigma_{j+2}'\\[1em]
\ \ \ \ \ \ \ \ \ \ \
=\eta^{-1}(\sigma_{j+2}'\dots\sigma_{2j+2}')
 c^{(2j)\uparrow 1}
 (\kappa_{2j+1}\dots\kappa_{j+2})(\kappa_1\dots\kappa_{j+1})
 (\sigma_{2j+2}'\dots\sigma_{j+2}')
\\[1em]
\ \ \ \ \ \ \ \ \ \ \
=(\sigma_{j+2}'\dots\sigma_{2j+2}')
c^{(2j+2)}\sigma_{2j+2}'\dots\sigma_{j+2}'\ .
\end{array}
%\lb{fupdo14}
$$
Here in the first equality we used the expression (\ref{kappa-i2}) for
the contractor; in the second equality we moved the element $\sigma_1'$
rightwards to the string $(\kappa_2\dots\kappa_{j+2})$; in the third
equality we transformed the underlined expression using the identity
(\ref{fupdo11}); in the fourth equality we moved the string
$(\sigma_2'\dots\sigma_{j+1}')$ leftwards to the contractor $c^{(2j)\uparrow 2}$; in the fifth equality we used the induction
assumption to transform the underlined expression; in the sixth
equality we transformed the underlined expression using the shift
${}^{\uparrow j+1}$ of the identity (\ref{fupdo13}); in the seventh
equality we rearranged terms and then used again the expression
(\ref{kappa-i2}) for the contractor in the eighth equality.

\smallskip
The following calculation establishes the identity (\ref{fup8}):
$$
c^{(2j)\uparrow 1}c^{(2j)}c^{(2j)\uparrow 1}=\eta^{-j} c^{(2j)\uparrow 1}
c^{(2j)}\sigma_{2j}\dots\sigma_{2}\sigma_{1}=\eta^{-2j} c^{(2j)\uparrow 1}
\ .
%\label{fupdo19}
$$
Here in the first equality we used the relation (\ref{fupdo16}) while in the second one we
used the relation (\ref{fupdo17}).

\medskip
To prove the identity (\ref{fup9}), it is enough to prove its particular case
\be c^{(2j)}\tau^{(2j)}=\mu^j c^{(2j)}\lb{fupdo20}\ee
since the element $c^{(2j)}$ is divisible by the element $c^{(2j-2k)\uparrow k}$ due to the
relations (\ref{idemp-c1}).

\smallskip
We shall need two identities. The first one is
\be c^{(2j+2)}\sigma_{1}\sigma_{2}\dots\sigma_{2j}=c^{(2j+2)}c^{(2j)\uparrow 1}
\sigma_{1}\sigma_{2}\dots\sigma_{2j}=\eta^j c^{(2j+2)}c^{(2j)}\ .\label{fupdo18}\ee
In the first equality we used the relations (\ref{idemp-c1}); in the second equality we used
the relations (\ref{fupdo15}) and again (\ref{idemp-c1}).

\smallskip
Here is the second identity:
\be\begin{array}{l}
c^{(2j+2)}c^{(2j)}\sigma_{2j+1}=c^{(2j+2)}\sigma_{1}c^{(2j)}
=c^{(2j+2)}\sigma_{2j-1}c^{(2j)}
=c^{(2j+2)}\sigma_{3}c^{(2j)}\\[1em]
\ \ \ \ \ \ \ \ \ \ \ \ \ \ \ \ \ \ \ \ \ \ \, =\dots =\mu c^{(2j+2)}c^{(2j)}\ .
   \end{array}\lb{fupdo21}\ee
In the first equality we moved the element $\sigma_{2j+1}$ leftwards through the
contractor $c^{(2j)}$ and then we replaced the combination $c^{(2j+2)}\sigma_{2j+1}$ by
$c^{(2j+2)}\sigma_{1}$ due to the relation (\ref{idemp-c2}); repeatedly using the
relation (\ref{idemp-c2}), we replaced the combination $\sigma_{1}c^{(2j)}$ by
$\sigma_{2j-1}c^{(2j)}$,
then $c^{(2j+2)}\sigma_{2j-1}$ by $c^{(2j+2)}\sigma_{3}$ {\it etc.} The index of the
element $\sigma$ jumps by 2; at one moment it becomes equal to either $j$ or $j+1$
and we use then the relation (\ref{idemp-c3}).

\smallskip
We now prove the relation (\ref{fupdo20}) by induction on $j$ (the base of induction,
$j=1$, is the relation (\ref{bmw2a})\,):
$$
\begin{array}{l}
c^{(2j+2)}\tau^{(2j+2)}=c^{(2j+2)}(\sigma_1\dots\sigma_{2j+1})\tau^{(2j+1)}
=\eta^j c^{(2j+2)}c^{(2j)}\sigma_{2j+1}\tau^{(2j+1)}\\[1em]
\ \ \ \ \
=\mu\eta^j c^{(2j+2)}c^{(2j)}\tau^{(2j+1)}
=\mu\eta^j c^{(2j+2)}c^{(2j)}\tau^{(2j)}(\sigma_{2j}\dots\sigma_1)\\[1em]
\ \ \ \ \
=\mu^{j+1}\eta^j c^{(2j+2)}c^{(2j)}(\sigma_{2j}\dots\sigma_1)
=\mu^{j+1}\eta^{2j} c^{(2j+2)}c^{(2j)}c^{(2j)\uparrow 1}\\[1em]
\ \ \ \ \
=\mu^{j+1}\eta^{2j} c^{(2j+2)}c^{(2j)\uparrow 1}c^{(2j)}c^{(2j)\uparrow 1}
=\mu^{j+1}c^{(2j+2)}\ .
\end{array}
%\lb{fupdo22}
$$
In the first equality we used the iterative definition of the
elements $\tau^{(i)}$ (it is different but equivalent to the one given in eq.(\ref{innalis})\,);
in the second equality we used the relation (\ref{fupdo18}); in the third equality we used the
relation (\ref{fupdo21}); in the fourth equality we used again the iterative definition
of the elements $\tau^{(i)}$; the fifth equality is the induction assumption; in the sixth
equality we used the relation (\ref{fupdo16}); in the seventh equality we used the relations
(\ref{idemp-c1}); finally, in the eighth equality we used the relation (\ref{fup8}).

\smallskip
The proof is finished.
\hfill$\blacksquare$

\begin{rem}{\rm We have also
$$
c^{(2j+2)}\tau^{(2j+1)}=c^{(2j+2)}(\sigma_1\dots\sigma_{2j})\tau^{(2j)}
=\eta^j c^{(2j+2)}c^{(2j)}\tau^{(2j)}=(\eta\mu )^j c^{(2j+2)}c^{(2j)}\ .
%\lb{fupdo23}
$$
In the first equality we used the iterative definition of the
elements $\tau^{(i)}$; in the second equality we used the relation (\ref{fupdo18});
in the third equality we used the identity (\ref{fup9}).}\end{rem}

%\appendix
\section{On twists in quasitriangular Hopf algebras}\label{append1}

\def\theequation{\thesection.\arabic{equation}}
%\makeatletter\@addtoreset{equation}{section}\makeatother

Here we shall  discuss universal (i.e., quasi-triangular Hopf algebraic) counterparts
of relations from subsections \ref{subsec3.2}, \ref{subsec3.2a}, especially from
proposition \ref{proposition3.6}:
we shall see, in  item {\bf  8} of the appendix, that these relations have a quite transparent
meaning, they reflect the properties of the twisted universal R-matrix.

\medskip
We do not give an introduction to the theory of quasitriangular Hopf algebras assuming that the reader
has some basic knowledge on the subject (see, e.g., \cite{CP}, the chapter 4).

\subsection{Generalities}

\paragraph{{\bf  1.}} Let ${\cal A}$ be a Hopf algebra; $m, \D ,\e$ and $S$ denote the multiplication, comultiplication,
counit and antipode, respectively.

\medskip
Assume that ${\cal A}$ is quasitriangular with a universal R-matrix $\cR =a\ot b$ (this is a symbolic
notation, instead of $\sum_i a_i\ot b_i$). One has $(S\ot S)\cR =\cR$. The universal R-matrix
$\cR$ is invertible, its inverse is related to $\cR$ by formulas~ $\cR^{-1}=S(a)\ot b$ ~or~
$(\id\ot S)(\cR^{-1})=\cR$.

\smallskip
{}For elements in ${\cal A}\ot {\cal A}$, the `skew' product $\odot$ is defined as the product in
${\cal A}^{\mathrm{op}}\ot {\cal A}$, where ${\cal A}^{\mathrm{op}}$ denotes the algebra with the opposite multiplication. In
other words, the skew product of two elements, $x\ot y$ and $\ti{x}\ot\ti{y}$ is
$(x\ot y)\odot (\ti{x}\ot\ti{y})=\ti{x}x\ot y\ti{y}$.
{}For a skew invertible element ${\cal X}\in {\cal A}\ot {\cal A}$, we shall denote its skew inverse by $\psi_{\cal X}$.
The universal R-matrix $\cR$ has a skew inverse, $\psi_\cR =a\ot S(b)$. The element $\psi_\cR$
is invertible, $(\psi_\cR)^{-1}=a\ot S^2(b)$. The element $\cR^{-1}$ is skew invertible as well,
its skew inverse is $\psi_{(\cR^{-1})} =S^2(a)\ot b$. All these formulas are present in \cite{D3}.
We shall see below that there are similar formulas for the twisting element $\cF$. However,
the properties of the twisting element $\cF$ and of the universal R-matrix $\cR$ are different,
for instance, the square of the antipode is given by $S^2(x)=u_{_\cR}\, x\, (u_{_\cR} )^{-1}$, where
$u_{_\cR} =S(b)a$, but there is no analogue of such formula for $\cF$. Because of this difference, we
felt obliged to give some proofs of the relations for $\cF$.

\medskip
Let $\rho$ be a representation of the algebra ${\cal A}$ in a vector space $V$. For an element ${\cal X}\in {\cal A}\ot {\cal A}$,
denote by $\hat{\rho}({\cal X})\in {\rm End}(V^{\otimes 2})$ an operator $\hat{\rho}({\cal X})=P\cdot (\rho\ot\rho)({\cal X})$
(recall that $P$ is the permutation operator). The skew product $\odot$ translates into the following
product $\hat{\odot}$ for elements of ${\rm End}(V^{\otimes 2})$:
\be
%\lb{skewprod}
(X\hat{\odot}Y)_{13}:=\tr_{(2)}(X_{12}Y_{23})\ .
\ee
In other words, if ${\cal X}\odot {\cal Y}={\cal Z}$ then
$\hat{\rho}({\cal X})\;\hat{\odot}\;\hat{\rho}({\cal Y})=\hat{\rho}({\cal Z})$. For an operator
$X\in {\rm End}(V^{\otimes 2})$, its skew inverse $\Psi_X$, in the sense explained in
subsection \ref{subsec3.1}, is presicely the inverse with respect to the product $\hat{\odot}$.

\paragraph{{\bf  2.}} The following lemma is well known (see, e.g., \cite{CP},
the chapter 4, and references therein).

\begin{lem}\lb{lemma3.2.0.1}
Consider an invertible element $\cF =\a\ot\b\in {\cal A}\ot {\cal A}$ (we use the symbolic notation,
$\a\ot\b =\sum_i \a_i\ot\b_i$, like for the universal R-matrix) and
let $\cF^{-1}=\g\ot\d$. Assume that the element $\cF$ satisfies
\be \cF_{12}\; (\D\ot\id)(\cF )=\cF_{23}\; (\id\ot\D )(\cF )\ .\lb{uu1}\ee
Assume also that
\be \e (\a )\,\b =\a\,\e(\b )=1\ .\lb{uu2'}\ee
Then an element $\vf =\a\, S(\b )$ is invertible, its inverse is
\be (\vf )^{-1}=S(\g )\,\d\ .\lb{vfin}\ee
One also has
\be S(\a )\, (\vf )^{-1}\,\b =1\ \ \ {\mathrm{and}}\ \ \ \g\,\vf\, S(\d )=1\ .\lb{uu14}\ee
Twisting the coproduct by the element $\cF$,
\be \D_\cF (a)=\cF\;\D (a)\;\cF^{-1}\ ,\ee
one obtains another quasitriangular structure on ${\cal A}$ with
\be \cR_\cF =\cF_{21}\;\cR\;\cF^{-1}\ \lb{anor}\ee
and
\be\ S_\cF (a)=\vf\; S(a)\; (\vf )^{-1}\ \lb{uu17}\ee
(the counit does not change).
\end{lem}

An element $\cF$, satisfying conditions (\ref{uu1}) and (\ref{uu2'}) is called {\em twisting} element.
We shall denote by ${\cal A}_\cF$ the resulting `twisted' quasitriangular Hopf algebra.

\begin{rem}{\rm On the representation level, the formula (\ref{anor}) transforms (compare with
eq.(\ref{R_f})$\,$) into $\hat{\rho}(\cR_\cF )=
\hat{\rho}(\cF )_{21}\hat{\rho}(\cR )_{21}\hat{\rho}(\cF )_{21}^{-1}$.
Below, when we talk about matrix counterparts of universal formulas, one should keep in mind this
difference in conventions.\lb{rrr1}}\end{rem}

\paragraph{{\bf  3.}} Assume, in addition to eq.(\ref{uu1}), that
\be (\D \ot\id )\; (\cF )=\cF_{13}\;\cF_{23}\ \lb{uu18}\ee
and
\be (\id\ot\D )\; (\cF )=\cF_{13}\;\cF_{12}\ .\lb{uu19}\ee

\smallskip
Now the conditions (\ref{uu2'}) follow from the relations (\ref{uu18}) and (\ref{uu19}) and the invertibility of
the twisting element $\cF$: applying $\e\ot\id\ot\id$ to the relation (\ref{uu18}), we find
$(\e\ot\id ) (\cF)=1$; applying $\id\ot\id\ot\e$ to the relation (\ref{uu19}), we find
$(\id\ot\e ) (\cF)=1$.

\medskip
Since $\D^{{\mathrm{op}}}(x)\cR =\cR\D (x)$ for any element $x\in {\cal A}$ (where $\D^{{\mathrm{op}}}$ is the
opposite comultiplication), it follows from the relation (\ref{uu18}) that
\be \cR_{12}\;\cF_{13}\;\cF_{23}=\cF_{23}\;\cF_{13}\;\cR_{12}\ .\lb{uu20}\ee

\smallskip
Similarly, the relation (\ref{uu19}) implies
\be \cR_{23}\;\cF_{13}\;\cF_{12}=\cF_{12}\;\cF_{13}\;\cR_{23}\ .\lb{uu22}\ee

\smallskip
When both relations (\ref{uu18}) and (\ref{uu19}) are satisfied, the relation (\ref{uu1}) is equivalent to
the Yang--Baxter equation for the twisting element $\cF$:
\be \cF_{12}\;\cF_{13}\;\cF_{23}=\cF_{23}\;\cF_{13}\;\cF_{12}\ .\lb{uu24}\ee

\begin{rem}{\rm One also has
$$
%\be
(\D_\cF\ot\id ) (\cF_{21})=\cF_{31}\;\cF_{32}\ \ {\mathrm{and}}\ \
(\id\ot\D_\cF )(\cF_{21})=\cF_{31}\;\cF_{21}\ .
%\ee
$$
Therefore, one can twist $\D_\cF$ again, now by the element $\cF_{21}$.

\smallskip
On the matrix level, this corresponds to the second conjugation of $\hat{\rho}(\cR )$ by
$\hat{\rho}(\cF )$,~~
$$\hat{\rho}\bigl(\, (\cR_\cF )_{\cF_{21}}\,\bigr)\,
=\,\hat{\rho}(\cF )^2\; \hat{\rho}(\cR )\;\hat{\rho}(\cF )^{-2}\ .$$
\lb{rrr2}}\end{rem}

\begin{rem}{\rm The element $\cF_{21}^{-1}$ satisfies the conditions (\ref{uu1}), (\ref{uu18}) and
(\ref{uu19}) if the element $\cF$ does. Thus, one can twist the coproduct $\D$ by the element $\cF_{21}^{-1}$ as well.\lb{rrr3}}\end{rem}

\paragraph{{\bf  4.}} The conditions (\ref{uu2'}), (\ref{uu18}), (\ref{uu19})
imply the invertibility and skew-invertibility of the element $\cF$. The formulas
for its inverse and skew inverse are similar to the corresponding formulas for the universal R-matrix
$\cR$ (in particular, we reproduce the standard formulas for $\cR$ since we can take $\cF =\cR$).

\begin{lem}\lb{lemma3.2.0.2}
Assume that the conditions (\ref{uu2'}) and (\ref{uu18})
are satisfied. Then the element $\cF$ is invertible, its inverse is
\be \cF^{-1}=S(\a )\ot\b \ .\lb{uu26}\ee
Assume that the conditions (\ref{uu2'}) and (\ref{uu19})
are satisfied. Then the element $\cF$ is skew invertible, with the skew inverse
\be \psi_\cF =\a\ot S(\b )\ .\lb{uu34}\ee
Assume that the conditions (\ref{uu2'}), (\ref{uu18}) and (\ref{uu19})
are satisfied. Then
\be (S\ot S)(\cF )=\cF\ .\lb{uu29}\ee
Moreover, the element $\psi_\cF $ is invertible, its inverse is
\be (\psi_\cF )^{-1}=\a\ot S^2(\b )\ \lb{uu35}\ee
and the element $\cF^{-1}$ is skew-invertible, its skew inverse reads
\be \psi_{(\cF^{-1})} =S^2(\a)\ot \b \ .\lb{uu34'}\ee
\end{lem}

\nin {\bf Proof.~} The calculations are similar to those, from textbooks, for the universal R-matrix.
We include this proof for a completness only.

\medskip
Applications of $m_{12}\circ S_1$ and $m_{12}\circ S_2$ to the relation (\ref{uu18}) imply
the formula (\ref{uu26})
(here $m_{12}$ is the multiplication of the first and the second tensor arguments; $S_1$ is an
operation of taking the antipode of the first tensor argument, {\it etc.}).

\smallskip
Applications of $m_{23}\circ S_2$ and $m_{23}\circ S_3$ to the relation (\ref{uu19}) establish
the formula (\ref{uu34}).

\medskip
Given the formula (\ref{uu34}), the statement, that the element $\psi_\cF$ is a  left skew
inverse of the element $\cF$, reads in components:
\be \a\a'\ot S(\b')\b =1\ ,\lb{uupsi}\ee
where primes are used to distinguish different
summations terms, the expression $\a\a'\ot S(\b')\b$ stands for $\sum_{i,j}\a_i\a_j\ot S(\b_j)\b_i$.
Applying $S_1$ to this equation, we find $(S(\a' )\ot S(\b'))\cdot (S(\a )\ot\b)=1$ which means
that the element $S(\a' )\ot S(\b')$ is the left inverse of the element $S(\a )\ot\b$.
However, the latter element is, by the formula (\ref{uu26}), the inverse of $\cF$. Therefore,
the relation (\ref{uu29}) follows.

\medskip
Applying $S_2$ to the equality (\ref{uupsi}), we find that the element $\a\ot S^2(\b )$ is the
right inverse of the element $\psi_\cF$.

\medskip
Applying $S_1^2$ to the equality (\ref{uupsi}) and using the relation (\ref{uu29}), we find that
$S^2(\a )\ot \b $ is a right skew inverse of the element $\cF^{-1}$.

\medskip
We shall not repeat details for the left inverse of the element $\psi_\cF$ and the left
skew inverse of the element $\cF^{-1}$, calculations are analogous.\hfill$\blacksquare$

\begin{rem}{\rm There is a further generalization of the formulas from lemma
\ref{lemma3.2.0.2}. Start with the element $\cF$ and alternate operations `take an inverse' and
`take a skew inverse'. Then the next operation is always possible, the result is always invertible
and skew invertible. One arrives, after $n$ steps, at $S^n(\a )\ot\b$ if the first operation was
`take an inverse'; if the first operation was `take a skew inverse' then one arrives
at $\a\ot S^n(\b )$ (see \cite{D3}, section 8). \lb{rrr4}}\end{rem}

\smallskip
{\it From now on, we shall assume that the twisting element $\cF$ is invertible
and satisfies the conditions (\ref{uu1}), (\ref{uu18}) and (\ref{uu19}).}

\subsection{Counterparts of matrix relations}

\paragraph{{\bf  5.}} We turn now to the Hopf algebraic meaning of relations from  subsections \ref{subsec3.2}, \ref{subsec3.2a}.

\medskip
The square of the antipode in an almost cocommutative Hopf algebra, with a universal R-matrix
$\cR =a\ot b$, satisfies the property $S^2(x)=u_{_\cR} x(u_{_\cR})^{-1}$, where $u_{_\cR} =S(b)a$,
for any element $x\in {\cal A}$.
In a matrix representation of an algebra ${\cal A}$, the element $u_{_\cR}$ maps to the matrix
$D_{\hat{\rho}({\cal R})}$
(and the element $S(u_{_\cR} )$ maps to the matrix $C_{\hat{\rho}({\cal R})}$), so an identity (which follows from
the relation (\ref{uu29})$\,$)
$$\begin{array}{rcl}(1\ot u_{_\cR})\;\cF^{-1}\; (1\ot (u_{_\cR})^{-1})&\equiv& (1\ot u_{_\cR} )
(S(\a )\ot\b )(1\ot (u_{_\cR})^{-1}=S(\a )\ot S^2(\b )
\\[1em]
&=&\a\ot S(\b )\equiv\psi_\cF\end{array}$$
becomes one of the relations from lemma \ref{lemma3.3}. In a similar manner, one can interpret other
relations from lemma \ref{lemma3.3}.

\medskip
Such an interpretation is not, however, unique. For instance, applying $m_{12}\circ S_2$ to the relation (\ref{uu24})
and using the formula (\ref{uu26}), one finds
$$ \vf\ot 1=\a'\vf S(\a )\ot\b\b'\ ,$$
which, after an application of $S_2$, becomes, due to the formulas (\ref{uu34}) and (\ref{uu29}),
\be \vf\ot 1=\psi_\cF\; (\vf\ot 1)\;\cF\ .\lb{uu39}\ee
Similarly, applying $(\id\ot S)\circ m_{23}\circ\tau_{23}\circ S_3$ (where $\tau$ is the flip,
$\tau (x\ot y)=y\ot x$) to eq.(\ref{uu24}) and using eqs.(\ref{uu29}) and (\ref{uu35}), one finds
$$ 1\ot\vf =\a\a'\ot S(\b')\vf S^2(\b )\ ,$$
which, after an application of $S_1$, becomes, with the help of eq.(\ref{uu29}),
\be 1\ot\vf =\cF\; (1\ot\vf )\;\psi_\cF\ .\lb{uu36}\ee
In the matrix picture, the relations (\ref{uu39}) and (\ref{uu36}) are also equivalent to particular cases of the
relations from lemma \ref{lemma3.3} -- but this time we did not use the fact that the square of the
antipode is given by the conjugation by the element $u_{_\cR}$.

\smallskip
Below we shall make use of another version of the formulas (\ref{uu39}) and (\ref{uu36}).

\smallskip
Writing the formulas (\ref{uu39}) and (\ref{uu36})
as $(\vf\ot 1)\cF^{-1}=\psi_\cF (\vf\ot 1)$ and
$ \cF^{-1}(1\ot\vf) =(1\ot\vf )\psi_\cF$, respectively, and using the expressions
for $\psi_\cF$, $(\psi_\cF )^{-1}$
and $ \cF^{-1}$ from lemma \ref{lemma3.2.0.2}, we find, in components:
\be\vf S(\a ) \ot\b =\a\vf\ot S(\b ) \ \lb{uu41}\ee
and, respectively,
\be S(\a )\ot\b\vf  =\a\ot\vf S(\b )\ .\lb{uu37}\ee

Applying $S_1$ or $S_2$ to eqs.(\ref{uu41}) and (\ref{uu37}), we obtain corresponding formulas with
$\vf$ replaced by $S(\vf )$. These formulas, together with eqs.(\ref{uu41}) and (\ref{uu37}), imply
\be\begin{array}{ccc} \cF\cdot (\vf S(\vf )\ot 1)&=&(\vf S(\vf )\ot 1)\cdot \cF\ ,\\[1em]
\cF\cdot (1\ot \vf S(\vf ))&=&(1\ot \vf S(\vf ))\cdot \cF\ .\end{array}\lb{fvsv}\ee

It follows, from a compatibility of the relations (\ref{uu39}) and (\ref{uu36}) (express the element $\psi_\cF$
in terms of $\cF$ and $\vf$ in two ways), that
\be \cF_{12}\cdot (\vf\ot\vf ) =(\vf\ot\vf )\cdot \cF_{12}\ .\lb{fvvvvf}\ee

The relations (\ref{fvsv}) and (\ref{fvvvvf}) are universal analogues of the matrix equalities
(\ref{FCD}) and (\ref{FCC}) (for certain choices of the compatible pairs of the R-matrices)
from the corollary \ref{corollary3.4}.

\paragraph{{\bf  6.}} We need some more information about the element $\vf$. The
inverse to the element $\vf$ is given by the formula (\ref{vfin}); it follows from lemma
\ref{lemma3.2.0.2} that $(\vf )^{-1}=S^2(\a )\b$.

\medskip
By eq.(\ref{uu29}), one has $S(\vf )=S(\b )\a$ and, then,
$S^2(\vf )=\vf $. Since $S^2(x)=u_{_\cR} x(u_{_\cR})^{-1}$ for any element $x\in {\cal A}$, we conclude that the element
$u_{_\cR}$ commutes with the element $\vf$ and, similarly, with the element $S(\vf )$.

\medskip
Making the flip in the relations (\ref{uu41}) and (\ref{uu37}), multiplying them out
and comparing, we find that the elements $\vf$ and $S(\vf )$ commute.

\begin{rem}{\rm In fact, more is true. Applying $\id\ot S^j$ to the relation (\ref{uu41}),
we obtain $\vf\a\ot S^{j-1}(\b )=\a\vf\ot S^{j+1}(\b )$ (we used the relation (\ref{uu29}) to
rearrange the powers of the antipode).
In a similar way, applying $S^{-j}\ot\id$ to
the relation (\ref{uu37}), we obtain $\a\ot S^{j-1}(\b )\vf =\a\ot\vf S^{j+1}(\b )$. Multiplying out and comparing
the right hand sides, we find that the element $\vf$ commutes with the elements $S^k(\a )\b$
$\forall$ $k\in {\Bbb Z}$.

\smallskip
The same procedure, applied to the flipped versions of the relations (\ref{uu41}) and
(\ref{uu37}) shows
that the element $\vf$ commutes with the elements $S^k(\b )\a$ $\forall$ $k\in {\Bbb Z}$.

\smallskip
Applying the antipode to these commutativity relations, we find that the element $S(\vf )$ commutes with
the elements $S^k(\a )\b$ and $S^k(\b )\a$ $\forall$ $k\in {\Bbb Z}$ as well.}\lb{rrrq}\end{rem}

\paragraph{{\bf  7.}} We shall now establish a Hopf algebraic counterpart of the relation
(\ref{rcdm}).

\medskip
There is a closed formula for the coproduct of the element $\vf$, again similar to the
standard formula for the coproduct of the element $u_{_\cR}$.

\begin{lem}\lb{lemma3.2.0.3} One has
\be \D (\vf )=\cF_{12}^{-1}\;\cF_{21}^{-1}\cdot (\vf\ot\vf )\ .\lb{deltavf}\ee
\end{lem}

\nin {\bf Proof.~} Together, eqs.(\ref{uu18}) and (\ref{uu19}) imply
$$
%\be
(\D\ot\D )(\cF )=\cF_{14}\cF_{13}\cF_{24}\cF_{23}\ .
%\lb{ddf}
%\ee
$$
Therefore, the coproduct of $\vf$ can be written in a form
\be
\lb{vvs1}
\D (\vf )=\a_{(1)} S(\b_{(2)})\ot\a_{(2)} S(\b_{(1)} )=\a\a' S(\b\b'')\ot\a'' \vf S(\b')\
\ee
(we use the Sweedler notation for the coproduct, $\D (x)=x_{(1)}\ot
x_{(2)}$ for an element $x\in {\cal A}$).

\smallskip
Using the relation (\ref{uu37}), we continue to rewrite the expression (\ref{vvs1}):
\be
\lb{vvs2}
\D (\vf )=\a S(\a') S(\b\b'')\ot\a''\b' \vf \ .
\ee
The relation (\ref{uu24}), in a form $\cF_{13}\cF_{23}\cF_{12}^{-1}=\cF_{12}^{-1}\cF_{23}\cF_{13}$,
reads, in components,
\be
\lb{cfmff}
\a S(\a')\ot\a''\b'\ot\b\b''=S(\a )\a''\ot\b\a'\ot\b'\b''\ .
\ee
Using eq.(\ref{cfmff}), we transform the right hand side of eq.(\ref{vvs2}) to a form
$$
%\be
%\lb{vvs3}
\D (\vf )=S(\a )\a''S(\b'')S(\b')\ot\b\a' \vf =S(\a )\vf S(\b')\ot\b\a' \vf\ .
%\ee
$$
Using again eq.(\ref{uu37}), we obtain
$$
%\be
%\lb{vvs4}
\D (\vf )=S(\a )\b'\vf\ot\b S(\a' ) \vf\ ,
%\ee
$$
which, by the formula (\ref{uu26}), is a component form of the relation (\ref{deltavf}).
\hfill$\blacksquare$

\medskip
Applying the flip to the relation (\ref{deltavf}), we find
$\D^{{\mathrm{op}}} (\vf )=\cF_{21}^{-1}\;\cF_{12}^{-1}\cdot (\vf\ot\vf )$. Since
$\D^{{\mathrm{op}}} (\vf )\;\cR =\cR\;\D (\vf )$, we conclude
\be (\cR_\cF )_{\cF_{21}}\; (\vf\ot\vf )=(\vf\ot\vf )\;\cR\ .\lb{twitwi}\ee
The translation of the equality (\ref{twitwi}) into the matrix language is equivalent to the
relation (\ref{rcdm}) (see the remarks \ref{rrr1} and \ref{rrr2}).

\begin{rem}{\rm It follows from the relation (\ref{deltavf}) that
\be\D (S(\vf ))=(S(\vf )\ot S(\vf ))\cdot\cF_{12}^{-1}\;\cF_{21}^{-1}\ .\lb{svf1}\ee
The relation (\ref{fvvvvf}), together with the relations (\ref{deltavf}) and (\ref{svf1}),
implies that an element
\be\varphi :=\vf S(\vf )^{-1}\lb{deffi}\ee
is group-like, $\D (\varphi )=\varphi\ot\varphi$. Therefore,
$S(\varphi )=\varphi^{-1}=S(\vf )(\vf )^{-1}$ but $S(\varphi )=S(\vf S(\vf )^{-1})=(\vf )^{-1}S(\vf )$, which shows
again that $\vf$ commutes with $S(\vf )$.\lb{rrr6}}\end{rem}

\paragraph{{\bf  8.}} The twisted Hopf algebra ${\cal A}_\cF$ is quasitriangular, so we can write the usual identities for its universal
R-matrix $\cR_\cF =\cF_{21}\cR\cF^{-1}$. The relations from
proposition \ref{proposition3.6} are the matrix counterparts of some of these identities.

\medskip
{}For the twisted Hopf algebra ${\cal A}_\cF$, one finds, with the help of the first relation
in eq.(\ref{uu14}),
that $u_{(\cR_\cF)}=\varphi\, u_{_\cR}$, where the element $\varphi$ is defined by the formula (\ref{deffi}) (on the matrix level, this becomes one of the relations (\ref{CDtwist})$\,$). In
particular,
\be (S_\cF)^2\, (x)=\varphi\; S^2(x)\;\varphi^{-1}\ .\lb{sfsf}\ee

\medskip
({\it i}) The relation (\ref{R_f-fin}) is a consequence of, for example, the identity
\be (\id\ot S_\cF )((\cR_\cF )^{-1})=\cR_\cF\ .\lb{ee1drf}\ee
We have
\be\begin{array}{rcl} \cR_\cF &=&(\id\ot S_\cF )((\cR_\cF )^{-1})=(\id\ot S_\cF )(\cF\cR^{-1}\cF_{21}^{-1})
=(\id\ot S_\cF )(\a S(a)\b'\ot\b b S(\a'))\\[1em]
&=& \a S(a)\b'\ot \vf S^2(\a')S(b)S(\b )(\vf )^{-1}=\a a\b'\ot \vf S^2(\a')b S(\b )(\vf )^{-1}\ .
\end{array}\lb{le36-1}\ee
Here we used eq.(\ref{uu17}) and the identities from lemma \ref{lemma3.2.0.2} for $\cF$ and $\cR$.
Applying $S^2\ot S$ to eq.(\ref{uu41}), we find
\be \vf S^2(\a )\ot\b =\a\vf\ot\b\ ,\lb{vsvs1}\ee
since $S^2(\vf )=\vf$.
Using the relation (\ref{vsvs1}) and the relation (\ref{uu37}) in a form
$S(\a )\ot (\vf)^{-1}\b =\a\ot S(\b )(\vf)^{-1}$, we rewrite the last expression in eq.(\ref{le36-1}):
$$ \cR_\cF =S(\a )a\b'\ot \a'\vf b (\vf )^{-1}\b \ $$
or
\be \cR_\cF =\cF_{21}\odot \Bigl( (1\ot\vf )\cR (1\ot\vf )^{-1}\Bigr)\odot\cF^{-1}\ ,\ee
which, on the matrix level, is equivalent to the relation (\ref{R_f-fin}).

\medskip
({\it ii}) Next,
$$\begin{array}{rcl} \psi_{(\cR_\cF )}&=&(\id\ot S_\cF)(\cR_\cF)
=(\id\ot S_\cF)(\cF_{21}\cR\cF^{-1})=(\id\ot S_\cF)(\b a S(\a')\ot\a b\b')\\[1em]
&=&\b a S(\a')\ot \vf S(\b')S(b)S(\a )(\vf )^{-1}
=\b a \a'\ot \vf \b' S(b)S(\a )(\vf )^{-1}\ \end{array}$$
or
\be (1\ot\vf )^{-1}\; \psi_{(\cR_\cF )}\;(1\ot\vf )=\cF\odot\psi_\cR\odot\cF_{21}^{-1}\ ,\ee
which, on the matrix level, is equivalent to the relation (\ref{Psi_R_f}).

\medskip
({\it iii}) To obtain another formula for $\psi_{(\cR_\cF )}$, we start with the identity
$\psi_{(\cR_\cF )}=(\id\ot (S_\cF )^2)((\cR_\cF )^{-1})$, which is a direct consequence of the identities from
lemma \ref{lemma3.2.0.2}:
\be\begin{array}{rcl}
\psi_{(\cR_\cF )}&=&(\id\ot (S_\cF )^2 )( \cF\cR^{-1}\cF_{21}^{-1})
=(\id\ot (S_\cF )^2 )\Bigl( \a S(a)\b'\ot\b b S(\a')\Bigr) \\[1em]
&=&\a S(a)\b'\ot\varphi S^2(\b )S^2(b)S^3(\a')\varphi^{-1}
=\a a\b'\ot\varphi S^2(\b )S(b)S^3(\a')\varphi^{-1}\\[1em]
&=&\a a\b'\ot S(\vf )^{-1}\b\vf S(b)(\vf)^{-1}S(\a')S(\vf )\ .\end{array}\lb{drfof}\ee
Here we used the identities from lemma \ref{lemma3.2.0.2}, relations
$\a\ot\vf S^2(\b )=\a\ot\b\vf$ and
$S^3(\a )(\vf )^{-1}\ot\b =(\vf )^{-1}S(\a )\ot\b$, which follow from eqs.(\ref{uu41}) and
(\ref{uu37}), and the formula (\ref{sfsf}) for the square of the twisted antipode.

\smallskip
Eq.(\ref{drfof}) can be rewritten as
\be \bigl( 1\ot S(\vf )\bigr)\; \psi_{(\cR_\cF )}\;\bigl( 1\ot S(\vf )^{-1}\bigr)
=\cF\; (1\ot \vf)\;\psi_\cR\; (1\ot(\vf)^{-1})\;\cF_{21}^{-1}\ ,\ee
which, in the matrix picture, is equivalent to eq.(\ref{Psi_R_f-another}).

\medskip
({\it iv}) The property $(S_\cF\ot S_\cF )(\cR_\cF )=\cR_\cF$ leads to
\be (\vf\ot\vf )\;\cF^{-1}\;\cR\;\cF_{21}= \cF_{21}\;\cR\;\cF^{-1}\; (\vf\ot\vf )\ .\lb{anom}\ee
Since the twisting element $\cF$ commutes with $\vf\ot\vf$, the formula (\ref{anom}) is another
manifestation of the relation (\ref{rcdm}).

\begin{rem}{\rm We conclude this appendix with several more properties of the group-like
element $\varphi$ defined in eq.(\ref{deffi}).

\medskip
We have
\be \cR\cdot (\varphi\ot\varphi )=(\varphi\ot\varphi )\cdot\cR\ .\lb{rgrelu}\ee
To see this, apply $S\otimes S$ to
the relation (\ref{twitwi}) and then compare with the same relation (\ref{twitwi}).

\smallskip
The matrix equivalent of the relation (\ref{rgrelu}) is the relation (\ref{rgrel}).

\medskip
Recall that a quasitriangular Hopf algebra ${\cal A}$ is called a ribbon Hopf algebra if it contains
a ribbon element $r$, that is, a central element such that $r^2=u_{_\cR} S(u_{_\cR} )$ and
$\D (r)=\cR_{12}^{-1}\cR_{21}^{-1}\cdot (r\ot r)$ (see \cite{Resh}, or \cite{CP}, the chapter 4).
The twisted algebra ${\cal A}_\cF$ is a ribbon Hopf algebra if the algebra ${\cal A}$ is; for the ribbon element
of the algebra ${\cal A}_\cF$, one can choose $r_{_\cF} = \varphi r$.
}\end{rem}

\section*{Acknowledgments}

The authors express their gratitude to
%Robert Coquereaux,
Dimitry Gurevich,
%Nikolai Iorgov,
Alexei Isaev, Alexander Molev, Andrei Mudrov
%Arun Ram
and Pavel Saponov for fruitful discussions and valuable remarks. The work of the first
author (O. O.) was supported by the Program of Competitive Growth of Kazan Federal University and by the grant RFBR 17-01-00585. The work of the second author (P. P.) was partially supported
by the grant of RFBR no.19-01-00726-a, and
by the Academic Fund Program at the HSE University  (grant no.282948 for the years 2020-2022) and the Russian Academic Excellence Project `5-100'.

%\section{vspom}

\bigskip\bigskip
%\newpage
\addtocontents{toc}{\contentsline {section}{\numberline {} References}{\pageref{refer}}}


\begin{thebibliography}{999}
\label{refer}

\bigskip

\bibitem[1]{Br} Brauer, R.: {\em `On algebras which are connected with the semisimple
continuous groups'}. Ann. Math. {\bf 38} (1937), pp. 854--872.

\bibitem[2]{BCC} O'Brien, D.M., Cant, A. and Carey, A.L.: {\em `On characteristic identities for Lie algebras'}.
Ann. Inst. Henry Poincar{\' e}, A {\bf 26} (1977), pp. 405--429.

\bibitem[3]{BGr} Bracken, A.J. and Green, H.S.: {\em `Vector operators and a polynomial identity
for ${\rm SO}(n)$'}. J. Math. Phys. {\bf 12} (1971), pp. 2099--2106.

\bibitem[4]{BW} Birman, J.S. and Wenzl, H.: {\em `Braids, link polynomials and a new algebra'}. Trans.
Amer. Math. Soc. {\bf 313}, no. 1, (1989), pp. 249--273.

\bibitem[5]{C} Cherednik, I.V.: Theor. Math. Phys. {\bf 61} no.1, (1984),  pp. 977--983.

\bibitem[6]{CP} Chari, V. and Pressley, A.: {\em A guide to quantum groups}.
Cambridge University Press, Cambridge, 1994.

\bibitem[7]{D1} Drinfel'd, V.G.: {\em `Quantum groups'}. Proceedings of the International Congress of
Mathematicians, Vol. 1, (Berkeley, California, USA, 1986), pp. 798--820, Amer. Math. Soc., Providence, RI, 1987.

\bibitem[8]{D2}  Drinfel'd, V.G.: {\em `Quasi-Hopf algebras'}.
Leningrad Math. J. {\bf 1} no. 6, (1990),  pp. 1419--1457.

\bibitem[9]{D3}  Drinfel'd, V.G.: {\em 'On almost cocommutative Hopf algebras'}.
Leningrad Math. J. {\bf 1} no. 6, (1990),  pp. 321--342.

\bibitem[10]{EOW} Ewen, H., Ogievetsky, O. and Wess, J.: {\em `Quantum matrices in two dimensions'}.
Lett. Math. Phys. {\bf 22} no. 4, (1991),  pp. 297--305.

\bibitem[11]{F} Formanek, E.: {\em 'The ring of generic matrices'}. J. Algebra {\bf 258} no. 1, (2002),   pp. 310--320.

%\bibitem[G--T]{G-T} Gelfand I. M., Krob D., Lascoux A., Leclerc B., Retakh V. S. and Thibon J.,
%{\em `Noncommutative symmetric functions'}. Adv. Math. {\bf 112} (1995), no. 2, pp. 218--348.
%ArXiv: hep-th/\-9407124.

%\bibitem[GIO]{GIO} Gorbounov V., Isaev A. and Ogievetsky O., {\em BRST Operator for Quantum
%Lie Algebras: Relation to Bar Complex}. Teoret. Mat. Fiz., {\bf 139} (1) (2004), pp. 29-44 (in
%Russian). English translation in: Theoret. and Math. Phys. {\bf 139} (1) (2004), pp. 473-485.

\bibitem[12]{Gou} Gould, M.D.: {\em `Characteristic identities for semi-simple Lie algebras'}. J.Austral. Math. Soc. B {\bf 26}, no.3,
(1985), pp. 257--283.

\bibitem[13]{Gr} Green, H.S.: {\em `Characteristic identities for generators of
${\rm GL}(n),\,{\rm O}(n)$ and ${\rm SP}(n)$'}. J. Math. Phys. {\bf 12} (1971), pp. 2106--2113.


%\bibitem[Gu]{Gu} Gurevich D. I., {\em 'Algebraic aspects of the quantum Yang-Baxter equation'}.
%Algebra i Analiz {\bf 2} (1990), no. 4, 119--148 (in Russian). English translation in:
%Leningrad Math. J. {\bf 2} (1991), no. 4, pp. 801--828.

\bibitem[14]{GPS1} Gurevich, D., Pyatov, P. and Saponov, P.: {\em `Hecke symmetries and characteristic
relations on reflection equation algebras'}. Lett. Math. Phys. {\bf 41} (1997), pp. 255--264.
% ArXiv: q-alg/9605048.

\bibitem[15]{GPS2} Gurevich, D., Pyatov, P. and Saponov, P.: {\em `Cayley-Hamilton theorem for quantum
matrix algebras of $GL(m|n)$ type'}. St. Petersburg Math. J. {\bf 17}, no.1, (2006), pp. 119-135.
%Algebra i Analiz, {\bf 17}, no. 1, (2005), pp. 160--181 (in Russian).
%ArXiv: math.QA/0412192.

\bibitem[16]{GPS3} Gurevich, D., Pyatov, P. and Saponov, P.: {\em `Quantum matrix algebras of the $GL(m|n)$ type:
The structure and spectral parameterization of the characteristic subalgebra'}.
Theor. Math. Phys. {\bf 147}, no.1, (2006), pp. 460–485.
%Preprint MPIM 2005--54. ArXiv: math.QA/0508506.

%\bibitem[GS]{GS} Gurevich D. and Saponov P., {\em `Geometry of non-commutative orbits related to Hecke
%symmetries'}. ArXiV: math.QA/0411579.

\bibitem[17]{GZB} Gould, M.D., Zhang, R.B. and Bracken A.J.: {\em `Generalized Gelfand invariants
and characteristic identities for quantum groups'}. J. Math. Phys. {\bf 32} no. 9, (1991),
pp. 2298--2303.

%\bibitem[Hay]{Hay}  Hayashi T., {\em `Quantum deformation of classical groups'}. Publ. Res. Inst.
%Math. Sci.  {\bf 28} (1992), no. 1, pp. 57--81.

\bibitem[18]{Hl} Hlavaty, L.: {\it `Quantized braided groups'}. J. Math. Phys. {\bf 35} (1994), pp. 2560--2569.
%ArXiv: hep-th/9210152.

\bibitem[19]{HSch} Heckenberger, I. and  Sch\"uler, A.: {\em `Symmetrizer and antisymmetrizer of the
Birman-Wenzl-Murakami algebras'}. Lett. Math. Phys. {\bf 50} (1999), pp. 45--51.
%ArXiv: math.QA/0002170.


\bibitem[20]{IWG} Isaac, P.S., Werry, J.L. and Gould, M.D.: {\em `Characteristic identities for Lie (super)algebras'}.
Journal of Physics: Conference Series {\bf 597} (2015) 012045.

\bibitem[21]{I}  Isaev, A.P.: {\em `Quantum groups and Yang-Baxter equations'}.
%MPIM Preprint 2004-132 (use http://www.mpim-bonn.mpg.de/html/preprints/preprints.html for uploads);
%previous version of this survey is published
Phys.\ Part.\ Nucl.\  {\bf 26} no.5, (1995), pp. 501--526.

\bibitem[22]{IMO}


\bibitem[23]{IOP} Isaev A., Ogievetsky O. and Pyatov P., {\em `Generalized Cayley-Hamilton-Newton
identities'}. Czech. Journ. of Physics {\bf 48} (1998), pp. 1369-1374. ArXiv: math.QA/9809047.

\bibitem[24]{IOP1} Isaev, A.P.,  Ogievetsky, O.V. and Pyatov, P.N.: {\em `On quantum matrix algebras
satisfying the Cayley-Hamilton-Newton identities'}. J. Phys. A: Math. Gen. {\bf 32} (1999), pp. L115--L121.
%ArXiv: math.QA/9809170.

\bibitem[25]{IOP2}  Isaev, A.P.,  Ogievetsky, O.V. and Pyatov, P.N.: {\em `Cayley-Hamilton-Newton
Identities and Quasitriangular Hopf Algebras'}. In Proc. of International Workshop `Supersymmetries and
Quantum Symmetries', July 27-31, 1999. Eds. E.Ivanov, S.Krivonos and A.Pashnev, JINR, Dubna
E2-2000-82, pp. 397--405. ArXiv: math.QA/9912197.

\bibitem[26]{IOP3} Isaev, A.P., Ogievetsky, O.V. and Pyatov P.N.: {\em `On R-matrix representations of
Birman-Murakami-Wenzl algebras'}.
%Proc. of Steklov Math. Inst. {\bf 246} (2004), pp. 147--153 (in Russian).
%English translation in:
Proc. Steklov Math. Inst. {\bf 246} (2004), pp. 134--141.
%ArXiv: math.QA/0509251

\bibitem[27]{IOPS} Isaev, A.P., Ogievetsky, O.V., Pyatov, P.N., and  Saponov P. A.: {\em `Characteristic
polynomials  for quantum matrices'}. In Proc. of International Conference in memory
of V.I.Ogievetsky `Supersymmetries and Quantum symmetries', (Dubna, Russia, 1997). Eds. J. Wess and
E. Ivanov, Lecture Notes in Physics, vol. {\bf 524}, pp. 322--330, Springer Verlag, 1998.


\bibitem[28]{It} Itoh M.: {\em `Capelli elements for the orthogonal Lie algebras'}. J. Lie Theory {\bf 10} (2000) pp. 463--489.

%\bibitem[JO]{JO} Jain V. and Ogievetsky O., {\em 'Classical isomorphisms for quantum groups'}. Modern Phys. Lett. A {\bf 7} (1992), no. 24, % pp. 2199--2209.

\bibitem[29]{JGr} Jarvis, P.D. and Green, H.S.: {\em `Casimir invariants and characteristic identities
for generators of the general linear, special linear and orthosymplectic graded Lie algebras'}.
J. Math. Phys. {\bf 20} no. 10,  (1979),  pp. 2115--2122.



\bibitem[30]{J} Jones, V.F.R.: {\em `On a certain value of the Kauffman polynomial'}. Comm. Math.
Phys. {\bf 125} (1989) pp. 459--467.




\bibitem[31]{KT1} Kantor, I. and Trishin, I.: {\em `On a concept of determinant in the supercase'}.
Communications in Algebra, {\bf 22} (1994), pp. 3679 -- 3739.

\bibitem[32]{KT2} Kantor, I. and Trishin, I.: {\em `On the Cayley-Hamilton equation in the supercase'}.
Communications in Algebra, {\bf 27} (1999), pp. 233 -- 259.

\bibitem[33]{KhO} Khoroshkin, S. and  Ogievetsky, O.:{\em `Diagonal reduction algebra and the reflection equation'}. 
Israel J.  Math. {\bf 221} no.2, (2017), pp.705--729.

\bibitem[34]{KS} Kulish and P.P., Sklyanin, E.K.: {\em `Algebraic structures related to reflection equations'}.
J. Phys. A  {\bf 25} no. 22, (1992),  pp. 5963--5975.
%ArXiv: hep-th/9209054.



\bibitem[35]{LR}  Leduc, R. and Ram, A.: {\em `A ribbon Hopf algebra approach to the irreducible
representations of centralizer algebras: the Brauer, Birman-Wenzl, and type A Iwahori-Hecke algebras'}.
Adv. Math.  {\bf 125} no. 1, (1997)  pp. 1--94.

\bibitem[36]{Mac} Macdonald, I.G.: {\em `Symmetric functions and Hall polynomials'}.
Oxford Mathematical Monographs, Oxford University Press, 1998.

\bibitem[37]{Man} Manin, Yu.I.: {\em `Notes on quantum groups and quantum de Rham complexes'}.
%Teor. Mat. Fiz. {\bf 92} (1992), no. 3, pp. 425--450 (in Russian). English translation in:
Theor. Math. Phys. {\bf 92} no. 3, (1992)  pp. 997--1023.

\bibitem[38]{Mol} Molev, A.: {\em `Sklyanin determinant, Laplace operators, and characteristic identities for classical Lie algebras'}.
J. Math. Phys. {\bf 36}, no.2, (1995), pp.923--943.

%\bibitem[Mol2]{Mol2} Molev A., {\em `Yangians and their applications'}.
% In `Handbook of Algebra', Vol.3, (M. Hazewinkel, Ed.), Elsevier, 2003.

\bibitem[39]{MRS} Molev, A.I., Ragoucy, E., Sorba, P.: {\em `Coideal subalgebras in quantum affine algebras'}.
Rev. Math. Phys. {\bf 15}, no.8, (2003), pp.789--822.

\bibitem[40]{Mudr} Mudrov, A. I.: {\em `Quantum conjugacy classes of simple matrix groups'}.
Comm. Math. Phys. {\bf 272}, no.3, (2007) pp. 635–-660.
%ArXiV: math. QA/0412538.

\bibitem[41]{M1} Murakami, J.: {\em `The Kauffman polynomial of links and representation theory'}.
Osaka J. Math. {\bf 24} (1987), pp. 745--758.

\bibitem[42]{NT} Nazarov, M. and Tarasov, V.: {\em `Yangians and Gelfand-Zetlin bases'}.
Publ. Res. Inst. Math. Sci. {\bf 30} no. 3, (1994),  pp. 459--478.
%ArXiV: hep-th/9302102.

\bibitem[43]{O} Ogievetsky, O.: {\em `Uses of quantum spaces'}. In Proc. of School
`Quantum symmetries in theoretical physics and mathematics' (Bariloche, 2000),  161--232,
Contemp. Math. {\bf 294} (2002), pp. 161--232.

\bibitem[44]{OP} Ogievetsky, O. and Pyatov, P.: {\em `Orthogonal and Symplectic Quantum Matrix Algebras and Cayley-Hamilton Theorem for them'}.
arXiv:math/0511618

\bibitem[45]{PS} Pyatov, P. and Saponov, P.: {\em `Characteristic relations for quantum matrices'}.
J. Phys. A: Math. Gen. {\bf 28} (1995), pp. 4415--4421.
%ArXiV: q-alg/9502012.

\bibitem[46]{Resh} Reshetikhin, N.Yu.: {\em Quasitriangular Hopf algebras and invariants of tangles}.
%Algebra i Analiz, {\bf 1}, no. 2, (1989), pp. 169--188 (in Russian). English translation in: 
Leningrad Math. J. {\bf 1} no. 2, (1990),  pp. 491--513.

\bibitem[47]{Resh2} Reshetikhin, N.Yu.: {\em `Multiparameter quantum groups and
twisted quasitriangular Hopf algebras'}. Lett. Math. Phys. {\bf 20} no. 4, (1990),  pp. 331--335.

\bibitem[48]{FRT} Reshetikhin, N. Yu., Takhtajan, L. A. and Faddeev, L. D.: {\em `Quantization of Lie groups
and Lie algebras'}.
%Algebra i Analiz, {\bf 1}, no. 1, (1989), pp. 178--206 (in Russian). English translation in:
Leningrad Math. J. {\bf 1} no. 1, (1990),  pp. 193--225.

\bibitem[49]{SchWZ} Schupp, P., Watts, P. and Zumino, B.: {\em `Bicovariant quantum algebras and quantum
Lie algebras'}. Comm. Math. Phys. {\bf 157} no. 2, (1993),  pp. 305--329.
%ArXiV: hep-th/9210150.

\bibitem[50]{TW} Tuba, I. and Wenzl, H.: {\em `On braided tensor categories of type BCD'}.
J. Reine Angew. Math.  {\bf 581}  (2005), pp. 31--69.
%ArXiv: math.QA/0301142.

\bibitem[51]{W} Wenzl, H.: {\em `Quantum groups and subfactors of type $B$, $C$, and $D$'}.
Comm. Math. Phys.  {\bf 133} no. 2, (1990),  pp. 383--432.

\bibitem[52]{Zh} Zhang, J.J.: {\em `The quantum Cayley-Hamilton theorem'}.
J. Pure Appl. Algebra {\bf 129} no. 1,  (1998),  pp. 101--109.

\end{thebibliography}
\end{document}